\documentclass[a4paper,10pt]{amsart}
\usepackage{epsf}  
\usepackage{tikz}
\usetikzlibrary{matrix,arrows,decorations.pathmorphing}
\usepackage{tikz-cd}
\usepackage{amsfonts, amssymb, amsthm, amsmath}  
\usepackage{hyperref}
\usepackage{enumerate}
\usepackage{pdfsync}
\usepackage{graphicx}
\usepackage{psfrag}
\newtheorem{thm}{Theorem}  

\newtheorem{cor}[thm]{Corollary}  
\newtheorem{lemma}[thm]{Lemma}  
\newtheorem{remark}[thm]{Remark}  
\newtheorem{defn}[thm]{Definition}  
\newtheorem{prop}[thm]{Proposition}

\numberwithin{thm}{section}  
\def\pf{\noindent\emph{Proof: }}  
\def\stop{\hfill$\square$}  
\providecommand{\Y}[1]{Y(#1)}

\providecommand{\tot}[1]{\ensuremath{[ #1]}}
\providecommand{\totl}[1]{\ensuremath{\lceil #1\rceil }}
\providecommand{\totb}[1]{\ensuremath{\underline{ #1}}}

\newcommand{\ex}{\bold}
\providecommand {\e}[1]{\mathfrak t^{#1}}

\providecommand{\C}[2]{\ensuremath {C^{#1,\underline{#2}}}}

\DeclareMathOperator{\id}{id}

\DeclareMathOperator{\coker}{coker}
\newcommand{\dbar}{\bar{\partial}}  

\newcommand{\dvert}{d_{\text vert}}
\providecommand{\ti}[1]{\tilde{\mathcal {#1}}}

\providecommand{\et}[2]{\ensuremath{\bold T^{#1}_{#2}}}
\providecommand{\lrb}[1]{\ensuremath{\left(#1\right)}}
\providecommand{\abs}[1]{\left\lvert #1\right\rvert}  
\providecommand{\norm}[1]{\left\lVert #1\right\rVert}

\author{Brett Parker }  
\email{brettdparker@gmail.com}
\title{Holomorphic curves in exploded manifolds: Regularity}

\begin{document}

\maketitle

\begin{abstract}The category of exploded manifolds is an extension of the category of smooth manifolds; for exploded manifolds,  some adiabatic limits appear as smooth families. This paper studies the $\dbar$ equation on variations of a given family of curves in an exploded manifold. 
Roughly, we prove that the $\dbar$ equation on variations of an exploded family of curves behaves as nicely as the $\dbar$ equation on variations of a smooth family of smooth curves, even though exploded families of curves allow the development of normal-crossing or log-smooth singularities.
The resulting regularity results are used in  \cite{uts, evc, vfc} to construct Gromov Witten invariants for exploded manifolds.

\end{abstract}

\tableofcontents

\section{Introduction}

\

Exploded manifolds, introduced in \cite{iec}, are a generalization of smooth manifolds akin to the generalization of complex manifolds given by log schemes, or the generalization of smooth manifolds given by Melrose's  b-structure  on manifolds with boundary and corners. 

A key feature of exploded manifolds is that normal-crossing or log-smooth degenerations  --- singular from the perspective of smooth manifolds  --- become smooth when viewed in the category of exploded manifolds.  Bubbling and node-formation on the domain of holomorphic curves seems singular from the perspective of smooth manifolds, but happens within smooth families of the exploded version of holomorphic curves. In this paper, we study the regularity of the $\dbar$ equation on such families of curves. The $\dbar$ equation is elliptic and well behaved when we restrict to the case of a single domain. For a smooth family of domains, the regularity of such an elliptic differential equation is also well understood.  Our goal is to understand the regularity of the $\dbar$ equation on apparently singular families of domains allowing bubbling and node-formation. 

Exploded manifolds provide a natural regularity to desire in families where bubbling and node-formation occurs; in this paper, we prove that the $\dbar$ equation has our desired regularity. In \cite{evc,vfc} this natural regularity  allows us to define Gromov--Witten invariants of suitable exploded manifolds, working on a moduli stack with a natural smooth structure.     Gromov--Witten invariants of suitable  exploded manifolds do not change in families, so using a degeneration that is smooth from the perspective of exploded manifolds,  Gromov--Witten invariants of a compact smooth symplectic manifold can be computed using an exploded manifold corresponding to a singular smooth manifold. The resulting computation is a gluing formula  involving  Gromov--Witten invariants relative normal crossing divisors. Such relative Gromov--Witten invariants are naturally defined within the category of exploded manifolds, and roughly equivalent to Ionel's GW invariants relative normal crossing divisors from \cite{IonelGW}, and log Gromov--Witten invariants defined by Gross and Siebert in \cite{GSlogGW} and Abramovich and Chen in \cite{Chen, acgw}. For comparisons of these different methods of defining Gromov--Witten invariants relative normal crossing divisors, see \cite{elc, tropicalIonel}.

One novel  feature of exploded manifolds is that each exploded manifold $\ex B$ is a set with several\footnote{Even though there are several topologies to keep in mind, we use the topology on $\ex B$ induced from  $\totl {\ex B}$ when referring to topological notions such as continuous functions or open subsets.} relevant topologies: a small scale, on which the exploded manifold generally looks like an infinite disjoint union of smooth manifolds, a large scale,  on which the exploded manifold looks like a union of integral affine polytopes, $\totb{\ex B}$, and a third topology, in which the exploded manifold looks like a stratified smooth space, $\totl{\ex B}$. On the large scale, holomorphic curves look like tropical curves --- piecewise linear maps of metric graphs, and the gluing formula referred to above is a sum of contributions corresponding to tropical curves.
This tropical gluing formula is proved in \cite{gfgw}; for a description of the formula in simple cases without using the language of exploded manifolds, see \cite{tpgf}, and for examples, see \cite{tec,3d}. 
  
  Despite the above strange features, differential geometry and topology works for exploded manifolds roughly like it does for smooth manifolds, so the reader may read the remainder of the introduction pretending that exploded manifolds are just manifolds. In the rest of the paper, some familiarity with the definitions in \cite{iec} or \cite{scgp} will be required.

\

\

This paper studies the $\dbar$ equation on a space of variations of a given family of curves in the category of exploded manifolds. In various guises, this is a key step in most constructions of Gromov--Witten invariants of general compact symplectic manifolds; \cite{Tian-Li, FO, McDuff, Ruanvirtual, Liu-Tian, CLW}. In section \ref{norms}, we define some Banach norms on this space of variations --- Banach norms involving many local choices. Then in section \ref{analysis}, the $\dbar$ equation is studied  on the corresponding Banach spaces of variations of our family of curves. Our final results are stated using the regularity $\C\infty1$, defined in section 7 of \cite{iec}.  (This regularity $\C\infty 1$ should be thought of as a generalization of smooth functions with an exponential decay condition of all weights $\delta<1$ at cylindrical ends.) With extra, less-analytic, work, the results stated in the remainder of this section are used in \cite{uts,evc,vfc} to describe the $\dbar$ equation on the moduli stack of $\C\infty1$ curves, and to construct the virtual fundamental class of the moduli space of holomorphic curves and Gromov Witten invariants of exploded manifolds.

\subsection{Acknowledgements}

Partially supported by SNF, No 200020-119437/1. Part of this work was also completed during the authorÕs stay at MIT,  UC Berkeley and the Mathematical Science Research Institute in Berkeley. Special thanks is due to an anonymous referee of a different paper, who went beyond the call of duty to provide valuable suggestions to improve this paper.

\subsection{Statement of results}

\

The technical heart of this paper  analyses   a family of (nonlinear) elliptic differential operators, $\dbar$, over a $\C\infty1$ family of curves\footnote{Curves are used in the sense of Definition 8.3 of \cite{iec}; in particular,  a curve is a complete $2$-dimensional exploded manifold $\ex C$ with a complex structure $j$, and a curve in an exploded manifold $\ex B$ is a map $f:\ex C\longrightarrow \ex B$. A family $\hat{\ex C}\longrightarrow \ex F$ is a kind of proper submersion, discussed in section 10 of \cite{iec}. See section \ref{lfs} for a brief discussion of the local structure of such a family of curves $\hat{\ex C}\longrightarrow \ex F$, including (co)tangent bundles and $\C\infty1$ regularity.  } $(\hat{\ex C}, j)\longrightarrow \ex F$. Let us first describe the form of the linearization, $D\dbar$, of $\dbar$. Both $\dbar$ and $D\dbar$ send sections of a complex vectorbundle $X$ to sections of $Y:=(T_{vert}^{*}\hat{\ex C}\otimes X)^{(0,1)}$.\footnote{For any family $\hat{\ex C}\longrightarrow \ex F$, let $T_{vert}\hat{\ex C}$ and $T_{vert}^{*}\hat{\ex C}$ respectively denote the vertical tangent and cotangent space. Our vector bundle $Y$ is the sub-bundle  of $T_{vert}^{*}\ex C\otimes_{\mathbb R}X$ fixed by $j\otimes J$, where $J$ is the complex structure on $X$.} The form of $D\dbar$ is
\begin{equation}\tag{$\star'$}\label{Ddbar form}D\dbar(\psi)=E'(\psi)+\frac 12(\nabla_{vert}\psi+J\circ\nabla_{vert}\psi\circ j)\end{equation}
where
\begin{itemize}\item $\nabla_{vert}$ is the restriction of a $\C\infty1$ connection\footnote{A $\C\infty1$ connection on a vectorbundle $V$ over an exploded manifold $\ex B$ is a linear map $\nabla$ from $\C\infty1$ sections of $V$ to $\C\infty1$ sections of $V\otimes T^{*}\ex B$ that satisfies the usual derivation condition, and so that $\nabla _{v}w=0$ for any $\mathbb R$-nil vector $v$ and $\C\infty1$ section $w$ of $V$. A $\mathbb R$-nil vector $v$ is one for which $vf=0$ for all smooth or $\C\infty1$ functions $f$. In particular, on a curve, the nonzero $\mathbb R$-nil vectors are the nonzero vectors on edges. See section 6 of \cite{iec} for a discussion of $T\ex B$.} $\nabla$ on $X$ to the vertical tangent space of $\hat{\ex C}$,  $T_{vert}\hat{\ex C}\subset T\hat{\ex C}$;
\item $J$ indicates the complex structure on $X$;

\item and \[E': X\longrightarrow Y\]
is a $\C\infty1$ map of  vectorbundles, vanishing on the edges\footnote{The tropical part of a curve $\ex C$ is a graph $\totb{\ex C}$.  As in Definition 8.3 of \cite{iec}, an edge of a curve refers to a stratum of  $\ex C$ over one of the edges of this graph $\totb{\ex C}$. The smooth part of a $\ex C$ is a nodal curve $\totl{\ex C}$ with marked points. Edges of $\ex C$ correspond to the nodes and marked points of this nodal curve $\totl{\ex C}$.} of curves in $\hat{\ex C}$.
\end{itemize}

Our linearized operator determines  a linear map 
\[D\dbar:X^{\infty,\underline 1}\longrightarrow Y^{\infty,\underline 1}\]
where   $X^{\infty,\underline 1}$ is the sheaf of $\C\infty1$ sections of $X$, and $Y^{\infty,\underline 1}$ is the sheaf of $\C\infty1$ sections of $Y$  vanishing on edges of curves in $\hat {\ex C}\longrightarrow \ex F$. The map $D\dbar$ above is then a map of sheaves of $\C\infty1(\ex F,\mathbb R)$--modules over $\ex F$. We can restrict $D\dbar$, $X^{\infty,\underline 1}$ and $Y^{\infty,\underline 1}$ to any open subset of $X$, but we can also restrict them to any individual curve $f\in \ex F$ to determine a map 
\[D\dbar(f):X^{\infty,\underline 1}(f)\longrightarrow Y^{\infty,\underline 1}(f)\]
which is an elliptic differential operator, satisfying the usual Fredholm and regularity properties. 

Although $X^{\infty,\underline 1}$ and $Y^{\infty,\underline 1}$ are not obviously vectorbundles over $\ex F$, because they are not obviously locally trivial,\footnote{This is related to the possibly-insurmountable  difficulty of putting a Banach orbifold structure on the moduli space of curves of a given regularity close to a nodal curve, and the need for different analytic structures such as polyfolds; see \cite{polyfold0, polyfold1, polyfold2, polyfoldint, polyfoldsc, polyfoldgw}.  We use standard analysis on Banach space completions of $X^{\infty,\underline 1}$ and $Y^{\infty,\underline 1}$ involving  unnatural choices, but obtain results independent of such choices. 
} they have a natural notion of a finite-rank sub-vectorbundle: 

\begin{defn}\label{sub-vectorbundle} A rank--$k$  sub-vectorbundle $K$ of $X^{\infty,\underline 1}$ (or $Y^{\infty,\underline 1}$) is a locally free, rank--$k$ subsheaf whose restriction to  $X^{\infty,\underline 1}(f)$ (or $Y^{\infty,\underline 1}(f)$) is a $k$-dimensional linear subspace $K(f)$ for all $f\in \ex F$.
\end{defn}

\begin{thm}\label{linear theorem} Given any curve $f\in \ex F$,  
$D\dbar(f)$
has closed image,  finite-dimensional kernel and cokernel, and index
 \[\dim \ker D\dbar(f)-\dim\coker D\dbar (f) =2c_{1}-2n(g-1)\] where $c_{1}$ is the integral of the first Chern class of $X$ over the curve $f$, $n$ is the complex-rank of $X$ and $g$ is the genus of the curve $f$.

 Given any finite-rank sub-vectorbundle $V\subset Y^{\infty,\underline 1}$ so that  $D\dbar(f)$ is transverse to $V(f)\subset Y^{\infty,\underline 1}(f)$, there  exists a neighborhood $\ex F'$ of $f\in \ex F$, on which $D\dbar$ surjects onto $Y^{\infty,\underline 1}(\ex F')/V$ with kernel  a finite-rank  sub-vectorbundle  $K=D\dbar^{-1}(V)$ of $X^{\infty,\underline 1}(\ex F')$. 

\end{thm}

The first part of
Theorem \ref{linear theorem} simply tells us that, restricted to any individual curve,  $D\dbar$ is Fredholm and has the expected index. The second part is a `gluing and regularity' result. It tells us that we can study $D\dbar$ with a finite-dimensional, $\C\infty1$ approximation $K\longrightarrow V$, where $K$ and $V$ are finite-rank $\C\infty1$ vectorbundles. A proof of Theorem \ref{linear theorem} is at the end of section \ref{linear section}.  The analogous theorem for the nonlinear operator $\dbar$  is that $\dbar^{-1}(V)$ is a $\C\infty1$ manifold of the expected dimension.


%
%
%

The nonlinear operators we consider are in the form
\begin{equation}\tag{$\star$}\label{dbar form}\begin{split}\dbar (\nu)&=E(\nu)+\frac 12\lrb{H'(\nu)\circ \nabla_{vert}\nu+J\circ H'(\nu)\circ \nabla_{vert}\nu\circ j}
\\ &:= E(\nu)+H(\nu)(\nabla_{vert}\nu)\end{split}\end{equation}
where
\begin{itemize}
\item $E$ is a (nonlinear) $\C\infty 1$ map so that the following diagram commutes
\[\begin{tikzcd}X\ar{dr}\rar{E}&Y\dar
\\ & \hat{\ex C}\end{tikzcd}\]
and so that, restricted to edges  of curves in $\hat {\ex C}$,  $E$ has image the zero-section of $Y$;
\item and $H'$ is a (nonlinear) $\C\infty1$ map to the space of invertible, $\mathbb R$-linear endomorphisms of $X$ so that the following diagram commutes
\[\begin{tikzcd}X\ar{dr}\rar{H'}& X\otimes_{\mathbb R} X^{*}\dar
\\ & \hat{\ex C}\end{tikzcd}\]
and so that $H'$ restricted to the zero-section is the identity.
\end{itemize}

Theorem \ref{smooth dbar} on page \pageref{smooth dbar} implies that $\dbar$ is continuously differentiable, at least when  $\ex F$ is extendable. The derivative of $\dbar$ at the zero-section is a family of linear elliptic differential operators in the form of (\ref{Ddbar form}). Moreover,  the derivative of $\dbar$ at any other section is in the form of (\ref{Ddbar form}) with a different $J$, so Theorem \ref{linear theorem} still applies.

%
 
 Our nonlinear regularity theorem becomes easier to state if $D\dbar$ is injective; we can achieve this by restricting the domain of $\dbar$. In particular, choose some  finite set of non-intersecting sections $s_{1},\dotsc,s_{m}$ of $\hat{\ex C}\longrightarrow \ex F$ that avoid  the edges of curves in $\hat {\ex C}$.  Now define $X^{\infty,\underline 1}$ to be the sheaf of $\C\infty1$ sections of $X$ that vanish at the image of $s_{i}$. Theorem \ref{linear theorem} still holds for $D\dbar$ with such a restricted domain, but the index of $D\dbar$ is now $2c_{1}-2n(g-1+m)$. By choosing enough such sections, we may assume that $D\dbar(f)$ is injective, and apply the following theorem.

\begin{thm}\label{nonlinear theorem}
Given a curve $f$ in $\ex F$
  and a finite-rank  sub-vectorbundle $V\subset Y^{\infty,\underline 1}$  so that $\dbar f \in V(f)$, and   $D\dbar(f)$ is injective and complementary to $V(f)\subset Y^{\infty,\underline 1}(f)$, there exists  a neighborhood $\ex F'$ of $f\in \ex F$ and a solution $\nu\in X^{\infty,\underline 1}(\ex F')$ to the equation 
\[\dbar\nu=0\mod V\ .\] 
 Moreover, there exists a neighborhood $O$ of $0\in X^{\infty,\underline 1}(\ex F')$ so that $\nu$ is the unique solution to $\dbar\nu=0\mod V$ within $O$ in the following sense. Given any curve $f'\in \ex F'$, let  $\nu(f')$ and $O(f')$ be the restriction of  $\nu$ and $O$   to $X^{\infty,\underline 1}(f')$. Then $\nu(f')$ is the unique solution to the equation $\dbar \nu(f')=0\mod V(f')$ within $O(f')$.
\end{thm}

So, when we parametrize curves by a fixed family of domains, the moduli space of solutions to $\dbar=0\mod V$ is regular when $D\dbar$ is transverse to $V$. Note that our fixed family can include bubbling and node-formation behavior, so this regularity theorem includes gluing analysis. When we no longer have a fixed family of domains, this theorem is key for proving  the regularity of the moduli stack of solutions to $\dbar=0\mod V$ --- for example, where $D\dbar$ is transverse to zero, the moduli stack of holomorphic curves is locally a $\C\infty1$ orbifold. For details and precise statements in the setting of the moduli stack of curves, see \cite{evc}.
Theorem \ref{nonlinear theorem} follows immediately from Theorem \ref{family regularity}, stated on page \pageref{family regularity}.

\subsection{Local model for families of curves}\label{lfs}

\

This section illustrates some  key differences between exploded manifolds and smooth manifolds; the reader already familiar with exploded manifolds may skip it.  Consider the standard local model for node formation.
\[\mathbb C^{2}\xrightarrow{z_{1}z_{2}} \mathbb C\]
Away from $0$, this map is a submersion with fiber a cylinder, but over $0$ the fiber is a pair of planes, joined at $0$. The corresponding local model for node formation in the category of exploded manifolds is 
\[(\ex T^{1}_{[0,\infty)})^{2}=\et 2{[0,\infty)^{2}}\xrightarrow {\tilde z_{1}\tilde z_{2}}\ex T^{1}_{[0,\infty)}\ .\]
Every exploded manifold $\ex B$ comes with a natural map to a stratified smooth space called its smooth part $\totl{\ex B}$. For example, an exploded curve $\ex C$ has smooth part $\totl {\ex C}$ a nodal Riemann surface, with strata given by nodes and marked points.\footnote{The curve $\ex C$ also has a tropical part $\totb{\ex C}$ consisting of a complete metric graph with a vertex for every component of the nodal curve $\totl{\ex C}$, an edge for every node of $\totl{\ex C}$, and a semi-infinite edge for every marked point of $\totl{\ex C}$. So, nodes and marked points of $\totl{\ex C}$ correspond to internal and external edges of $\totb{\ex C}$, and we will often refer to such strata of $\ex C$ as `edges'. } The smooth part of the above  exploded  node-formation model is the standard node-formation model (stratified by where coordinate functions are 0). In particular $\totl{\et 1{[0,\infty)}}=\mathbb C$, and the smooth part $\totl{\tilde z}$ of the standard coordinate on $\et 1{[0,\infty)}$ is the standard coordinate on $\mathbb C$. In fact, $\mathbb R$-valued smooth or continuous functions on an exploded manifold $\ex B$ factor through its smooth part $\ex B\longrightarrow \totl{\ex B}$, so the smooth or continuous functions on our exploded node-formation model are the same as those on the standard node-formation model. Moreover,  although there are several topologies on $\ex B$, by default, we use the topology induced by the map $\ex B\longrightarrow \totl{\ex B}$.

The tangent sheaf, however is different. The tangent space of $\et 1{[0,\infty)}$ has a standard trivialization as $\mathbb R^{2}\times \et 1{[0,\infty)}$, with basis  the real and imaginary parts of  $\tilde z\frac\partial{\partial{\tilde z}}$; similarly, a basis for the cotangent space is the real and imaginary parts of $\tilde z^{-1}d\tilde z$.  The derivative of our exploded node-formation model is constant in this standard basis, sending $a\tilde z_{1}\frac\partial{\partial\tilde z_{1}}+b\tilde z_{1}\frac\partial{\partial \tilde z_{2}}$ to $(a+b)\tilde z\frac\partial{\partial \tilde z}$. So, unlike the standard node-formation model, our exploded node-formation model is a submersion. 

A second difference is what replaces $0\in\mathbb C$ and the singular fiber of the standard node-formation model. On $\et1{[0,\infty)}$, our coordinate $\tilde z$ takes values in $\mathbb C^{*}\e {[0,\infty)}\subset \mathbb C^{*}\e {\mathbb R}$. As a group, $\mathbb C^{*}\e{\mathbb R}$ is $(\mathbb C^{*},\times)\times(\mathbb R,+)$ --- this group action preserves our standard basis for the tangent space. The smooth part map sends points $c\e{0}$ to $c$, and everything else  to $0$. In some sense,  $0\in \mathbb C$ is replaced by an infinite set of cylinders, $\mathbb C^{*}$, one for every point in $(0,\infty)$.  Accordingly, the singular fiber of the standard node-formation model is replaced by many exploded manifold fibers, one for each $c\e{a}\in\mathbb C^{*}\e{(0,\infty)}$. One way to think of this $c\e{a}$ is as an infinitesimal gluing parameter; the fiber over $c\e{a}$ has $2$ natural coordinates $\tilde z_{1}$ and $\tilde z_{2}$, related by $\tilde z_{2}=c\e{a}\tilde z_{1}^{-1}$.

In this simple case, let us consider the natural $\C\infty1$ regularity used in this paper. We can think of our exploded manifold $\et 1{[0,\infty)}$ replacing $\mathbb  C\setminus 0$ with a manifold with a cylindrical end, and replacing $0\in \mathbb  C$ with lots of `cylinders at infinity'. A standard  setup for analysis on manifolds with cylindrical ends uses Sobolov spaces with exponential weights at ends, however, as cylindrical coordinates correspond to $\log z$, an exponential weight $\delta$ corresponds to using a weight $\abs z^{-\delta}$. The regularity we will achieve in this paper is $\C\infty1$. Let us describe  $\C\infty1$ functions on $\et 1{[0,\infty)}$ using the coordinate $z=\totl{\tilde z}$ on the smooth part of $\et 1{[0,\infty)}$. A continuous $\mathbb R$-valued function $h(z)$ is $\C\infty1$ if $(h(z)-h(0))\abs z^{-\delta}$ extends to a continuous function for all $\delta<1$, and the same applies to any number of derivatives of $h$, using the real or imaginary part of $z\frac \partial{\partial z}$. How about $\C\infty1$ functions on $\et 2{[0,\infty)^{2}}$? Such functions $h(z_{1},z_{2})$ have 
\[(h(z_{1},z_{2})-h(0,0))(\abs z_{1}+\abs z_{2})^{-\delta},\]
\[(h(z_{1},z_{2})-h(z_{1},0))\abs {z_{2}}^{-\delta}\ ,\] 
\[(h(z_{1},z_{2})-h(0,z_{2}))\abs {z_{1}}^{-\delta}\ ,\] and
\[(h(z_{1},z_{2})-h(z_{1},0)-h(0,z_{2})+h(0,0))\abs{z_{1}z_{2}}^{-\delta}\] 
extending to continuous functions for all $\delta<1$, and the same holds for all derivatives of $h$ using our standard basis consisting of the real and imaginary parts of $z_{i}\frac \partial{\partial z_{i}}$. The different conditions  above correspond to different sets of strata of $[0,\infty)^{2}$. In this case, the first expression above is redundant, but the fourth expression might seem excessive. It ensures that each fiberwise-constant $\C\infty1$ function on $\et 2{[0,\infty)^{2}}$ is the pullback of a $\C\infty1$ function on $\et 1{[0,\infty)}$. This is true more generally: given a family of exploded manifolds, $\C\infty1$ functions on the base correspond to fiberwise-constant $\C\infty1$ functions on the total space; moreover,  any $\C\infty1$ function on a fiber extends to a $\C\infty1$ function on the total space. 

Most `singular' behavior in our family of curves $\hat{\ex C}\longrightarrow \ex F$ is captured in the above exploded node-formation model. Locally on $\hat{\ex C}$, $\hat{\ex C}\longrightarrow \ex F$  is either locally trivial, or locally modeled on a base-change of our node-formation model, given by a fiber product-diagram.
\begin{equation}\label{node local model}\begin{tikzcd}\mathbb R^{k}\times\et {m+1}Q\dar\rar & \et 2{[0,\infty)^{2}}\dar
\\ \mathbb R^{k}\times \et mP\rar & \et 1{[0,\infty)}\end{tikzcd}\end{equation}
When $P=[0,\infty)^{n}$, $\et mP$ is just $(\et 1{[0,\infty)})^{m}$; if $P$ is a more complicated polytope, cut from $[0,\infty)^{M}$ by integral-affine equations $\sum\alpha_{i} x_{i}=c$, the exploded manifold $\et mP$ is the subset of $(\et 1{[0,\infty)})^{M}$ cut out by the corresponding monomial equations $\prod\tilde z_{i}^{\alpha_{i}}=\e{c}$, and the $\C\infty1$ functions on $\et mP$ consist of the restriction of $\C\infty1$ functions from $(\et 1{[0,\infty)})^{M}$. 

\subsection{Geometric setup}

\

This section explains how theorems \ref{linear theorem} and \ref{nonlinear theorem} apply to the $\dbar$ equation on the space of variations of a given family of curves $\hat f$  in a smooth family of exploded manifolds $\hat{\ex B}\longrightarrow \ex B_{0}$.  To discuss holomorphic curves in such a family $\hat{\ex B}$, we require an almost complex structure $J$ on each fiber.\footnote{Almost complex structures are discussed in section 8 of \cite{iec}. For this paper, we shall not require that $J$ is civilized in the sense of Definition 8.2 of \cite{cem}, and also not require the stronger condition of being $\dbar$-log compatible, Definition 3.1 of \cite{cem}. This stronger condition is used in \cite{cem} to prove compactness for the moduli stack of holomorphic curves, and if it is assumed,  the results of the present paper could conceivably be improved, with `smooth' replacing $\C\infty1$. Such an improvement would not be particularly important, because the category of $\C\infty 1$ exploded manifolds behaves similarly to the category of smooth exploded manifolds. }

 A family of curves $\hat f$ in a family of exploded manifolds $\hat{\ex B}\longrightarrow \ex B_{0}$ is a commutative diagram,
\begin{equation}\label{curves}\begin{tikzcd}(\hat{\ex C},j)\dar{\pi_{\ex F}}\rar{\hat f} &(\hat{\ex B},J)\dar
\\ \ex F\rar&\ex B_{0}\end{tikzcd}\end{equation}
where $\pi_{\ex F}$ is a family of exploded manifolds, and $j$ is a family of almost complex structures so that each fiber is a curve; see Definition 8.3 of \cite{iec}. Such a family has regularity $\C\infty1$ if $j$, and all maps, are $\C\infty1$. By variations of $\hat f$, we mean variations of the map $\hat f$ so that all other maps in the above commutative diagram are fixed.  Locally, we can identify an open neighborhood of $\hat f$ in the space of such variations with $\C\infty1$ sections of the vectorbundle $\hat f^{*}T_{vert}\hat{\ex B}$. The complex vectorbundle $X$ in the previous section is $\hat f^{*}T_{vert}\hat{\ex B}$ with its almost complex structure, $J$.

Given such a family of curves, use 
\[\dvert \hat f: T_{vert}\hat {\ex C}\longrightarrow T_{vert}\hat{\ex B}\]
to indicate $d\hat f$ restricted to the vertical tangent space, $T_{vert}\hat {\ex C}\subset T\hat{\ex C}$. This  family $\hat f$ is a family of  holomorphic curves if $\dvert\hat f\circ j=J\circ\dvert\hat f$.
\  Define a map
\[\dbar \hat f: T_{vert}\hat {\ex C}\longrightarrow T_{vert}\hat{\ex B}\]
\[\text{by }\ \ \ \dbar \hat f:=\frac 12\lrb{ \dvert f +J\circ \dvert f\circ j}\ .\]
We shall also consider $\dbar\hat f$ as a section of the vectorbundle,
 \begin{equation}\label{Ydef}\Y{\hat f}:=\lrb{{T_{vert}^{*}\hat{\ex C}}\otimes {\hat f^{*}T_{vert}\hat{\ex B}}}^{(0,1)}\end{equation}
 which is the sub-bundle of $T_{vert}^{*}\hat{\ex C}\otimes_{\mathbb R} {\hat f^{*}T_{vert}\hat{\ex B}}$ fixed by $j\otimes J$.  In fact, $Y(\hat f)$ is the pullback
 of a  vectorbundle over the universal curve over the moduli stack of $\C\infty1$ curves in $\hat{\ex B}\longrightarrow \ex B_{0}$\footnote{See section 11 of \cite{iec} for a short discussion of the moduli stack of $\C\infty1$ curves. Roughly speaking, this moduli stack is a category with objects consisting of $\C\infty1$ families of curves. A more detailed analysis of this moduli stack is contained in \cite{evc}.}; given any map  of families of curves $\hat f\longrightarrow \hat g$, there is a corresponding map of vectorbundles $\Y{\hat f}\longrightarrow \Y{\hat g}$. 

\

When we locally identify variations of $\hat f$ with sections of $\hat f^{*}T_{vert}\hat{\ex B}$, we can also choose a  trivialization of  $Y$ over this space of variations. Then the $\dbar$ equation takes the form (\ref{dbar form}), and we can apply our regularity theorems.

\begin{defn}\label{trivialization def} Given a $\C\infty1$ family of curves, 
\[\begin{tikzcd}\hat {\ex C}\rar{\hat f}\dar& \hat{\ex B}\dar
\\ \ex F\rar &\ex B_0\end{tikzcd}\]
 a  trivialization is:
 \begin{enumerate}
 \item
 a $\C\infty1$ map, $\mathcal F$, so that the following diagram commutes
 \[\begin{tikzcd} \hat f^{*}{T_{vert}\hat{\ex B}}\rar{\mathcal F}\dar &\hat{\ex B} \dar
\\ \ex F\rar &\ex B_0
 \end{tikzcd}\]
 and so that \begin{enumerate}\item $\mathcal F$, restricted to the zero-section,  equals  $\hat f$, 
 \item $T\mathcal F$ restricted to the natural inclusion  $\hat f^{*}{T_{vert}\hat{\ex B}}\subset T(\hat f^{*}T_{vert}\hat{\ex B})$ over the zero-section is  the identity,
 \item  $T\mathcal F$ restricted to the vertical tangent space at any point of $\hat f^{*}T_{vert}\hat{\ex B}$ is injective; 
 \end{enumerate} 
\item and a $\C\infty1$ vectorbundle map, $\Phi$, that fits into the commutative diagram below:
 \[\begin{tikzcd}\mathcal F^{*}{T_{vert}\hat{\ex B}}\rar{\Phi}\dar& \hat f^{*}{T_{vert}\hat{\ex B}}\dar{\pi}
 \\ \hat f^{*}{T_{vert}\hat{\ex B}}\rar{\pi} &\hat{\ex C}\end{tikzcd}\]
and which is a $J$-preserving isomorphism on each fiber and is the identity when the vectorbundle  $\mathcal F^{*}{T_{vert}\hat{\ex B}}\longrightarrow \hat f^{*}{T_{vert}\hat{\ex B}}$ is restricted to the zero-section of $\hat f^{*}T_{vert}\hat{\ex B}$.
\end{enumerate}
 \end{defn}

For example, we may construct a trivialization by extending $\hat f$ to a map $\mathcal F$ satisfying the above conditions (for instance by choosing a smooth connection on $T_{vert}\hat {\ex B}$ and reparametrising the exponential map on a neighborhood of the zero-section in  $f^{*}T_{vert}\hat {\ex B}$), and letting $\Phi$ be given by parallel transport along a linear path to the zero-section using a smooth $J$-preserving connection on $T_{vert}\hat{\ex B}$.

\begin{defn}
 A trivialization allows us to define $\dbar$ of a section 
 \[\nu :\hat{\ex C}\longrightarrow\hat f^{*}{T_{vert}\hat{\ex B}}\] as follows:\\ \\$\mathcal F\circ \nu$ is a family of curves 
 \[\mathcal F\circ\nu:\hat {\ex C}\longrightarrow \hat {\ex B}\] so $\dbar (\mathcal F\circ \nu)$ is a section of $\Y{\mathcal F\circ\nu}=\lrb{{T_{vert}^{*}\hat{\ex C} }\otimes {(\mathcal F\circ \nu)^{*}T_{vert}\hat{\ex B}}}^{(0,1)}$. 
 Applying the map $\Phi$ to the second component of $\dbar(\mathcal F\circ \nu)$ defines a section $\dbar \nu$  of $\Y{\hat f}$.

\end{defn}

Lemma \ref{local form} below implies that $\dbar$ as defined above is in the form (\ref{dbar form}).

\

 Given a trivialization for a family $\hat f$, we can define the following `simple perturbation' of the $\dbar$ equation on variations of $\hat f$.

\begin{defn}\label{simple perturbation}
Given a trivialization for  $\hat f$, a simple perturbation of $\dbar$ is a map $\dbar'$ from sections  $\hat{\ex C}\longrightarrow\hat f^{*}{T_{vert}\hat{\ex B}}$ to sections of $\Y{\hat f}$ so that 
\[\dbar'\nu=\dbar(\nu)+\Psi(\nu)\]
where $\Psi$ is a (possibly nonlinear)  $\C\infty1$ map so that the following diagram commutes 
\[\begin{tikzcd}\hat f^{*}{T_{vert}\hat {\ex B}}\rar{\Psi}\ar{dr}& \Y{\hat f}\dar
\\ &\hat{\ex C} \end{tikzcd}\]
and so that the image of $\Psi$ restricted to edges of curves in $\hat {\ex C}$ is contained in the zero-section of $Y(\hat f)$.

\
%

\end{defn}

\begin{lemma}\label{local form} 
Let $\hat f$ be a $\C\infty1$ family of curves with a choice of trivialization. Then any simple perturbation $\dbar'$ of $\dbar$ is in the form (\ref{dbar form}). In particular,  in any coordinate chart $U$ on $\hat{\ex C}$ where $\hat f^{*}T_{vert}\hat{\ex B}$ is trivialized as $\mathbb R^{n}\times U$, consider a section as a map $\nu:U\longrightarrow \mathbb R^{n}$. In such a coordinate chart, there is the following formula for $\dbar'$.
\[\dbar' \nu (u)= E(\nu(u),u)+H(\nu(u),u)(\dvert\nu)\]
where
\[H(\dvert\nu):=\frac 12(H'\circ \dvert\nu +J\circ H'\circ \dvert\nu\circ j) \]
and $E$, $H$ and $H'$ are as described in (\ref{dbar form}). 
\end{lemma}

\

\pf

The map $E(x,u)$ is  $\dbar'$ of the constant section $u\mapsto(x,u)$.  This is a $\C\infty1$ map, because $\hat f$, $J$, and $j$ are $\C\infty1$, and it vanishes on the edges of fibers of $\hat {\ex C}\longrightarrow \ex F$ because every  $\C\infty1$ curve is automatically holomorphic restricted to any edge, and Definition \ref{simple perturbation} tells us that simple perturbations do not modify $\dbar$ on edges of curves. 

The tensor $H'$ is given by the formula
\[H'(x,u) :=\Phi(x,u)\circ d_{vert}\mathcal F(x,u)\]
  where $\Phi$ and $\mathcal F$ are as in Definition  \ref{trivialization def}. To interpret $H'(x,u)$ as an endomorphism of $\hat f^{*}T_{vert}\hat{\ex B}$ at $x$, identify the vertical tangent space of $\hat f^{*}T_{vert}\hat{\ex B}\longrightarrow \hat{\ex C}$ at $(x,u)$ with $\hat f^{*}T_{vert}\hat{\ex B}$ at $u$. Both $\Phi(x,u)$ and  the vertical derivative of $\mathcal F$  are injective, and both are equal to the identity at the zero-section, so  the same holds for $H'$. Because $\Phi$ and $\mathcal F$ are $\C\infty1$, $H'$ is $\C\infty1$ too. Direct computation gives that  $\dbar' \nu$ obeys the above formula.

Given any $\C\infty1$ connection, $\nabla$ on $\hat f^{*}T_{vert}\hat{\ex B}$, $\nabla_{vert}-d_{vert}$ is a $\C\infty1$ section of $T^{*}_{vert}\hat{\ex C}\otimes\hat f^{*} T_{vert}\hat {\ex B}\otimes \hat{f}^{*}T_{vert}^{*}\hat{\ex B}$ which vanishes on edges of curves in $\hat{\ex C}\longrightarrow \ex F$ because $\nabla_{v}w$ always vanishes when $v$ is a $\mathbb R$-nil vector, and the vectors on edges of a curve are $\mathbb R$-nil.

It follows that $E'(\nu,u):=H(\nu,u)(\nabla_{vert}\nu-d_{vert}\nu)$ is a $\C\infty1$ map satisfying the conditions on $E$, so exchanging $d_{vert}$ for $\nabla_{vert}$ gives  (\ref{dbar form}). It also follows that an operator is in the form above on coordinate charts if and only if it is in the form (\ref{dbar form}).

\stop

\section{Some norms on sections of vectorbundles}

In what follows, we define some norms on the space of sections of a real vectorbundle $ V\longrightarrow\hat{\ex C}$ over a family of curves $\hat{\ex C}\xrightarrow{\pi_{\ex F}}\ex F$. Applied to the vectorbundles $X$ and $Y$, these norms will determine a series of Banach-Space completions $X_{k,\delta}$ and $Y_{k,\delta}$ of (bounded subspaces) in $X^{\infty,\underline 1}$ and $Y^{\infty,\underline 1}$ so that $\dbar$ determines a well-behaved, continuously differentiable map $X_{k,\delta}\longrightarrow Y_{k,\delta}$, and so that convergence in $\C\infty1$ is implied by convergence in $\norm{\cdot}_{k,\delta}$  for all natural numbers $k$ and $\delta<1$. We will prove our desired $\C\infty1$ regularity by proving regularity in each of these norms. 

\subsection{$e_{S}$ and $\Delta_{S}$}

\

For defining and working with our norms, we need to recall some basic notions from \cite{iec}. In particular, we need the notation  for a  coordinate chart on an exploded manifold from section 3 of \cite{iec}, the notion of a stratum of $P$ or $\et mP$ from section 4, and the following definitions from section 7.

\begin{defn}[The operator $e_{S}$]
Given any $C^0$ function  $h$ on $\mathbb R^n\times \et mP$ and a stratum $S\subset P$, define 
\[e_S(h)(x,\tilde z_{1},\dotsc,\tilde z_{m}):=h(x,\tilde z_{1}\e {(a_{1}-\totb{\tilde z_{1}})/2},\dotsc,
\tilde z_{m}\e {(a_{m}-\totb{\tilde z_{m}})/2})\]
where $(a_1,\dotsc,a_m)$ is any point  in $S$, and $\totb{c\e x}:= x$. 
%

\end{defn}

So $e_{S}h(x,\tilde z)$ samples the function $h$ at a point with tropical part half way between $\totb{\tilde z}$ and the point $ a\in S$. Note that  $e_{S}h$ does not depend on the choice of the point $a\in S$, because a different choice of $a\in S$ would ask us to sample $h$ at a topologically indistiguishable point.

\

For example consider $\et 22:=\et 2{[0,\infty)^{2}}$. The polytope $[0,\infty)^{2}$ has a zero-dimensional stratum, $S_{0}$,  two one-dimensional strata 
\[S_{1}:=(0,\infty)\times 0 \ \ \ \ \ S_{2}:= 0\times(0,\infty)\] and one two-dimensional stratum $S_{3}:=(0,\infty)^{2}$. 

 If we have a function $h\in C^0(\et 22)$, then $e_{S_{0}}h=h$, but
\[ e_{S_{1}}h(z_{1},z_{2})=h(0,z_{2})\ \ \ \ \ \ \ e_{S_{2}}h(z_{1},z_{2})=h(z_{1},0) \ \ \ \ \ e_{S_{3}}h(z_{1},z_{2})=h(0,0)\ .\]

As a second example, consider $\et1{[0,1]}$. Smooth or continuous functions on $\et 1{[0,1]}$ are generated by $\zeta_{1}=\totl{\tilde z}$ and $\zeta_{2}=\totl{\e{1}\tilde z^{-1}}$. There are three strata of  $[0,1]$ to consider: $0$, $1$, and $(0,1)$.
\[e_{0}\zeta_{1}=\zeta_{1}\ \ \ \  e_{0}\zeta_{2}=0\ \ \ \ e_{0}h(\zeta_{1},\zeta_{2})=h(\zeta_{1},0)\]
\[e_{1}\zeta_{1}=0\ \ \ \ e_{1}\zeta_{2}=\zeta_{2}\ \ \ \ e_{1}h(\zeta_{1},\zeta_{2})=h(0,\zeta_{2})\]
\[e_{(0,1)}\zeta_{1}=0\ \ \ \ e_{(0,1)}\zeta_{2}=0\ \ \ \ e_{(0,1)}h(\zeta_{1},\zeta_{2})=h(0,0)\]
We can consider $\et 1{[0,1]}$ as the subset of $\et 22$ where $\tilde z_{1}\tilde z_{2}=\e 1$. From this perspective we can relate the above two examples by  $e_{0}=e_{S_{1}}$ and $e_{1}=e_{S_{2}}$.

\

 In general, the smooth or continuous functions on $\et mP$ are generated by the smooth monomials: functions $\zeta_{i}$ of the form $\totl{\e {a_{i}}\tilde z^{\alpha^{i}}}:=\totl{\tilde \zeta_{i}}$. For any stratum $S\subset P$ one of the following two options hold:
 \begin{itemize}
 \item $e_{S}\zeta_{i}=0$, $\zeta_{i}$ vanishes on the stratum of $\et mP$ corresponding to $S$ and $\totb{\tilde \zeta_{i}}> 0$ on $S$,\footnote{For this paper, the tropical part of a function is $\mathbb R$-valued instead of $\e{\mathbb R}$-valued, so $\totb{c\e a}:=a$. One complication is that  $\e {}$ is infinitesimal, so $a>0$ is equivalent to $\e a<\e 0$.}
 
 \item or $e_{S}\zeta_{i}=\zeta_{i}$ and  $\zeta_{i}$ is nowhere $0$ on the stratum of $\et mP$ corresponding to S, and $\totb{\tilde \zeta_{i}}=0$ on $S$.
 \end{itemize}
 
  The operation $e_{S}$ on a continuous function $h$ on $\et mP$ is then  \[e_{S}h(\zeta_{1},\dotsc, \zeta_{n})=h(e_{S}\zeta_{1},\dotsc, e_{S}\zeta_{n})\ .\] Of course, this implies that if $h$ is smooth or continuous,  $e_{S}h$ is too.

 Note that the operations $e_{S_{i}}$ commute and $e_{S_{i}}e_{S_{i}}=e_{S_{i}}$. More generally, $e_{S_{i}}e_{S_{j}}=e_{S'}$ where $S'$ is the smallest stratum of $P$ whose closure contains both $S_{i}$ and $S_{j}$.
 
 \begin{defn}If $I$ is  a nonempty collection of strata $\{S_{1},\dotsc, S_{n}\}$, use the following notation.
\[e_{I}h:=e_{S_{1}}\lrb{e_{S_{2}}\lrb{\dotsb e_{S_{n}}h}}\]
\[\Delta_{I}h:=\lrb{\prod_{S_{i}\in I}(\id-e_{S_{i}})}h\]
For an empty collection of strata, define $e_{\emptyset}h=h$ and $\Delta_{\emptyset}h=h$.
\end{defn}
For example,

\[\begin{split}\Delta_{\{S_{1},S_{2}\}}h(z_{1},z_{2})&:=(1-e_{S_{1}})(1-e_{S_{2}})h(z_1,z_2)
\\&:=h(z_1,z_2)-h(0,z_2)-h(z_1,0)+h(0,0)\end{split}\]
Note that if $S\in I$, $e_{S}\Delta_{I}=0$, so $\Delta_I$ us an operator for constructing functions that vanish on all strata in $I$. In the above example, this corresponds to \[\Delta_{\{S_{1},S_{2}\}}h(z_{1},0)=0=\Delta_{\{S_{1},S_{2}\}}h(0,z_{2})\ .\]

\subsection{Allowable coordinate charts}

\

We will define our norms using a class of  allowable coordinate charts on a family of curves. These coordinate charts are in a rather rigid form so that $\Delta_{S}$ is defined, and to simplify gluing analysis and facilitate defining a global version $\tilde \Delta_{S}$ of $\Delta_{S}$. 

\begin{defn}\label{allowable}
An allowable  coordinate chart on a family of exploded curves $\hat {\ex C}\xrightarrow{\pi_{\ex F}} \ex F$ is a coordinate chart on $\hat {\ex C}\xrightarrow {\pi_{\ex F}}\ex F$ satisfying the following requirements:

\begin{itemize}
\item The coordinate chart on $\ex F$ is some open subset $U\subset\mathbb R^{k}\times\et mP$ so that $e_{S}h$  is defined on $U$ for all strata $S$ if $h$ is defined on $U$, and so that $U$ is contained in a compact\footnote{Topological notions refer to the topology induced from the smooth part of an exploded manifod; see \cite{iec}.} subset of $\mathbb R^{k}\times \et mP$. 
\item The coordinate chart on $\hat{\ex C}$ and the map $\pi_{\ex F}$ is some restriction of a standard projection 
\[\mathbb R^{k}\times \et {m+1}Q\xrightarrow{\pi}\mathbb R^{k}\times \et m P\]
so that 
\begin{itemize}\item in standard coordinates, $\pi$ is given by
\[(x,\tilde z_{1},\dotsc,\tilde z_{m+1})=(x,\tilde z_{2},\dotsc,\tilde z_{m+1})\]
\item the polytope $Q=\totb{\mathbb R^{k}\times \et{m+1}Q}$ is defined by the equations 
\[(\totb{\tilde z_{2}},\dotsc,\totb {\tilde z_{m+1}})\in P\]
\[\totb{\tilde z_{1}}\geq 0\]
and possibly 
\[\totb{\e a\tilde z^{(-1,\beta_{2},\dotsc \beta_{m+1})}}\geq 0\ ;\]
in this case, let $\tilde z^{\beta}:=\e a\tilde z^{(-1,\beta_{2},\dotsc \beta_{m+1})}$.
\end{itemize}
\item The coordinate chart  $\tilde U$ on $\hat {\ex C}$ is in one of the following forms:
\begin{enumerate}[I]
\item \label{interesting} If the polytope $Q$ is not a product of $ P$ with $[0,\infty)$, then 
\[\tilde U:=\pi^{-1}U\cap\{\{\abs{\tilde z_{1}}<c, \abs{\tilde z^{\beta}}<c\} \}\]
and  $\abs{\tilde z_{1}\tilde z^{\beta}}<c^{2}/16$ on $U$.
\item \label{boring} If the polytope $Q$ is  $P\times [0,\infty)$, $\tilde U$ is in one of the following forms: 
\begin{enumerate}
\item $\tilde U$ is the intersection of $\pi^{-1}U$ with the set where the coordinate  $\tilde z_{1}$ on $\et 1{[0,\infty)}$ takes values in some open ball compactly contained in $\mathbb C^{*}\e 0$;
\item or $\tilde U$ is the intersection of $\pi^{-1}U$ with the set where $\abs{\tilde z_{1}}<c$. 
\end{enumerate}
\end{enumerate}

\end{itemize}

\end{defn}  

The nodes of curves in our family are covered by charts of type \ref{interesting}. We will be doing gluing analysis in these charts --- the rigid conditions on the shape of such a chart are artifacts of our chosen gluing construction; for example, the condition $\abs{\tilde z_{1}\tilde z^{\beta}}<c^{2}/16$ is there to ensure that some cutoff functions do not interact.
Such an allowable coordinate chart always exists around any point in $\ex C\longrightarrow \ex F$. The condition on its tropical part follows from the requirement that $Q\longrightarrow P$ be surjective, and have derivative (restricted to each stratum of $Q$) surjective on integral vectors; see Definition 10.1 of \cite{iec}.

\subsection{Weight functions $w_{S}$ and $w_{0}$}

\

Our norms  generalize the notion (on manifolds with cylindrical ends)  of Sobolev spaces with exponential weights. On a manifold with a cylindrical end, we only need a single weight function; this situation corresponds to a coordinate chart modeled on $\mathbb R^{k}\times \et 1{[0,\infty)}$. For more general coordinate charts $\et mP$, we need lots of weight functions. In particular, we need a weight function $w_{\mathcal S}$ for any set of strata $\mathcal S$ of $P$. This weight function should vanish on all these strata, but decay at least as slowly as any smooth function vanishing on the same strata. We need to use these weight functions when analysing a family of vertical differential operators on an allowable coordinate chart, $\tilde U\longrightarrow  U$; for this we must control terms in $w_S$ that are not lifted from $U$ by using the algorithm below.

\begin{defn} Given any allowable coordinate chart $\tilde U$, and  any collection of strata $\mathcal S:=\{S_{j}\subset \totb{\tilde U}\subset Q\}$, define a weight function $w_{\mathcal S}$ as follows: choose a set of generators  $\{\zeta_{i}\}$ for the monomial ideal\footnote{This monomial ideal is an ideal within the monoid of smooth monomials on $\et {m+1}Q$; the set $\{\zeta_i\}$ generate this monomial ideal in the sense that  any smooth monomial $\zeta$ with $e_{S_j}\zeta=0$ for all $S_j\in \mathcal S$ is the product of some $\zeta_i$ with  another smooth monomial.} of smooth monomials $\zeta$ on $\et {m+1}Q$ so that $e_{S_{j}}\zeta=0$ for all $S_{j}\in\mathcal S$.  Choose these generators for our monomial ideal using the following algorithm: Choose a finite set of generators for the monoid of smooth monomials   on  $\et {m+1}Q$, a set consisting  of $\totl{\tilde z_{1}}, \totl{\tilde z^{\beta}}$ and lifts of smooth monomials on $\et mP$. Include in our set $\{\zeta_{i}\}$ all the generators which are in our monomial ideal, then all products of two unused generators  in our monomial ideal. Continue this, adding all products of $k$ generators in our monomial ideal, but with no  sub-product in our monomial ideal.

Now define
\[w_{\mathcal S}:=\lrb{\sum \abs{\zeta_{i}}}\ .\]

 \end{defn}
 
 This weight function $w_{\mathcal S}$ vanishes on all  strata in $ \mathcal S$, and has the following property: If $h$ is any smooth function, then for any $\delta<1$,  the function $w^{-\delta}\Delta_{\mathcal S}h$  extends to a continuous function that vanishes on the strata in $\mathcal S$. So,  the smooth function $\Delta_{\mathcal S} h$ vanishes on the strata in $\mathcal S$, and decays faster than $w^{\delta}_{\mathcal S}$ for any $\delta<1$; see Lemma 7.4 of \cite{iec}.
 
 \
 
 A special case of $w_{\mathcal S}$: let $\mathcal S_{0}$ be the set of strata  not projecting homeomorphically onto strata of $\totb U$. These are the strata corresponding to the edges of the exploded curves in our family.
 
For example, in the case  of coordinate charts in the form of (\ref{interesting}) above,  then
\[w_{\mathcal S_{0}}:=\lrb{\abs{\totl{\tilde z_{1}}}+\abs{\totl{\tilde z^{\beta}}}+\sum_{i}\abs{\zeta_{i}}}\]  where $\{\zeta_{i}\}$ is some finite set of generators for the monomial ideal of smooth monomials on $U$ that vanish on all strata on which $\totl{\tilde z_{1}\tilde z^{\beta}}$ vanishes. 
 This weighting function $w_{\mathcal S_{0}}$ satisfies the following: if $h$ is a smooth or $\C\infty1$ function on $\et {m+1}Q$ vanishing where $\totl{\tilde z_{1}}=\totl{\tilde z^{\beta}}=0$, then for all $\delta<1$, $hw_{\mathcal S_{0}}^{-\delta}$ extends to a continuous function on $\et {m+1}Q$, also vanishing where $\totl{\tilde z_{1}}=\totl{\tilde z^{\beta}}=0$. Another way to say this is as follows: If  $e_{ S}h=0$ for all $S$ in $\mathcal S_{0}$, then $hw_{\mathcal S_{0}}^{-\delta}$ extends to a continuous function so that $e_{ S}(hw_{\mathcal S_{0}}^{-\delta})=0$.

\

For uniform control over the norm of a cutting map used in gluing analysis,  we replace $w_{\mathcal S_{0}}$ with a simpler weighting function  $w_{0}$ that does not include these extra $\abs{\zeta_{i}}$ terms.

\begin{defn} On allowable coordinate charts of type \ref{interesting}, define
\[w_{0}:=\abs {\totl{\tilde z_{1}}}+\abs{\totl{\tilde z^{\beta}}}\]
and on allowable coordinate charts of type \ref{boring}, define
\[w_{0}:=\abs {\totl{\tilde z_{1}}}\ .\]

\end{defn}

The following lemma ensures that repacing $w_{\mathcal S_{0}}$ with $w_{0}$ is not too dangerous. 
\begin{lemma}\label{k comparison}Given any allowable coordinate chart, there exists a positive number ${q} $,  and  a constant $c$ so that
\[cw_{0}\geq  w_{\mathcal S_{0}}^{ {q}}\ .\]\end{lemma}

\pf

On allowable coordinate charts of type (\ref{boring}), $w_{0}=w_{\mathcal S_{0}}$, so we can reduce to the case of charts of type (\ref{interesting}). On such a chart, $w_{\mathcal S_{0}}:=\lrb{\abs{\totl{\tilde z_{1}}}+\abs{\totl{\tilde z^{\beta}}}+\sum_{i}\abs{\zeta_{i}}}$, where $\{\zeta_{i}\}$ is some finite generating set  for the  smooth monomials vanishing  where $\totl{\tilde z_{1}\tilde z^{\beta}}=0$.  Each such $\zeta_{i}$ is the smooth part of some exploded monomial $\tilde \zeta_{i}$ whose tropical part is a nonnegative integral affine function $\totb{\tilde \zeta_{i}}:P\longrightarrow [0,\infty)$, strictly positive on every strata where $\totb{\tilde z_{1}\tilde z^\beta}$ is strictly positive.  If there are no strata on which $\totb{\tilde \zeta_{i}}=0$, then $\zeta_{i}$ is identically $0$, and can be safely discarded. Otherwise, there exists some positive integer ${q}$ so that 
\begin{equation}\label{kest} {q}\totb{\tilde \zeta_{i}}\geq \totb{\tilde z_{1}\tilde z^\beta}\ .\end{equation}
Such an estimate holds on any ray within $P$ starting at a point on the largest stratum on which $\totb{\tilde \zeta_{i}}$ vanishes, and holds with a uniform constant because the set of rays normal to this stratum at a given point is compact. As $P$ is contained within the span of vectors within this strata and such normal rays, it follows that the above estimate (\ref{kest}) holds for a uniform ${q}$ on all of $P$. As $\{\zeta_{i}\}$ is finite, we can choose ${q}$ so that (\ref{kest}) holds for all $i$. It follows that  the product of any ${q}$ such $\tilde \zeta_{i}$ will have tropical part bounded by $\totb{\tilde z_{1}\tilde z^\beta}$, and therefore be the product of $\tilde z_{1}\tilde z^\beta$ with a monomial $\tilde \zeta$ whose tropical part is nonnegative  on $P$. The corresponding product of $\zeta_{i}$ will therefore be $\totl{\tilde \zeta}\totl{\tilde z_{1}}\totl{\tilde z^\beta}$, and on our coordinate chart, its absolute value will be bounded by a constant times $\abs{\totl{\tilde z_{1}}}$. Similarly any term appearing in $w_{\mathcal S_{0}}^{{q}}$ will be bounded by a constant times either $\abs{\totl{\tilde z_{1}}}$ or $\abs{\totl{\tilde z^\beta}}$. As there are only finitely many such terms, our required estimate follows.

\[cw_{0}\geq  w_{\mathcal S_{0}}^{ {q}}\]

\stop

Lemma \ref{k comparison} implies that given any smooth or $\C\infty1$ function $h$ and $\delta'<\frac 1{{q}}$, $w_{0}^{-\delta'}\Delta_{\mathcal S_{0}}h$ extends to a continuous function vanishing on all strata in $\mathcal S_{0}$.

\subsection{Norms in allowable coordinate charts}\label{norms}

\

Given a finite collection of allowable coordinate charts,  we now define a series of norms on vector-valued functions  on these coordinate charts $\tilde U\subset \mathbb R^{n}\times \et {m+1}Q$. In all that follows, use the standard metric in which $\frac \partial{\partial x_{i}}$ and the real and imaginary parts of $\tilde z_{i}\frac \partial{\partial \tilde z_{i}}$ are orthonormal. In the case of a chart on a single curve, these norms $\norm\cdot_{k,\delta}$ and $\norm\cdot_{k-1,\delta}^{1}$ are analogous to $L_{k}^{p}$ with exponential weights on cylindrical ends. In the more general case of a chart on a family of curves, these norms are also analogous weighted $L_{k}^{p}$ norms, however we treat vertical directions specially, because we are analyzing a family of vertical differential operators. In particular, the norm $\norm\nu_{k,\delta}^{1}$ controls $k$ derivatives of $\nu$ in addition to one vertical derivative. 

\begin{enumerate}
\item Choose some exponent $p>2$ and weight $0<\delta<1$. Also,  choose some number $K$  so  that, for each coordinate chart, Lemma \ref{k comparison} holds for some  ${q}<K$. Then let 
\[\delta':=\frac 1K(1-\delta)\]
and define the following norm.
\[\norm\nu_{\delta'}:=\sup_{x\in\pi (\tilde U)} \lrb{\int_{\pi^{-1}(x)\cap \tilde U}\abs{w_{0}^{-\delta'}\nu}^{p}}^{\frac 1p}\]
So $\norm\nu_{\delta'}$ is the supremum, over all fibers,  of a weighted $L^{p}$ norm of $\nu$ restricted to a fiber.
This is only defined when $\Delta_{\mathcal S_{0}}\nu=\nu$ --- so $e_{S}\nu=0$ for all $S\in\mathcal S_{0}$, or equivalently $\nu$ vanishes on edges of the curves in our family.

\item 

\[\norm\nu_{0,\delta}:=\norm{\nu}_{\delta'}+ \max_{\mathcal S} \lrb{\norm {w_{\mathcal S}^{-\delta} \Delta_{\mathcal S}\nu }_{\delta'}}\]

The maximum is taken over all collections $\mathcal S$ of substrata of $\totb{\tilde U}$. (Where $w_{\mathcal S}$ and $\Delta _{\mathcal S}\nu$ both vanish, we use the convention that $w_{\mathcal S}^{-\delta} \Delta_{\mathcal S}\nu =0$ for the purposes of calculating these norms.)
 \item
 
 \[\norm \nu_{k, \delta}:= \norm \nu_{k-1, \delta}+\norm{d\nu}_{k-1, \delta}\]
 
 These norms should be thought of as suitable generalizations of $L^{p}_{k}$ with exponential weights.
 \item For this last collection of  norms, we shall use the notation $\dvert\nu$ to refer to $d\nu$ restricted to the vertical tangent space --- spanned by the real and imaginary parts of $\tilde z_{1}\frac\partial{\partial\tilde z_{1}}$.

 \[\norm\nu^{1}_{ \delta'}:=\sup \abs\nu +\norm{\dvert\nu}_{ \delta'}\]
 \[\norm\nu^{1}_{0, \delta}:= \norm\nu^{1} _{\delta'}+\max_{\mathcal S} \lrb{\sup\abs {w_{\mathcal S}^{-\delta} \Delta_{\mathcal S}\nu }+\norm{w_{\mathcal S}^{-\delta} \Delta_{\mathcal S}\dvert\nu }_{\delta'}}\]
\[\norm\nu^{1}_{k, \delta}:=\norm \nu^{1}_{k-1, \delta}+\norm{d\nu}^{1}_{k-1, \delta} \]
 
\end{enumerate}

\

\

Recall the following definition of $C^{k,\delta}$ for any $0<\delta<1$.

\begin{defn}[$C^{k,\delta}$ and $\C\infty\delta$ regularity]\label{Ckdelta}
Define $C^{0,\delta}$ to be the same as $C^0$. A sequence of smooth functions $f_i\in C^\infty(\mathbb R^n\times\et mP)$ converge to a continuous function $f$ in  $C^{k,\delta}(\mathbb R^n\times \et mP)$ if the following conditions hold:
\begin{enumerate}
\item Given any collection $I$ of at most $k$ nonzero strata, the sequence of  functions
\[\abs{w_{I}^{-\delta}\Delta_{I}(f_i-f)}\]
converges to $0$ uniformly on compact subsets of $\totl{\mathbb R^{n}\times\et mP}$ as $i\rightarrow\infty$. (This includes the case where our collection of strata is empty and $f_{i}\rightarrow f$ uniformly on compact subsets.)
\item For any smooth vectorfield $v$, $v(f_i)$ converges to some function $vf$ in $C^{k-1,\delta}$.
\end{enumerate} 
Define $C^{k,\delta}(\mathbb R^n\times \et mP)$ to be the closure of $C^\infty$ in $C^0$ with this topology.\footnote{This notion of convergence induces a well-behaved topology because it is equivalent to convergence in some countable sequence of norms. In particular, a subset is closed in the induced topology if and only if it is sequentially closed.} \\Define $C^{\infty,\delta}$ to be the intersection of $C^{k,\delta}$ for all $k$.\\ Define $\C\infty\delta$ to be the intersection of $C^{\infty,\delta'}$ for all $\delta'<\delta$.
\end{defn}

Note that  any $\C\infty1$ function $\nu$ restricted to any  compact subset will have $\norm\nu^{1}_{k, \delta}$ finite, and if $\Delta_{\mathcal S_{0}}\nu=\nu$, then $\norm\nu_{k,\delta}$ will also be finite. Standard Sobolev estimates imply that if $\norm{\nu}_{k,\delta}$ or $\norm\nu_{k,\delta}^{1}$ is finite for all $k$, and all $\delta<1$ then $\nu$ is $\C\infty1$.

\

We shall often need the following observations about the behavior of products in our norms:
\begin{lemma}\label{product lemma}

\

\begin{enumerate}[A]
\item\label{product1} For $I_1$ and $I_2$ any collection of strata of an  allowable coordinate chart,   $w_{I_{1}}w_{I_{2}}w^{-1}_{I_{1}\cup I_{2} }$ is bounded. 
\item\label{product2}The following product formula for $\Delta_{\mathcal S}$ holds.\footnote{Throughout this paper, the symbol $\subset$ means $\subseteq$. Proper subsets will be indicated by the symbol $\subsetneq$.} \[\Delta_{\mathcal S}(\phi\psi)=\sum_{I\subset\mathcal S}\lrb{ e_{\mathcal S-I}\Delta_{I}\phi}\Delta_{\mathcal S-I}\psi\]

\item \label{product3}
On any allowable coordinate chart and for any $k\in\mathbb N$, there exists some constant $c$, depending only on $k$ and our choices defining the above norms, so that, if at least one of $\phi$ or $\psi$ is a real valued function,
 \[\norm{\phi\psi}_{k, \delta}\leq c\norm\phi_{k, \delta}\norm{\psi}_{k, \delta}^{1}\ .\]
\item\label{product4}
On any  allowable coordinate chart and and for any $k\in\mathbb N$, there exists some constant $c$, depending only on $k$ and our choices defining the above norms, so that, if at least one of $\phi$ or $\psi$ is a real valued function,
 \[\norm{\phi\psi}^{1}_{k, \delta}\leq c\norm\phi^{1}_{k, \delta}\norm{\psi}_{k, \delta}^{1}\ .\]

\end{enumerate}
\end{lemma}
\pf

To see item \ref{product1}, note that the $w_{I_{i}}$ is a finite sum of absolute values of smooth monomials vanishing on strata in $I_{i}$, so $w_{I_{1}}w_{I_{2}}$ is a finite sum of absolute values of smooth monomials  vanishing on strata in ${I_{1}\cup I_{2}}$. The weight function $w_{I_{1}\cup I_{2}}$ is a sum of absolute values of generators for the monomial ideal of such smooth monomials, so on $\mathbb R^{n}\times \et mQ$, every item in the former sum is a continuous function times an element of the latter sum. As allowable coordinate charts are always compactly contained inside $\mathbb R^{n}\times \et mQ$, it follows that $w_{I_{1}}w_{I_{2}}w^{-1}_{I_{1}\cup I_{2} }$ is bounded.

\

To prove item \ref{product2}, note first that if $S$ denotes a single stratum, 
\[\Delta_{S}\phi\psi=\phi\psi-(e_{S}\phi)(e_{S}\psi)=(\Delta_{S}\phi)\psi+(e_{S}\phi)\Delta_{S}\psi\]
so the required identity holds if $\mathcal S$ consists of a single stratum. Suppose now that $\mathcal S=\mathcal S'\cup\{S\}$, and the required identity holds for $\mathcal S'$. Then,
\[\begin{split}\Delta_{\mathcal S}\phi\psi&=\Delta_{S}\lrb{\sum_{I\subset\mathcal S'}\lrb{ e_{\mathcal S'-I}\Delta_{I}\phi}\Delta_{\mathcal S'-I}\psi}
\\&=\sum_{I\subset\mathcal S'}\lrb{ e_{\mathcal S'-I}\Delta_{S}\Delta_{I}\phi}\Delta_{\mathcal S'-I}\psi
+\lrb{ e_{S}e_{\mathcal S'-I}\Delta_{I}\phi}\Delta_{S}\Delta_{\mathcal S'-I}\psi
\\&=\sum_{I\subset\mathcal S}\lrb{ e_{\mathcal S-I}\Delta_{I}\phi}\Delta_{\mathcal S-I}\psi
\ .\end{split}\]
So by induction, the required identity holds for any set of strata $\mathcal S$. 

\

To prove item \ref{product3}, note
\[\norm{\phi\psi}_{\delta'}\leq \norm{\phi}_{\delta}\sup \abs\psi\leq\norm{\phi}_{\delta'}\norm{\psi}^{1}_{\delta'}\ .\]
Next, we use item \ref{product2}, then item \ref{product1}, then the above observation to show the following.

\[\begin{split}\norm{w_{\mathcal S}^{-\delta}\Delta_{\mathcal S}(\phi\psi)}_{\delta'}&\leq \sum_{I\subset \mathcal S}
\norm{w_{\mathcal S}^{-\delta}\lrb{e_{\mathcal S-I}\Delta_{I}\phi}\Delta_{\mathcal S-I}\psi}_{\delta'}
\\&\leq c\sum_{I\subset \mathcal S}
\norm{w_{I}^{-\delta}\lrb{e_{\mathcal S-I}\Delta_{I}\phi}w^{-\delta}_{\mathcal S-I}\Delta_{\mathcal S-I}\psi}_{\delta'}
\\&\leq c\sum_{I\subset \mathcal S}\norm{w_{I}^{-\delta}e_{\mathcal S-I}\Delta_{I}\phi}_{\delta'}
\sup\abs{w^{-\delta}_{\mathcal S-I}\Delta_{\mathcal S-I}\psi} 
\end{split}\]

The constant $c$ above depends on the collection of strata $\mathcal S$, but there are only a finite number of strata, and because $\Delta_{S}\Delta_{S}=\Delta_{S}$, we need only consider sets of distinct strata. Therefore, for a different constant $c$,\footnote{Throughout this paper, constants $c$ in such inequalities will depend on the choice of coordinate chart, norm and maybe a collection of strata, but will they will never depend on $\phi$ or $\psi$. There are many inequalities, and we will reuse $c$, even though it refers to a different constant in different inequalities.}
\[\norm{\phi\psi}_{0,\delta}\leq c\norm{\phi}_{0,\delta}\norm{\psi}_{0,\delta}^{1}\ .\]
 Now we can use induction on the number of derivatives: Suppose that the required inequality holds for $k-1$ derivatives.
\[\begin{split}\norm{\phi\psi}_{k,\delta}&=\norm{\phi\psi}_{k-1,\delta}+\norm{d(\phi\psi)}_{k-1,\delta}
\\&\leq c'(\norm\phi_{k-1,\delta}\norm{\psi}_{k-1,\delta}^{1}+\norm{d\phi}_{k-1,\delta}\norm{\psi}^{1}_{k-1,\delta}
+\norm{\phi}_{k-1,\delta}\norm{d\psi}^{1}_{k-1,\delta})
\\&\leq 3c'\norm\phi_{k,\delta}\norm\psi_{k,\delta}^{1}
\end{split}\]
Therefore, by induction the required inequality holds for all $k$.

\

We prove item \ref{product4} similarly as follows:
\[\sup\abs{\phi\psi}\leq \sup\abs\phi\sup\abs\psi\]
so as above, we may estimate using item \ref{product2} and item \ref{product1}.

\[\begin{split}\sup\abs{w_{\mathcal S}^{-\delta}\Delta_{\mathcal S}(\phi\psi)}&\leq \sum_{I\subset\mathcal S}
\sup\abs{w_{\mathcal S}^{-\delta}\lrb{e_{\mathcal S-I}\Delta_{I}\phi}\Delta_{\mathcal S-I}\psi}
\\&\leq c\sum_{I\subset \mathcal S}
\sup\abs{w_{I}^{-\delta}\lrb{e_{\mathcal S-I}\Delta_{I}\phi}w^{-\delta}_{\mathcal S-I}\Delta_{\mathcal S-I}\psi}
\\&\leq c\sum_{I\subset \mathcal S}\sup\abs{w_{I}^{-\delta}e_{\mathcal S-I}\Delta_{I}\phi}
\sup\abs{w^{-\delta}_{\mathcal S-I}\Delta_{\mathcal S-I}\psi} 
\end{split}\]

We also have from item \ref{product3},

\[\begin{split}\norm{\dvert(\phi\psi)}_{k,\delta}&\leq \norm{(\dvert\phi)\psi}_{k,\delta}+\norm{\phi\dvert\psi}_{k,\delta}
\\&\leq c\norm{\dvert\phi}_{k,\delta}\norm{\psi}_{k,\delta}^{1}
+c\norm\phi_{k,\delta}^{1}\norm{\dvert\psi}_{k,\delta} 
\\&\leq 2c\norm{\phi}_{k,\delta}^{1}\norm{\psi}_{k,\delta}^{1}\ .
\end{split}\]

Therefore, for some new constant $c$ we get the required inequality.
\[\norm{\phi\psi}_{0,\delta}^{1}\leq c\norm\phi_{0,\delta}^{1}\norm\psi_{0,\delta}^{1} \]
The general case now follows by induction because if it holds for $k-1$,

\[\begin{split}\norm{\phi\psi}^{1}_{k,\delta}&=\norm{\phi\psi}^{1}_{k-1,\delta}+\norm{d(\phi\psi)}^{1}_{k-1,\delta}
\\&\leq c'(\norm\phi^{1}_{k-1,\delta}\norm{\psi}_{k-1,\delta}^{1}+\norm{d\phi}^{1}_{k-1,\delta}\norm{\psi}^{1}_{k-1,\delta}
+\norm{\phi}^{1}_{k-1,\delta}\norm{d\psi}^{1}_{k-1,\delta})
\\&\leq 3c'\norm\phi^{1}_{k,\delta}\norm\psi_{k,\delta}^{1}\ .
\end{split}\]

\stop

\

We use the adjective extendable to grant non-compact sets some compact-set superpowers. In particular:

\begin{defn}\label{extendable}
\begin{enumerate}

\item An extendable open subset of a topological space is an open subset  contained inside a compact subset of that space.
\item An extendable allowable coordinate chart is an allowable coordinate chart  that is also an extendable subset of some larger allowable coordinate chart.
\item An extendable function is a function which is defined on an extendable subset, and which is the restriction of a function defined on a larger, compact domain.
\item An extendable vectorbundle is a vectorbundle which is the restriction of some vectorbundle to an extendable subset. 
\item An extendable function on an extendable vectorbundle is the restriction of some function to an extendable vectorbundle.  
 \end{enumerate}

\end{defn}

\begin{lemma}\label{composition}
If, on some extendable allowable coordinate chart $\tilde U$, $\nu$ is a $\mathbb R^{n}$--valued function with $\norm\nu_{k,\delta}^{1}$ finite and $E$ is a $\C\infty1$ extendable function on the extendable vectorbundle  $\mathbb R^{n}\times \tilde U$, then $\norm{E(\nu)}^{1}_{k,\delta}$ is bounded. If $\Delta_{\mathcal S_{0}} E=E$, then $\norm{E(\nu)}_{k,\delta}$ is bounded. These bounds can be chosen to depend continuously on $\nu$ in the $\norm\cdot^{1}_{k,\delta}$ topology.\end{lemma}

\pf

For a single stratum, $\Delta_{S}(E(\nu))=(\Delta_{S}E)(\nu)+\Delta_{S}\lrb{(e_{S}E)(\nu)}$, so if $\mathcal S$ is a collection of strata, 
 \begin{equation}\label{comp1}\Delta_{\mathcal S}(E(\nu))=\sum_{I\subset\mathcal S}\Delta_{I}\lrb{\lrb{e_{I}\Delta_{\mathcal S-I}E}(\nu)}\ .
 \end{equation} 

Because $\norm\nu_{k,\delta}^{1}$ is bounded, $\sup\abs{\nu}$ is bounded, so, for $E(\nu)$, we may restrict attention to 
a fiberwise-bounded subset of $\mathbb R^{n}\times \tilde U\longrightarrow \tilde U$. Here, we have bounds on $E$ and all its derivatives, because $E$ and $\tilde U$ are extendable.  In the following, let $D_{I}$ indicate the derivative with respect to $\frac \partial{\partial t_{i}}$ for all $i\in I$.
\begin{equation}\label{comp2}\Delta_{I}\lrb{(e_{I}E)(\nu)}=\int_{0}^{1}\dotsb\int_{0}^{1} D_{I}\lrb{(e_{I}E)\lrb{\prod_{i\in I}(e_{S_{i}}+t_{i}\Delta S_{i})\nu}}\prod_{i\in I} dt_{i} \end{equation}
To estimate $\Delta_{I}(e_{I}E)(\nu)$, estimate this integrand. Use the notation $\phi_{I}:=\prod_{i\in I}(e_{S_{i}}+t_{i}\Delta S_{i})\nu$.
\begin{equation}\label{comp3}\begin{split}D_{I}(e_{I}E)(\phi_{I})&=\sum_{\coprod_{j=1}^{n}I_{j}=I}(D^{n}e_{I}E)(D_{I_{1}}\phi_{I})\dotsb(D_{I_{n}}\phi_{I})
\\&=\sum_{\coprod_{j=1}^{n}I_{j}=I}(D^{n}e_{I}E)(\Delta_{I_{1}}\phi_{I-I_{1}})\dotsb(\Delta_{I_{n}}\phi_{I-I_{n}})
\end{split}\end{equation}

The sum above is over all partitions of $I$. The above equations (\ref{comp2}) and (\ref{comp3}) also hold with $\Delta_{\mathcal S-I}E$ replacing $E$; note also that $D^{n}e_{I}\Delta_{\mathcal S-I}E=\Delta_{S-I}e_{I}D^{n}E$.  Using the equations (\ref{comp1}), (\ref{comp2}) and (\ref{comp3})  and Lemma \ref{product lemma} part \ref{product1}, we get that  
$\sup\abs{w_{\mathcal S}^{-\delta}\Delta_{\mathcal S}(E(\nu))}$ is bounded  by a constant times the following expression.
\begin{equation}\sum_{I\subset\mathcal S}\sum_{\coprod_{j=1}^{n}I_{j}=I}\sup\abs{w_{\mathcal S-I}^{-\delta}\Delta_{\mathcal S-I}e_{I}(D^{n}E)(\phi_{I})}\prod_{j=1}^{n}\sup\abs{w_{I_{j}}^{-\delta}
\Delta_{I_{j}}\phi_{I-I_{j}}}\end{equation}

The first term in each of the  above summands is finite because $E$ is extendable and in $\C\infty1$, and the other terms are bounded by $\norm{\nu}_{0,\delta}^{1}$. Our estimate for each of these terms can be chosen to depend continuously on $\nu$ in the $\norm\cdot_{0,\delta}^{1}$ topology. (The estimate of the first term can be chosen continuous in the supremum topology, which is weaker than the $\norm{\cdot}_{0,\delta}^{1}$ topology.) 

Similar to the above, if $\Delta_{S_{0}}E=E$, we may bound $\norm{w_{\mathcal S}^{-\delta}\Delta_{\mathcal S}(E(\nu))}_{\delta'}$ by a constant times
\begin{equation}\label{gngn}\sum_{I\in\mathcal S}\sum_{\coprod_{j=1}^{n}I_{j}=I}\sup\abs{w_{0}^{-\delta'-\epsilon}w_{\mathcal S-I}^{-\delta}\Delta_{\mathcal S-I}e_{I}(D^{n}E)(\phi_{I})}\prod_{j=1}^{n}\sup\abs{w_{I_{j}}^{-\delta}
\Delta_{I_{j}}\phi_{I-I_{j}}}\end{equation}

for some $\epsilon>0$. The first term is bounded for $\epsilon$ small enough using Lemma \ref{k comparison}: 
\[\begin{split}\sup\abs{w_{0}^{-\delta'-\epsilon}w_{\mathcal S-I}^{-\delta}\Delta_{\mathcal S-I}(D^{n}E)}
&\leq c\sup\abs{w_{\mathcal S_{0}}^{-{q}(\delta'+\epsilon)}w_{\mathcal S-I}^{-\delta}\Delta_{\mathcal S-I}(D^{n}E)}
\\ &\ \ \  =c\sup\abs{w_{\mathcal S_{0}}^{-{q}(\delta'+\epsilon)}w_{\mathcal S-I}^{-\delta}\Delta_{\mathcal S_{0}\cup (\mathcal S-I)}(D^{n}E)}
\\ &\ \ \  \leq c'\sup\abs{w_{\mathcal S_{0}\cup(\mathcal S-I)}^{-{q}(\delta'+\epsilon)-\delta}\Delta_{\mathcal S_{0}\cup (\mathcal S-I)}(D^{n}E)}
\end{split}\]

We may choose $\epsilon$ and ${q}$ so that the above inequality holds and ${q}(\delta'+\epsilon)+\delta<1$, so the final term in the inequality above is bounded. The other terms on the  righthand side of (\ref{gngn}) are bounded by $\norm{\nu}_{0,\delta}^{1}$. Again, the  bounds can be chosen continuous in the $\norm\cdot_{0,\delta}^{1}$ topology. Therefore, in this case $\norm{E(\nu)}_{0,\delta}$ is bounded, and the bound can be chosen to depend  continuously on $\nu$ in the $\norm\cdot_{0,\delta}^{1}$ topology.

In the case that $\Delta_{\mathcal S_{0}}E$ is not necessarily equal to $E$, we already know that  $\sup\abs{w_{\mathcal S}^{-\delta}\Delta_{\mathcal S}(E(\nu))}$ is bounded, and we must show that $\norm{\dvert(E(\nu))}_{0,\delta}$ is bounded. To this end, note that $\dvert(E(\nu))=DE(\nu)(\dvert\nu)+(\dvert E)(\nu)$. The second term can be dealt with by observing that $\dvert E$ is in $\C\infty1$ and $\Delta_{S_{0}}\dvert E=\dvert E$, therefore, as argued above, $\norm{(\dvert E)(\nu)}_{0,\delta}$ is bounded. The first term can be dealt with in the same way as the product was dealt with in Lemma \ref{product lemma} part \ref{product3}. Note that in that argument, the only part of the norm $\norm\cdot_{0,\delta}^{1}$ used was the part involving the supremum. Following this argument, we can bound $\norm{DE(\nu)(\dvert\nu)}_{0,\delta}$ by the product of $\norm{\dvert\nu}_{0,\delta}$ with $\sum_{\mathcal S}\sup\abs{w_{\mathcal S}^{-\delta}\Delta_{\mathcal S}(DE(\nu))}$, which as argued above is bounded. Again, our bounds may be chosen to depend on $\nu$ continuously in the $\norm{\cdot}_{0,\delta}^{1}$ topology.

We have now shown that if $\norm\nu_{0,\delta}^{1}$ is bounded, $\norm{E(\nu)}_{0,\delta}^{1}$ is bounded, and if $\Delta_{\mathcal S_{0}}E=E$, then $\norm{E(\nu)}_{0,\delta}$ is bounded, and these  bounds can be chosen to  depend continuously on $\nu$ in the $\norm\cdot_{0,\delta}^{1}$ topology. For induction, suppose that the equivalent statement holds for $\norm{\cdot}_{k,\delta}^{1}$ and $\norm{\cdot}_{k,\delta}$. Note that $d(E(\nu))(\cdot)=DE(\nu)(0,d\nu(\cdot))+(DE)(\nu)(\cdot,0)$. If $\norm{\nu}_{k+1,\delta}^{1}$ is bounded, the first term is a composition of the $\C\infty1$ function $DE(\cdot)(0,\cdot)$ and  $(\nu,d\nu)$ which has $\norm{(\nu,d\nu)}_{k,\delta}^{1}$ bounded, and the second term is a composition of the $\C\infty1$ function $DE(\cdot)(\cdot,0)$ with $\nu$. By our inductive assumption,  the $\norm{\cdot}_{k,\delta}^{1}$ norm of both these terms is bounded, and $\norm{d(E(\nu))}_{k+1,\delta}^{1}$ is bounded.  Similarly, $\norm{d(E(\nu))}_{k+1,\delta}$ is bounded if $\norm{\nu}_{k+1,\delta}^{1}$ is bounded and $\Delta_{\mathcal S_{0}}E=E$, as in that case $\Delta_{\mathcal S_{0}}DE=DE$. All these bounds can be chosen to depend continuously on $\nu$ in the $\norm\cdot_{k+1,\delta}^{1}$ topology. By induction, we have proved the lemma for all $k$.

\stop

\begin{cor}\label{continuous composition}
If  $E$ is an extendable $\C\infty1$ function on $\mathbb R^{n}\times\tilde U$,  then  $E(\nu)$ in the $\norm\cdot^{1}_{k,\delta}$ topology depends continuously on $\nu$ in the $\norm{\cdot}^{1}_{k,\delta}$ topology. If $\Delta_{\mathcal S_{0}} E=E$, then  $E(\nu)$ in the $\norm{\cdot}_{k,\delta}$ topology depends continuously on $\nu$ in the $\norm\cdot_{k,\delta}^{1}$ topology.

\end{cor}

\

\pf \[E(\nu_{1})-E(\nu_{2})=\int_{0}^{1}DE(\nu_{1}+t(\nu_{2}-\nu_{1}))(\nu_{2}-\nu_{1})dt\]
so using Lemma \ref{product lemma} part \ref{product4},
\[\norm{E(\nu_{1})-E(\nu_{2})}^{1}_{k,\delta}\leq c \norm{\nu_{1}-\nu_{2}}_{k,\delta}^{1}\int_{0}^{1}\norm{DE(\nu_{1}+t(\nu_{2}-\nu_{1}))}^{1}_{k,\delta}dt\ .\]
As $DE$ is $\C\infty1$, Lemma \ref{composition} tells us that $\norm{DE(\nu')}^{1}_{k,\delta}$ can be bounded uniformly for $\nu'$ in a $\norm{\cdot}_{k,\delta}^{1}$--neighborhood of $\nu_{1}$, therefore, if $\nu_{2}$ is in this neighborhood,  
\[\norm{E(\nu_{1})-E(\nu_{2})}^{1}_{k,\delta}\leq c \norm{\nu_{1}-\nu_{2}}_{k,\delta}^{1}\ .\]
So, $E(\nu)$ in the $\norm\cdot^{1}_{k,\delta}$ topology depends continuously on $\nu$ in the $\norm\cdot_{k,\delta}^{1}$ topology.

Similarly, using Lemma \ref{product lemma} part \ref{product3},
\[\norm{E(\nu_{1})-E(\nu_{2})}_{k,\delta}\leq c\norm{\nu_{1}-\nu_{2}}_{k,\delta}^{1}\int_{0}^{1}\norm{DE(\nu_{1}+t(\nu_{2}-\nu_{1}))}_{k,\delta}dt\ .\]
As $DE$ is $\C\infty1$  and $\Delta_{\mathcal S_{0}}DE=DE$ if $\Delta_{\mathcal S_{0}}E=E$, $\norm{DE(\nu')}_{k,\delta}$ can be bounded uniformly for $\nu'$ in a $\norm{\cdot}_{k,\delta}^{1}$ neighborhood of $\nu_{1}$ if $\Delta_{\mathcal S_{0}}E=E$. Therefore, if $\nu_{2}$ is in this neighborhood,  
\[\norm{E(\nu_{1})-E(\nu_{2})}_{k,\delta}\leq c \norm{\nu_{1}-\nu_{2}}_{k,\delta}^{1}\ .\]
So, if $\Delta_{\mathcal S_{0}}E=E$,  $E(\nu)$ in the $\norm\cdot_{k,\delta}$ topology depends continuously on $\nu$ in the $\norm\cdot_{k,\delta}^{1}$ topology.
\stop

\subsection{Equivalent norms using lifted sets of strata}

\

We  now describe equivalent norms for $\norm\cdot_{0,\delta}$ and $\norm{\cdot}_{0,\delta}^{1}$. These new norms  only involve weighting functions $w_{\mathcal S}$ with $\dvert w_{\mathcal S}=0$ --- making estimates involving vertical derivatives easier.  For this, we  need the following concepts:

\begin{defn}
Given an allowable coordinate chart $\tilde U\xrightarrow{\pi}U$ and a set of strata $\mathcal S$ of $\totb{U}$, the lift, $\tilde{\mathcal S}$, of $\mathcal S$ is a set of strata of $\totb{\tilde U}$  defined as follows.
\[\tilde {\mathcal S}:=\{S\text{ so that }\totb{\pi}(S)\in\mathcal S\}\] 
For our purposes, two sets of strata $\mathcal S$ and $\mathcal S'$ will act identically if $\Delta_{\mathcal S}=\Delta_{\mathcal S'}$, so define a lifted set of strata $\mathcal S$ to be a set of strata of $\totb{\tilde U}$ with the following property: if for some strata $T$, the projection $\totb{\pi}( T)\in\totb\pi(\mathcal S)$, then $\Delta_{T}\Delta_{\mathcal S}=\Delta_{\mathcal S}$.

For any set of strata $\mathcal S$ in $\tilde U$, define the complement $\mathcal S^{c}$ to be the set of strata $S'$ so that $\Delta_{S'}\Delta_{\mathcal S}\neq\Delta_{\mathcal S}$, and $\totb{\pi}(S')=\totb{\pi}(S)$ for some $S\in\mathcal S$.
 
Use $\mathcal S^{c}=\emptyset$ when $\mathcal S$ is a lifted set of strata.
 \end{defn}

If $\mathcal S^{c}=\emptyset$, the vertical derivative $\dvert w_{\mathcal S}=0$ because $w_{\mathcal S}$ is the lift of some function from $U$. Note also that  $(\mathcal S\cup\mathcal S^{c})^{c}=\emptyset$.


\begin{lemma}\label{sobolev} On any allowable coordinate chart, given a stratum $S_{j}$ and $\bar S_{j}\in S_{j}^{c}$, there exists a constant $c>0$ so that
\[\abs{e_{\bar S_{j}}\Delta_{S_{j}}\phi}\leq c\norm{e_{\bar S_{j}}\dvert \phi}_{\delta}e_{\bar S_{j}}w_{0}^{\delta}\ .\]

\end{lemma}
\pf

The lefthand side of the above inequality is $\abs{e_{\bar S_{j}}\phi-e_{\bar S_{j}}e_{S_{j}}\phi}$. This  equals  the difference between the $\phi$ on some fiber of the coordinate chart and   $\phi$ on the edge contained in the same fiber. We shall bound this difference using a standard Sobolev estimate on this fiber. We will then get a uniform bound using the bounded geometry of allowable coordinate charts. 

Without losing generality, we may assume that the part of this fiber of interest  has coordinate $\tilde z_{1}$, and the smooth component of the fiber of interest is $\{0<\abs{\totl{\tilde z_{1}}}<c\}\subset \mathbb C$, where $c$ is a constant depending only on the coordinate chart, and not the particular fiber. Use cylindrical coordinates $\totl{\tilde z_{1}}=e^{t+i\theta}$.  Denote $\phi$ restricted to this smooth component of this fiber simply as $\phi$. 

We are interested in bounding $\abs{\phi(t,\theta)-\phi(-\infty,\theta)}$ in terms of $\norm{e_{\bar S_{j}}\dvert \phi}_{\delta}e_{\bar S_{j}}w_{0}^{\delta}$. In the case we have restricted ourselves to, $e_{\bar S_{j}}w_{0}=\abs {\totl{\tilde z_{1}}}=e^{t}$.  Using this observation and the definition of $\norm\cdot_{\delta}$, we have 
\[\norm{e_{\bar S_{j}}\dvert \phi}_{\delta}\geq \lrb{\int_{\{t<\log c\}}(e^{-\delta t}\abs{d\phi})^{p}dt d\theta}^{\frac 1p}\ .\]
(The inequality sign above is there because the  righthand side is considering only one fiber. The $\phi$ on the  righthand side is  $\phi$ restricted to this fiber, so $e_{\bar S_{j}}\dvert \phi=d\phi$.)
So long as $p>2$, a Sobolev estimate implies that there exists some constant $c_{1}>0$ so that 
\[\lrb{\int_{\{x-1<t<x\}}\abs{d\phi}^{p}dtd\theta}^{\frac 1p}\geq c_{1}\sup _{\{x-1<t_{i}<x\}}\abs{\phi(t_{1},\theta_{1})-\phi(t_{2},\theta_{2})}\ .\]
Therefore,
\[\sup _{\{x-1<t_{i}<x\}}\abs{\phi(t_{1},\theta_{1})-\phi(t_{2},\theta_{2})}\leq c_{1}e^{\delta x}\norm{e_{\bar S_{j}}\dvert \phi}_{\delta}\]
so
\[\sup _{\{t_{i}<x\}}\abs{\phi(t_{1},\theta_{1})-\phi(t_{2},\theta_{2})}\leq c_{1}\frac{e^{\delta x}}{1-e^{-\delta}}\norm{e_{\bar S_{j}}\dvert \phi}_{\delta}\ .\]
The above implies the required estimate:
\[\abs{\phi(t,\theta)-\phi(-\infty,\theta)}\leq c_{1}\frac{e^{\delta t}}{1-e^{-\delta}}\norm{e_{\bar S_{j}}\dvert \phi}_{\delta}\leq c'\norm{e_{\bar S_{j}}\dvert \phi}_{\delta}e_{\bar S_{j}}w_{0}^{\delta}\]

\stop

\

\begin{defn}Define the norm $\norm\nu_{\delta\&\delta'}$ by using the $\norm{\cdot}_{\delta'}$ norm on smooth manifold fibers of $\pi:\tilde U\longrightarrow U$ and the norm $\norm\cdot_{\delta+\delta'}$ on fibers with tropical part not equal to a point:

\[\norm\nu_{\delta\&\delta'}:=\max\begin{cases}\sup_{\totb{\pi^{-1}(x)}=\text{point}} \lrb{\int_{\pi^{-1}(x)\cap \tilde U}\abs{w_{0}^{-\delta'}\nu}^{p}}^{\frac 1p}\\ \sup_{\totb{\pi^{-1}(x)}\neq \text{point}} \lrb{\int_{\pi^{-1}(x)\cap \tilde U}\abs{w_{0}^{-\delta-\delta'}\nu}^{p}}^{\frac 1p}
\end{cases}\]

Define the norm $\norm\nu_{\delta\&\delta'}^{1}=\sup\abs\nu+\norm{\dvert\nu}_{\delta\&\delta'}$.

\end{defn}

We shall see that $\norm{\nu}_{\delta\&\delta'}$ is equivalent to the norm displayed in equation (\ref{equivalent norm2}) below.

\begin{lemma}\label{equivalent norm}

\begin{enumerate}
\item On any  allowable coordinate chart, the norm $\norm{\nu}_{0,\delta}$ is equivalent to the following norm using the lifts $\tilde{\mathcal S}$ of sets of strata $\mathcal S$ in $\totb{U}$:
\[\norm\nu_{\delta'}+\max_{ \mathcal S}\norm{w_{\tilde{\mathcal S}}^{-\delta}\Delta_{\tilde{\mathcal S}}\nu}_{\delta\&\delta'}\]

\item On any  allowable coordinate chart, the norm $\norm{\nu}^{1}_{0,\delta}$ is equivalent to the following norm:
\[\norm{\nu}^{1}_{\delta'}+\max_{ \mathcal S}\norm{w_{\tilde{\mathcal S}}^{-\delta}\Delta_{\tilde{\mathcal S}}\nu}^{1}_{\delta\&\delta'}\]
\end{enumerate}
\end{lemma}

\

\pf

The strata on which the higher weight $\delta+\delta'$ is used for the $\norm\cdot_{\delta\&\delta'}$ norm are the strata $T$ so that $T^{c}\neq\emptyset$. Therefore, the norm $\norm\nu_{\delta\&\delta'}$ is equivalent to the following norm
\begin{equation}\label{equivalent norm2}\norm\nu_{\delta'}+\max_{T^{c}\neq 0}\norm{e_{T}\nu}_{\delta+\delta'}\end{equation}
because  $\norm{e_{T}\nu}_{\delta+\delta'}$ is equivalent to taking the supremum involved in the $\norm{\cdot}_{\delta+\delta'}$ norm just over the stratum $T$. 

Therefore, to show that $\norm{\nu}_{0,\delta}$ and $\norm{\nu}^{1}_{0,\delta}$ dominate the two new norms above, it suffices to show that $\norm{\nu}_{0,\delta}$  and $\norm{\nu}^{1}_{0,\delta}$ dominate $\norm{e_{T}w_{\ti S}^{-\delta}\Delta_{\ti S}\nu }_{\delta'+\delta}$ and $\norm{e_{T}w_{\ti S}^{-\delta}\Delta_{\ti S}\nu }^{1}_{\delta'+\delta}$ respectively. 

\

{\em Claim}: if  $T^{c}\neq\emptyset$,  there exists a constant $c$ so that \begin{equation}\label{claim}e_{T}w_{\ti S\cup\mathcal S_{0}}\leq c e_{T}w_{\ti S}w_{0}\ .\end{equation} (Here $\ti S$ is a lifted set of strata and $\mathcal S_{0}$ indicates the set of  strata on which $w_{0}$ disappears.) 

\

The above claim holds  trivially  if $T\in\mathcal S_{0}$ or $e_{T}\Delta_{\ti S}=0$, as then both $e_{T}w_{\ti S\cup\mathcal S_{0}}$ and $e_{T}w_{\ti S}w_{0}$ are $0$. As $T^{c}\neq \emptyset$, assuming that $T\notin\mathcal S_{0}$ gives that  $e_{T}w_{0}=\abs{\zeta_{0}}$, where $\zeta_{0}$ is a smooth monomial either equal to $\totl{\tilde z_{1}}$ or $\totl{\tilde z^{\beta}}$. Suppose now that $\zeta$ is one of the smooth monomials used in the definition of $w_{\ti S\cup\mathcal S_{0}}$. For any smooth monomial, either $e_{S}\zeta=0$ or $e_{S}\zeta=\zeta$. Because $\zeta=0$ on $\mathcal S_{0}$, if $e_{T}\zeta\neq0$, then $\zeta$ must be  $\zeta_{0}^{n}\zeta'$ where $\zeta'$ is some lifted smooth monomial and $n\geq 1$. Moreover, $\zeta'$ must disappear on $\ti S$; this is because both $\ti S$ and  $\zeta'$ are lifted, so if $\zeta'$ is nonzero  on some stratum $S\in\ti S$, it is nonzero on all $S'\in S^{c}\subset \ti S$, but  $\zeta_{0}\neq 0$ on one such $S'$ and therefore $\zeta\neq 0$ on $S'$ --- a contradiction because $\zeta_{0}^{n}\zeta'=0$ on all $S'\in \ti S$.   As $\zeta'$ vanishes on all strata in $\ti S$, it is bounded by some constant times $w_{\ti S}$, so $\zeta=e_{T}\zeta$ is bounded by some constant times $e_{T}w_{\ti S}w_{0}$. Applied to all monomials  $\zeta$ appearing in $w_{\ti S\cup S_{0}}$, this proves the above claim  that $e_{T}w_{\ti S\cup\mathcal S_{0}}$ is bounded by a constant times $e_{T}w_{\ti S}w_{0}$.

\

Our inequality (\ref{claim}) implies that  $\norm{e_{T}w_{\ti S}^{-\delta}\Delta_{\ti S}\nu }_{\delta'+\delta}$ is dominated by $\norm{w^{-\delta}_{\ti S\cup\mathcal S_{0}}\Delta_{\ti S\cup\mathcal S_{0}}\nu}_{\delta'}$ and $\norm{e_{T}w_{\ti S}^{-\delta}\Delta_{\ti S}\nu }^{1}_{\delta'+\delta}$ is dominated by $\norm{w^{-\delta}_{\ti S\cup\mathcal S_{0}}\Delta_{\ti S\cup\mathcal S_{0}}\dvert \nu}_{\delta'}+\sup\abs{w^{-\delta}_{\ti S}\Delta_{\ti S} \nu}$ so our new norms are  dominated by $\norm\nu_{0,\delta}$ and $\norm\nu_{0,\delta}^{1}$ respectively.

\

We  must show that $\norm\nu_{0,\delta}$ and $\norm\nu_{0,\delta}^{1}$ are dominated by the new norms above.
For  $I$ an arbitrary collection of strata, we need to bound $\norm{w_{I}^{-\delta}\Delta_{I}\nu}_{\delta'}$ with these new norms. Let $\ti S$ be the largest lifted collection of strata so that $\Delta_{\ti S}\Delta_{I}=\Delta _{I}$. Then $\Delta_{I}=\Delta_{\ti S}\prod\Delta_{S}$, where the product is over strata $S \in I$ not contained in $\ti S$, (which implies that $S^{c}\neq\emptyset$). Simple computation gives the following identity:

\[\Delta_{S}=\Delta_{ S\cup S^{c}}+\sum_{\bar S\in S^{c}} e_{\bar S}\Delta_{S}\]  
Therefore,
\begin{equation}\label{Delta_{I}}\Delta_{I}=\sum_{I'\subset (I-\ti S)}\lrb{\prod_{S\in I-\ti S-I'}\sum_{\bar S\in S^{c}}e_{\bar S}\Delta_{S}}\Delta_{\ti S\cup I'\cup (I')^{c}}\ .\end{equation}

Note that if $\bar S\in S^{c}$, then $e_{\bar S}\Delta_{S}\Delta_{\mathcal S_{0}}=e_{\bar S}\Delta_{\mathcal S_{0}}$.
If $\norm\nu_{\delta'}$ is finite, $\Delta_{\mathcal S_{0}}\nu=\nu$, so if $\norm\nu_{\delta'}$ is finite we can replace $e_{\bar S}\Delta_{S}$ with $e_{\bar S}$ in the above equation (\ref{Delta_{I}}), getting
\begin{equation}\label{de2}\Delta_{I}\nu=\sum_{I'\subset (I-\ti S)}\lrb{\prod_{S\in I-\ti S-I'}\sum_{\bar S\in S^{c}}e_{\bar S}}\Delta_{\ti S\cup I'\cup (I')^{c}}\nu \ \ \ \ \text{ if }\Delta_{\ti S_{0}}\nu=\nu\ .\end{equation}

 To bound $\norm{w_{I}^{-\delta}\Delta_{I}\nu}_{\delta'}$ and $\norm{w_{I}^{-\delta}\Delta_{I}\nu}^{1}_{\delta'}$ using the above decomposition of $\Delta_{I}$, we also need the following estimate:
\begin{equation}\label{cw_{I}}cw_{I}\geq  \lrb{\prod_{S\in I-\ti S-I'}e_{\bar S} }w_{\ti S\cup I' \cup(I')^{c} }w_{0}\end{equation}
 To prove the above inequality, suppose that $\zeta$ is one of the smooth monomials involved in the  righthand side, so $\zeta$ vanishes on all strata in $ \mathcal S_{0}\cup\ti S\cup I' \cup (I')^{c}$, but does not vanish on $\bar S$. In particular, $\zeta=\totl{\tilde \zeta}$ where  $\totb{\tilde \zeta}= 0$ on some $\bar S\in S^{c}$, and $\totb{\tilde \zeta}>0$ on $S^{0}\in\mathcal S_{0}$, where $\totb\pi (S)=\totb\pi (S^{0})$. Therefore, $\totb{\tilde \zeta}> 0$ on $S$. This is valid for all $S\in I-\ti S-I'$. As all other strata in $I$ appear in $\ti S\cup I' \cup (I')^{c}$, it follows that $\zeta$ vanishes on $I$, so is bounded by a constant times $  w_{I}$. The above inequality (\ref{cw_{I}}) follows. 
 
 Using the decomposition (\ref{de2}), the above inequality (\ref{cw_{I}}) and the fact that $w_{I}\geq c w_{I\cup I^{c}}$, we get that the following inequality holds, so long as $\Delta_{\mathcal S_{0}}\nu=\nu$.

\[
\begin{split}\norm{w_{I}^{-\delta}\Delta_{I}\nu}_{\delta'}&\leq c\norm{w^{-\delta}_{I\cup I^{c}}\Delta_{I\cup I^{c}}\nu}_{\delta'}
\\& \ \ \ \ \ +c\sum_{I'\subsetneq(I-\ti S)}\norm{\lrb{\prod_{S\in I-\ti S- I'}\sum_{\bar S\in S^{c}}e_{\bar S}}w^{-\delta}_{\ti S\cup I'\cup (I')^{c}}w^{-\delta}_{0}\Delta_{\ti S\cup I'\cup (I')^{c}}\nu}_{\delta'}
\\&= c\norm{w^{-\delta}_{I\cup I^{c}}\Delta_{I\cup I^{c}}\nu}_{\delta'}
\\& \ \ \ \ \ +c
\sum_{I'\subsetneq(I-\ti S)}\norm{\lrb{\prod_{S\in I-\ti S-I'}\sum_{\bar S\in S^{c}}e_{\bar S}}w^{-\delta}_{\ti S\cup I'\cup (I')^{c}}\Delta_{\ti S\cup I'\cup (I')^{c}}\nu}_{\delta'+\delta}
\end{split}\]
 
 (In the above inequality, we also used the fact that $e_{S}w_{I}\leq w_{I}$.) All the terms above are dominated by  the first new norm, and we have proved that $\norm\nu_{0,\delta}$ is dominated by the first new norm.

\

To dominate $\norm{w^{-\delta}_{I}\Delta\dvert \nu}_{\delta'}$ by the second new norm, we can use the above inequality. Now dominate $\sup\abs{w^{-\delta}_{I}\Delta\nu}$ using (\ref{Delta_{I}}) and (\ref{cw_{I}}) as follows: 
\begin{equation}\label{boundabs}
\begin{split}\sup &\abs{w_{I}^{-\delta}\Delta_{I}\nu}\leq c\sup\abs{w^{-\delta}_{I\cup I^{c}}\Delta_{I\cup I^{c}}\nu}
\\&  +c\sum_{I'\subsetneq(I-\ti S)}\sup\abs{\sum \lrb{\prod_{S\in I-\ti S- I' }e_{\bar S}}
w^{-\delta}_{\ti S\cup I'\cup(I')^{c}}w_{0}^{-\delta}\Delta_{I-\ti S- I'}\Delta_{\ti S\cup I'\cup(I')^{c}}\nu}
\end{split}\end{equation}
The unspecified sum (above  and in the equation below) is over all choices of $\bar S\in S^{c}$ for all $S\in I-\ti S- I'$.
We must bound the terms appearing in the sum above in (\ref{boundabs}). Lemma \ref{sobolev} gives that that $\abs{e_{\bar S}\Delta_{S}\phi}\leq c\norm{e_{\bar S}\dvert \phi}_{\delta+\delta'}e_{\bar S}w_{0}^{\delta+\delta'}$. Use this inequality for a single $S\in I-\ti S-I'$, and remove the other terms in $\Delta_{I-\ti S-I'}$ by observing that $\norm{\Delta_{S'}\phi}\leq 2\norm{\phi}$. 
\[\begin{split}\sup\abs{\sum \lrb{\prod_{S\in I-\ti S- I' }e_{\bar S}}
w^{-\delta}_{\ti S\cup I'\cup(I')^{c}}w_{0}^{-\delta}\Delta_{I-\ti S- I'}\Delta_{\ti S\cup I'\cup(I')^{c}}\nu}\hspace{-7cm}&
\\&\leq c\norm{\lrb{\prod_{S\in I-\ti S- I'}\sum_{\bar S\in S^{c}}e_{\bar S}}w^{-\delta}_{\ti S\cup I'\cup(I')^{c}}\Delta_{\ti S\cup I'\cup(I')^{c}}\dvert \nu}_{\delta+\delta'}
\end{split}\]

The last term is bounded by our second new norm, so  this completes the proof that that $\norm{\nu}_{0,\delta}^{1}$ is dominated by the second new norm.

\stop

\subsection{Norms on an allowable collection of coordinate charts}

\

\

We are now ready to define our norms on a vectorbundle $V$ over a family of curves $\pi_{\ex F}:\hat{\ex C}\longrightarrow\ex F$. To avoid the problems that arise if  $\ex F$ is not  compact, we  define the notion of an `allowable' family, which is extendable and can be covered by an `allowable' collection of extendable allowable coordinate charts. We shall also need a version of this definition when there is a collection of marked points on our family.

\begin{defn}\label{allowable collection}
An allowable collection of coordinate charts on a vectorbundle $V$ over a family $\pi_{\ex F}:\hat{\ex C}\longrightarrow \ex F$ is a local trivialization of $V$ over a finite collection of extendable allowable coordinate charts $\pi:U_{\alpha,i}\longrightarrow U_{\alpha}$ satisfying the following additional conditions:

 \begin{enumerate}

\item \[\pi_{\ex F}^{-1}\lrb{U_{\alpha}}=\bigcup_{i} U_{\alpha, i}\]
\item The restriction of the vectorbundle $V$ to $U_{\alpha,i}$ is $\mathbb R^{n}\times U_{\alpha,i}$ with the obvious projection. 
\item Coordinate change maps and intersections between $U_{\alpha, i}$ and $U_{\alpha,j}$ satisfy the following:
\begin{enumerate}
\item If $U_{\alpha, i}$ and $U_{\alpha,j}$ are coordinate charts of type \ref{interesting} from the Definition  \ref{allowable} (these are the charts that cover  an internal edge of an exploded curve), and $i\neq j$, then they do not intersect.
\item If $U_{\alpha,i}$ and $U_{\alpha,j}$ are coordinate charts of type \ref{boring}, they are the product of  an  open subset  of $\et 11$ with $U_{\alpha}$. In this case, their  intersection in either of these coordinate charts is also  the product of  an open subset of ${\et 11}$ with $U_{\alpha}$, and the  coordinate change map is a product of some map between these subsets of $\et 11$ with the identity map on $U_{\alpha}$. The coordinate change map between the vector bundle trivializations over these charts is independent of position in $U_{\alpha}$ and depends only on $\et 11$.
\item If $U_{\alpha,i}$ is a coordinate chart of type \ref{interesting} and $U_{\alpha,j}$ is a coordinate chart of type \ref{boring}, their intersection is as follows: In $U_{\alpha,j}$, it is  the product of a subset  $O\subset \et 11$ with $U_{\alpha}$. If $U_{\alpha,i}$ is given by $\{\abs {\tilde z_{1}}<c,\abs {\tilde z^{\beta}}<c\}$, then the intersection with $U_{\alpha}$ is  a subset of the form $\{\tilde z_{1}\in O'\}$ or $\{\tilde z^{\beta}\in O'\}$ where $O'\subset \{\frac c2<\abs z<c\}\in\mathbb C$. In either case, we identify this subset with the product $O'\times U_{\alpha}$. The transition map in this case is  the product of a diffeomorphism between $O$ and $O'$ and the identity on $U_{\alpha}$. Again, the coordinate change map between the vector bundle trivializations over these charts is independent of position in $U_{\alpha}$ and depends only on $O$.
\end{enumerate}

\end{enumerate}

If our family also has a collection of non-intersecting marked-points  sections $\ex F\longrightarrow \hat {\ex C}$  not intersecting edges of curves in our family, then an  allowable collection of  coordinate chart is an allowable collection as above with the extra conditions that no marked points are inside coordinate charts of type \ref{interesting}, and in the  coordinate charts of type \ref{boring} containing marked points, the sections corresponding to these marked points are constant sections $U_{\alpha}\longrightarrow \et 11\times U_{\alpha}$. 
\end{defn}

Outside a neighborhood of internal edges, every family of exploded curves is locally trivial (ignoring complex structures), so it is easy to verify that any single curve $\ex C\neq \ex T$ in a family of curves (with marked points as described above) is covered by an allowable collection of coordinate charts --- first cover internal edges with separated charts $U_{i}\longrightarrow U$ of type \ref{interesting} that avoid marked points. On the subset of $U_{i}$ where $c/2<\tilde z_{1}<c$ or $c/2<\tilde z^{\beta}<c$ there is a canonical map to  $U_{i}\cap \ex C$ given by keeping all coordinates apart from $\tilde z^{\beta}$ or $\tilde z_{1}$ constant; after possibly  shrinking $c$, there is a contractible choice of extension of this map to $\pi_{\ex F}^{-1}U\setminus\bigcup _{i} U_{i}\longrightarrow \ex C$ trivializing the family $\pi_{\ex F}^{-1}(U)\longrightarrow U$ (along with its marked points). We can then choose the remaning charts and transition maps compatible with this trivilization.   Moreover, by shrinking and subdividing these coordinate charts if necessary, we may extend these allowable coordinate charts to any  finite-rank vector bundle. Note that the almost complex structure on fibers can not be assumed to be  the standard one from the allowable coordinate charts.

\begin{defn}\label{allowable family}
An allowable family $\hat{\ex C}\xrightarrow {\pi_{\ex F}}\ex F$  is a subset of a family of curves $\hat{\ex C}\xrightarrow {\pi_{\ex F}}\ex F'$ covered, along with the vectorbundle $X$, by an allowable collection of coordinate charts. 

In particular, an allowable family of curves 
\[\begin{tikzcd}\hat {\ex C}\rar{\hat f}\dar &\hat{\ex B}\dar
\\ \ex F \rar& \ex B_0
\end{tikzcd}\]
 is a family so that $\hat {\ex C}\longrightarrow \ex F$ with the vectorbundle $\hat f^{*}T_{vert}\hat {\ex B}$ is covered by an allowable collection of coordinate charts.
\end{defn}

If $\hat{\ex C}\longrightarrow \ex F$ is an allowable family of curves, $T^{*}_{vert}\hat {\ex C}$ is canonically an allowable vectorbundle, trivialized  using the standard basis for the cotangent space in coordinates consisting of $dx_{i}$ and the real and imaginary parts of $\tilde z^{-1}_{i}d\tilde z_{i}$. So, $T^{*}_{vert}\hat {\ex C}\otimes X$ or $T^{*}_{vert}\hat{\ex C}\otimes \hat f^{*}T_{vert}\hat {\ex B}$ is also canonically an allowable vectorbundle. For the purposes of measuring the norm of sections of $Y:=(T^{*}_{vert}\hat{\ex C}\otimes X)^{(0,1)}$, we consider $Y$ as a sub-bundle of this allowable vectorbundle. Note that $Y$ may not be a constant sub-bundle using this trivialization as  $j$ may not be constant in our local trivialization of $T_{vert}^{*}\hat{\ex C}$. 

\begin{defn}\label{mixed definition}
On an allowable family $\hat f$ with vectorbundle $V$ covered by the allowable collection of  coordinate charts $U_{\alpha,i}\longrightarrow U_{\alpha}$, define the norms

\[\norm\nu_{**}:=\sum _{U_{\alpha,i}}\norm{\nu\rvert_{U_{\alpha,i}}}_{**}\]
where $**$ stands for the different possible labels for norms defined so far. 
For example, $\norm\nu^{1}_{k,\delta}:=\sum _{U_{\alpha,i}}\norm{\nu\rvert_{U_{\alpha,i}}}^{1}_{k,\delta}$. 

\

Given a choice of weight $\delta_{\alpha,i}$ for each coordinate chart $U_{\alpha,i}$, define the following norms:

\[\norm\nu_{\text{mixed } \delta}:=\sum _{U_{\alpha,i}}\norm{\nu\rvert_{U_{\alpha,i}}}_{\delta_{\alpha,i}}\]

\[\norm\nu^{1}_{\text{mixed } \delta}:=\sum _{U_{\alpha,i}}\norm{\nu\rvert_{U_{\alpha,i}}}^{1}_{\delta_{\alpha,i}}\]

\end{defn}

One problem with the norms $\norm\cdot_{k,\delta}$ is that $w_{\mathcal S}\Delta_{\mathcal S}\nu$ is only defined on coordinate patches and not globally defined on a fiber. The following definition provides a means of remedying this.

\

\begin{defn}\label{tilde delta}
Given an allowable collection of coordinate charts $U_{\alpha,i}\longrightarrow U_{\alpha}$, and a collection $\mathcal S$ of strata in $\totb{U_{\alpha}}$, define $\tilde \Delta_{\mathcal S}\nu$ on $\bigcup_{i}U_{\alpha,i}$ as follows:

\begin{enumerate}
\item Choose a smooth cutoff function $\rho: \et 11\longrightarrow [0,1]$ so that 
\[\rho(\tilde z)=\begin{cases}1\text{ if }\abs{\tilde z}\leq \frac 12 
\\ 0\text{ if }\abs{\tilde z}\geq 1\end{cases}\]
\item On coordinate charts $U_{\alpha,i}$ of type \ref{interesting}, we have coordinates $\tilde z_{1}$ and $\tilde z_{\beta}$ so that $\abs{\tilde z_{1}}<c$,  $\abs{\tilde z_{\beta}}<c$ and $\abs{\tilde z_{1}\tilde z^{\beta}}<\frac {c^{2}}{16}$. Use the notation $\tilde {\mathcal S}$ for the lift of $\mathcal S$ to this chart, $\mathcal S^{+}$ for the collection of strata in $\tilde{\mathcal S}$ so that $\totb{\tilde z_{1}}=0$, and $\mathcal S^{-}$ for the collection of strata in $\tilde {\mathcal S}$ so that $\totb{\tilde z^{\beta}}=0$.  

Define $\tilde \Delta_{\mathcal S}\nu$ on this chart by:
\[\begin{split}\tilde \Delta_{\mathcal S}\nu:=&\rho\lrb{\frac{2\tilde z_{1}}c}\rho\lrb{\frac{2\tilde z^{\beta}}c}\Delta_{\tilde{\mathcal S}}\nu
\\&+
\lrb{1-\rho\lrb{\frac{2\tilde z_{1}}c}} \Delta_{\mathcal S^{+}}\nu
\\&+\lrb{1-\rho\lrb{\frac{2\tilde z^{\beta}}c}} \Delta_{\mathcal S^{-}}\nu\end{split}\] 

\item On all other  coordinate charts $U_{\alpha_{i}}$ of type \ref{boring}, let $\tilde {\mathcal S}$ denote the lift of $\mathcal S$  to this coordinate chart, and define
\[\tilde \Delta_{\mathcal S}\nu:=\Delta_{\tilde{\mathcal S}}\nu \ . \]
\end{enumerate}

\end{defn}

\

It follows from the types of transition functions allowed for allowable collections of coordinate charts that $\tilde \Delta_{\mathcal S} \nu$ is well defined on $\bigcup_{i}U_{\alpha,i}$ and smooth or $\C\infty1$ if $\nu$ is. We can now state a fiberwise-global version of Lemma \ref{equivalent norm}.

\

\begin{lemma}\label{global norm}
Restricted to a collection of allowable coordinate charts $U_{\alpha,i}$ over a single chart $U_{\alpha}$, the norm $\norm\nu_{0,\delta}$  is equivalent to the norm

\[\norm\nu_{\delta'}+\max_{\mathcal S}\norm{w^{-\delta}_{\mathcal S}\tilde \Delta_{\mathcal S} \nu}_{\delta\&\delta'}\]

and the norm $\norm\nu_{0,\delta}^{1}$ is equivalent to the norm

\[\norm\nu^{1}_{\delta'}+\max_{\mathcal S}\norm{w^{-\delta}_{\mathcal S}\tilde \Delta_{\mathcal S} \nu}^{1}_{\delta\&\delta'} \ .\]

In the above, the maximum is taken over all collections of strata $\mathcal S$ in $\totb{U_{\alpha}}$, and $w_{S}$ indicates the lift of the weighting function  $w_{\mathcal S}$ on $U_{\alpha}$, (which  equals  $w_{\tilde {\mathcal S}}$ in each of the coordinate charts $U_{\alpha, i}$). The norm $\norm\cdot_{\delta\&\delta'}$ is defined on page \pageref{equivalent norm}.

\end{lemma}

\pf

On all coordinate charts of type \ref{boring}, this lemma follows immediately from Lemma \ref{equivalent norm} as in this case $\tilde\Delta_{\mathcal S}=\Delta_{\tilde{\mathcal S}}$.

On a coordinate chart of type \ref{interesting}, where $\abs{\tilde z_{1}}<c$ and $\abs{\tilde z^{\beta}}<c$, using the notation from Definition \ref{tilde delta} we have that 
\[\begin{split}\Delta_{\tilde {\mathcal S}}\nu-\tilde \Delta_{\mathcal S}\nu =&\lrb{1-\rho\lrb{\frac{2\tilde z_{1}}c}\rho\lrb{\frac{2\tilde z^{\beta}}c}}\Delta_{\tilde{\mathcal S}}\nu
\\&-
\lrb{1-\rho\lrb{\frac{2\tilde z_{1}}c}} \Delta_{\mathcal S^{+}}\nu
\\&-\lrb{1-\rho\lrb{\frac{2\tilde z^{\beta}}c}} \Delta_{\mathcal S^{-}}\nu\ .\end{split}\] 

Use the notation $\check\rho^{+}$ for the function $(1-\rho(\frac {2\tilde z_{1}}c) )$ and $\check\rho^{-}$ for $(1-\rho(\frac {2\tilde z^{\beta}}c) )$. As $\rho(\tilde z)=1$ when $\abs{\tilde z}\leq \frac 12$, and our coordinate chart has $\abs{\tilde z_{1}\tilde z^{\beta}}<c^{2}/{16}$, the region where $\rho(\frac {2\tilde z_{1}}c)\neq 1$ is disjoint from the region where $\rho(\frac {2\tilde z^{\beta}}c)\neq 1$, so we can rewrite the above equation as follows.

\begin{equation}\label{move tilde}\begin{split}\Delta_{\tilde {\mathcal S}}\nu-\tilde \Delta_{\mathcal S}\nu =&
\check\rho^{+}\lrb{\Delta_{\tilde{\mathcal S}}\nu- \Delta_{\mathcal S^{+}}\nu}
\\&+\check\rho^{-}\lrb{\Delta_{\tilde{\mathcal S}}\nu- \Delta_{\mathcal S^{-}}\nu}\end{split}\end{equation}

Given any stratum $T$ for which $T^{c}\neq \emptyset$, $e_{T}\totl{\tilde z_{1}\tilde z^{\beta}}=0$. We then obtain the following dichotomy.
\[\text{ If }e_{T}\totl{\tilde z_{1}}=0,\text{ then }e_{T}\check\rho^{+}=0\text{ and }e_{T}\Delta_{\tilde{\mathcal S}}=e_{T}\Delta_{\mathcal S^{-}},\] 
\[\text{ and if }e_{T}\totl{\tilde z^{\beta}}=0,\text{ then }e_{T}\check\rho^{-}=0\text{ and }e_{T}\Delta_{\tilde{\mathcal S}}=e_{T}\Delta_{\mathcal S^{+}}.\]
 Therefore, for any stratum $T$ so that $T^{c}\neq \emptyset$, in our coordinate chart
\begin{equation}\label{interior equivalence}e_{T}\tilde \Delta_{\mathcal S}=e_{T}\Delta_{\tilde{\mathcal S}}\ .\end{equation} 
The above equation (\ref{interior equivalence}) holds on any fiber of our coordinate chart $U_{\alpha, i}\longrightarrow U_{\alpha}$ with  nontrivial  tropical part, so (\ref{interior equivalence}) implies that 
\begin{equation}\label{ed}\norm{w^{-\delta}_{\mathcal S}\tilde \Delta_{\mathcal S}\nu-w^{-\delta}_{\mathcal S}\Delta_{\tilde {\mathcal S}}\nu}_{\delta\&\delta'}=\norm{w^{-\delta}_{\mathcal S}\tilde \Delta_{\mathcal S}\nu-w^{-\delta}_{\mathcal S}\Delta_{\tilde {\mathcal S}}\nu}_{\delta'}\end{equation}
and 
\begin{equation}\label{e1d}\norm{w^{-\delta}_{\mathcal S}\tilde \Delta_{\mathcal S}\nu-w^{-\delta}_{\mathcal S}\Delta_{\tilde {\mathcal S}}\nu}^{1}_{\delta\&\delta'}=\norm{w^{-\delta}_{\mathcal S}\tilde \Delta_{\mathcal S}\nu-w^{-\delta}_{\mathcal S}\Delta_{\tilde {\mathcal S}}\nu}^{1}_{\delta'}\ .\end{equation}

Note that $\Delta_{\ti S}-\Delta_{\mathcal S^{+}}$ is a sum of terms $\pm e_{T}\Delta_{\mathcal S^{+}}$ for $T\in \ti S- \mathcal S^{+}$,  so equation (\ref{move tilde}) implies that there exists some constant $c'$ so that on this coordinate chart,
\[\begin{split}\norm{w^{-\delta}_{\mathcal S}\tilde \Delta_{\mathcal S}\nu-w^{-\delta}_{\mathcal S}\Delta_{\tilde {\mathcal S}}\nu}_{\delta'}
\leq& c'\sum_{T\in \tilde {\mathcal S}-\mathcal S^{+}}\norm{w^{-\delta}_{\mathcal S}\check\rho^{+}e_{T}\Delta_{\mathcal S^{+}}\nu}_{\delta'}
\\&+c'\sum_{T\in \tilde {\mathcal S}-\mathcal S^{-}}\norm{w^{-\delta}_{\mathcal S}\check \rho^{-}e_{T}\Delta_{\mathcal S^{-}}\nu}_{\delta'}\ .\end{split}\]

Use the notation  $\hat{\mathcal S}$ to indicate $\mathcal S^{+}\cap \mathcal S^{-}$. This is the set of strata within $\ti S$ where $\totb{\tilde z_{1}}=0=\totb{\tilde z^{\beta}}$, or equivalently the set of strata $S\in\ti S$ with $S^{c}=\emptyset$. If $T\in\tilde{\mathcal S}-\mathcal S^{+}$, then \[e_{T}\Delta_{\mathcal S^{+}}=e_{T}\Delta_{\hat{\mathcal S}}\Delta_{\mathcal S_{0}}\] where $\mathcal S_{0}$ is the set of strata where $\totl{\tilde z_{1}}=0=\totl{\tilde z^{\beta}}$. To see this, note that if $T'\in\mathcal S^{+}$, then either $T'\in \hat{\mathcal S}$ or  $e_{T}e_{T'}\Delta_{\mathcal S_{0}}=0$ because $e_{T}\totl{\tilde z_{1}}=0$ and $e_{T'}\totl{\tilde z^{\beta}}=0$. Therefore, $e_{T}\Delta_{\mathcal S^{+}}e_{T}\Delta_{\hat{\mathcal S}}\Delta_{\mathcal S_{0}}=e_{T}\Delta_{\hat{\mathcal S}}\Delta_{\mathcal S_{0}}$. Similarly, if $S\in \mathcal S_{0}$, then $e_{S}e_{T}\Delta_{\mathcal S^{+}}=0$, because if $T^{+}$ indicates the strata in $\mathcal S^{+}$ with the same projection as $T$, then $e_{S}e_{T^{+}}=e_{S}e_{T}$. Therefore, $e_{T}\Delta_{\mathcal S^{+}}e_{T}\Delta_{\hat{\mathcal S}}\Delta_{\mathcal S_{0}}=e_{T}\Delta_{\mathcal S^{+}}$, so $e_{T}\Delta_{\mathcal S^{+}}=e_{T}\Delta_{\hat{\mathcal S}}\Delta_{\mathcal S_{0}}$ as required.

 Similarly, if $T\in\tilde{\mathcal S}-\mathcal S^{-}$, then $e_{T}\Delta_{{\mathcal S^{-}}}=e_{T}\Delta_{\hat{\mathcal S}}\Delta_{\mathcal S_{0}}$. Therefore, we get

\[\begin{split}\norm{w^{-\delta}_{\mathcal S}\tilde \Delta_{\mathcal S}\nu-w^{-\delta}_{\mathcal S}\Delta_{\tilde {\mathcal S}}\nu}_{\delta'}
\leq& c'\sum_{T\in \tilde {\mathcal S}-\mathcal S^{+}}\norm{w^{-\delta}_{\mathcal S}\check\rho^{+}e_{T}\Delta_{\hat{\mathcal S}\cup \mathcal S_{0}}\nu}_{\delta'}
\\&+c'\sum_{T\in \tilde {\mathcal S}-\mathcal S^{-}}\norm{w^{-\delta}_{\mathcal S}\check \rho^{-}e_{T}\Delta_{\hat{\mathcal S}\cup\mathcal S_{0}}\nu}_{\delta'}\ .\end{split}\]

\emph{Claim}: if $T\in \tilde {\mathcal S}-\mathcal S^{+}$, then $w_{\mathcal S}^{-\delta}\check \rho^{+}$ is bounded by a constant times  $e_{T}w^{-\delta}_{\hat{\mathcal S}}w_{0}^{-\delta}$. To prove this claim, it suffices to show that $\check\rho^{+}e_{T}w_{\hat{\mathcal S}}w_{0}$ vanishes on all the strata in $\tilde{\mathcal S}$. If $S$ is a stratum in $\ti S$, then either $S\in\hat{\mathcal S}$,  $S\in\mathcal S_{0}$,  $e_{S}\totl{\tilde z_{1}}=0$, or $e_{S}\totl{\tilde z^{\beta}}=0$. In the first two cases, $e_{S}w_{\hat{\mathcal S}}w_{0}=0$. If $e_{S}\totl{\tilde z_{1}}=0$, then $e_{S}\check\rho^{+}=0$. If $e_{S}\totl{\tilde z^{\beta}}=0$, then as $e_{T}\totl{\tilde z_{1}}=0$, we get $e_{S}e_{T}w_{0}=0$. Therefore $e_{S}\check\rho^{+}e_{T}w_{\hat{\mathcal S}}w_{0}=0$, and the above claim follows. 

 Similarly, if $T\in\tilde {\mathcal S}-\mathcal S^{-}$, then $w_{\mathcal S}^{-\delta}\check \rho^{-}$ is bounded by a constant times  $e_{T}w^{-\delta}_{\hat{\mathcal S}}w_{0}^{-\delta}$. Therefore,

\[\begin{split}\norm{w^{-\delta}_{\mathcal S}\tilde \Delta_{\mathcal S}\nu-w^{-\delta}_{\mathcal S}\Delta_{\tilde {\mathcal S}}\nu}_{\delta'}
\leq& c''\sum_{T\in \tilde {\mathcal S}-\hat {\mathcal S}}\norm{e_{T}w^{-\delta}_{\hat{\mathcal S}}w_{0}^{-\delta}\Delta_{\hat{\mathcal S}\cup \mathcal S_{0}}\nu}_{\delta'}
\ .\end{split}\]
As all the strata $T$ in the above inequality satisfy $T^{c}\neq \emptyset$, on the strata $T$ the norm $\norm{\cdot}_{\delta\&\delta'}$ always uses the higher weight $w_{0}^{-\delta-\delta'}$. Therefore, we get that
\begin{equation}\label{ineq1}\norm{w^{-\delta}_{\mathcal S}\tilde \Delta_{\mathcal S}\nu-w^{-\delta}_{\mathcal S}\Delta_{\tilde {\mathcal S}}\nu}_{\delta'}
\leq c'''\norm{w^{-\delta}_{\hat{\mathcal S}}\Delta_{\hat{\mathcal S}}\nu}_{\delta\&\delta'}\ .
\end{equation}

The set of strata $\hat {\mathcal S}$ is the lift of some set of strata which we shall again call $\hat {\mathcal S}$, because each stratum lifts uniquely. With this slight abuse of notation $\hat{\mathcal S}^{+}=\hat{\mathcal S}=\hat{\mathcal S}^{-}$, so $\tilde\Delta_{\hat {\mathcal S}}=\Delta_{\hat{\mathcal S}}$, and we get

\begin{equation}\label{ineq2}\norm{w^{-\delta}_{\mathcal S}\tilde \Delta_{\mathcal S}\nu-w^{-\delta}_{\mathcal S}\Delta_{\tilde {\mathcal S}}\nu}_{\delta'}
\leq c'''\norm{w^{-\delta}_{\hat{\mathcal S}}\tilde\Delta_{\hat{\mathcal S}}\nu}_{\delta\&\delta'}\ .
\end{equation}

The above inequalities (\ref{ineq1}) and (\ref{ineq2})   together with equation (\ref{ed}) and  Lemma \ref{equivalent norm} proves that the norm $\norm\nu_{0,\delta}$  is equivalent to the norm

\[\norm\nu_{\delta'}+\max_{\mathcal S}\norm{w^{-\delta}_{\mathcal S}\tilde \Delta_{\mathcal S} \nu}_{\delta\&\delta'}\ .\]

As we already have the required estimates for the part of the $\norm\nu ^{1}_{\delta}$ norm involving $\norm{\dvert\nu}_{\delta}$, it remains to estimate 
\[\sup\abs{w^{-\delta}_{\mathcal S}\tilde \Delta_{\mathcal S}\nu-w^{-\delta}_{\mathcal S}\Delta_{\tilde {\mathcal S}}\nu}\]
on our coordinate chart. As argued for the $\norm\cdot_{\delta'}$ norm above, we get 
\[\sup\abs{w^{-\delta}_{\mathcal S}\tilde \Delta_{\mathcal S}\nu-w^{-\delta}_{\mathcal S}\Delta_{\tilde {\mathcal S}}\nu}\leq c'\sum_{T\in \tilde {\mathcal S}-\hat {\mathcal S}}\sup\abs{e_{T}w^{-\delta}_{\hat{\mathcal S}}w_{0}^{-\delta}\Delta_{\hat{\mathcal S}\cup \mathcal S_{0}}\nu}\ .\]
We can estimate the  righthand side of the above inequality with Lemma \ref{sobolev}, as $S_{0}$ will contain some stratum whose complement contains $T$. Therefore, 
\[\sup\abs{w^{-\delta}_{\mathcal S}\tilde \Delta_{\mathcal S}\nu-w^{-\delta}_{\mathcal S}\Delta_{\tilde {\mathcal S}}\nu}\leq c''\sum_{T\in \tilde {\mathcal S}-\hat {\mathcal S}}\norm{e_{T}w^{-\delta}_{\hat{\mathcal S}}\Delta_{\hat{\mathcal S}}\dvert\nu}_{\delta}\ .\]

So, as argued above, we get the two inequalities 
\begin{equation}\label{ineq3}\sup\abs{w^{-\delta}_{\mathcal S}\tilde \Delta_{\mathcal S}\nu-w^{-\delta}_{\mathcal S}\Delta_{\tilde {\mathcal S}}\nu}\leq c\norm{w^{-\delta}_{\hat{\mathcal S}}\tilde\Delta_{\hat{\mathcal S}}\dvert\nu}_{\delta\&\delta'}
\end{equation}
and 
\begin{equation}\label{ineq4}\sup\abs{w^{-\delta}_{\mathcal S}\tilde \Delta_{\mathcal S}\nu-w^{-\delta}_{\mathcal S}\Delta_{\tilde {\mathcal S}}\nu}\leq c\norm{w^{-\delta}_{\hat{\mathcal S}}\Delta_{\hat{\mathcal S}}\dvert\nu}_{\delta\&\delta'}\ .
\end{equation}

This together with Lemma \ref{equivalent norm} completes the proof of our lemma.

\stop

\

The above proof contains the inequalities  (\ref{ineq1}), (\ref{ineq2}), (\ref{ineq3}) and (\ref{ineq4}), which together with equations (\ref{ed}) and (\ref{e1d}) imply the following estimates which will be useful later:

\begin{lemma}\label{tilde estimates}
Let $\hat{\mathcal S}$ denote the subset of $\ti S$ consisting of all strata $T\in\ti S$ so that $T^{c}=\emptyset$, and also let $\hat {\mathcal S}$ denote the corresponding subset of $\mathcal S$. Then there exist some  constant $c$ so that the following inequalities hold for all $\nu$ so that both sides are defined.
\[\norm{w^{-\delta}_{\mathcal S}\tilde \Delta_{\mathcal S}\nu-w^{-\delta}_{\ti S}\Delta_{\ti S}\nu}_{\delta\&\delta'}\leq
c\norm {w^{-\delta}_{\hat{\mathcal S}}\Delta_{\hat {\mathcal S}}\nu }_{\delta\&\delta'}=
c\norm {w^{-\delta}_{\hat{\mathcal S}}\tilde\Delta_{\hat {\mathcal S}}\nu }_{\delta\&\delta'}\] 

\[\norm{w^{-\delta}_{\mathcal S}\tilde \Delta_{\mathcal S}\nu-w^{-\delta}_{\ti S}\Delta_{\ti S}\nu}^{1}_{\delta\&\delta'}\leq
c\norm {w^{-\delta}_{\hat{\mathcal S}}\Delta_{\hat {\mathcal S}}\nu }^{1}_{\delta\&\delta'}=
c\norm {w^{-\delta}_{\hat{\mathcal S}}\tilde \Delta_{\hat {\mathcal S}}\nu }^{1}_{\delta\&\delta'}\] 
\end{lemma}

\section{Analysis of $\dbar$ equation in families}\label{analysis}

\subsection{$\dbar$ is continuously differentiable} 

\

Throughout this section, we shall be considering an operator $\dbar$  in the form (\ref{dbar form}) described on page \pageref{dbar form}.

\

Recall that an allowable family of curves as defined on page \pageref{allowable family} is an extendable family covered by an allowable collection of coordinate charts in which the vectorbundle $X$ or $\hat f^{*}T_{vert}\hat{\ex B}$ is trivialized. So, on the space of sections of $X$, we may use our norms from Definition \ref{mixed definition}. Similarly, each of our allowable coordinate charts comes with a canonical trivialization of $T_{vert}^{*}\hat{\ex C}$, and therefore a trivialization of $T_{vert}^{*}\hat{\ex C}\otimes X$, so we may define norms of sections of $T_{vert}^{*}\hat{\ex C}\otimes X$ and therefore norms on sections of the sub bundle $Y:=(T_{vert}^{*}\hat{\ex C}\otimes X)^{(0,1)}$.  Note that any family of curves may be shrunk to a compactly contained sub-family to obtain an allowable family of curves.

\begin{defn}\label{defX} Consider an allowable  family of curves with some (possibly empty)  collection   of nonintersecting marked-point  sections $\ex F\longrightarrow\hat{\ex C}$ avoiding the edges of the curves in $\hat{\ex C}\longrightarrow \ex F$. 

 Define the Banach space $X_{**}$ to be the $\norm\cdot^{1}_{**}$--completion of  the subspace of $X^{\infty,\underline 1}$  with finite $\norm\cdot^1_{**}$ norm,\footnote{A section in $X^{\infty,\underline 1}$ will have finite $\norm{\cdot}^1_{**}$--norm if it is extendable, but not in general because we don't have control near the boundary of $\ex F$.} where $X^{\infty,\underline 1}$ is the space of  $\C\infty1$ sections of the vector bundle $X$ vanishing at marked points, and  $**$ stands for the different labels for norms used in section \ref{norms}.
 
 Define $Y_{**}$ to be the $\norm\cdot_{**}$--completion of the subspace of $Y^{\infty,\underline 1}$ with finite norm, where $Y^{\infty,\underline 1}$ is  the space of   $\C\infty1$ sections of $Y$ vanishing on all edges of curves in $\hat{\ex C}\longrightarrow \ex F$.

Given a curve $f$ in $\hat{\ex C}\longrightarrow\ex F$,  we can restrict our  data to the domain $\ex C$ of $f$. Define $X_{\delta}(f)$ to be the corresponding Banach space $X_{\delta}$ with this restricted data.
If $\nu\in X_{\delta}$, define $\nu(f)\in X_{\delta}(f)$ to be the restriction of $\nu$ to the domain $\ex C$ of $f$. 
We can similarly define $Y_{\delta}(f)$.

\end{defn}

 Note that the norms used for $X_{**}$ control one extra vertical derivative when compared to the norms used for $Y_{**}$. We  consider $\dbar$ as a map from $X_{**}$ to $Y_{**}$. 

\begin{thm}\label{smooth dbar} If $\hat f$ is an allowable $\C\infty1$ family,   then 
$\dbar$ defines a $C^{1}$ map from  $X_{k, \delta}$ to $Y_{k, \delta}$.
\end{thm}
\pf
Our local form for $\dbar$ is 
\[\dbar \nu = E(\nu)+H(\nu)(\dvert\nu)\ .\]

By using Lemma \ref{product lemma} part \ref{product3}, we may estimate in this coordinate chart:

\[\norm{\dbar\nu}_{k,\delta}\leq \norm{E(\nu)}_{k,\delta}+c\norm{H(\nu)}^{1}_{k,\delta}\norm{\dvert\nu}_{k,\delta}\]

Lemma \ref{composition}  together with Lemma \ref{local form}   imply that the terms above $\norm{E(\nu)}_{k,\delta}$ and $\norm{H(\nu)}^{1}_{k,\delta}$  are bounded if $\norm{\nu}^{1}_{k,\delta}$ is bounded. Similarly, Corollary \ref{continuous composition} implies that $\dbar\nu$ depends continuously on $\nu$, therefore $\dbar$ gives a well defined map from $X_{k,\delta}$ to $Y_{k,\delta}$.
 
 \
 
The derivative is given by the following formula.
\begin{equation}\label{derivative formula}D\dbar(\nu)(\psi)=DE(\nu)(\psi(u))+DH(\nu)(\psi)(\dvert\nu)+H(\nu)(\dvert\psi)\end{equation}

Using Lemma \ref{product lemma} part \ref{product3}, estimate this derivative as follows.
\[\norm{D\dbar(\nu)(\psi)}_{k,\delta}\leq c\lrb{\norm{DE(\nu)}_{k,\delta}+\norm {DH(\nu)}_{k,\delta}^{1}\norm{\dvert \nu}_{k,\delta}+\norm{H(\nu)}_{k,\delta}^{1}}\norm{\psi}_{k,\delta}^{1}\]
The terms $\norm{H(\nu)}_{k,\delta}^{1}$, $\norm{DH(\nu)}_{k,\delta}^{1}$ and $\norm{DE(\nu)}_{k,\delta}$ are bounded by Lemma \ref{composition}, therefore  $D\dbar(\nu)$ defines a bounded  linear map from $X_{k,\delta}$ to $Y_{k,\delta}$. We must prove that $D\dbar(\nu)$ is continuous in $\nu$.
\[\begin{split}D\dbar(\nu_{1})(\psi)-D\dbar(\nu_{2})(\psi)&= \lrb{DE(\nu_{1})-DE(\nu_{2})}(\psi)+
\lrb{H(\nu_{1})-H(\nu_{2})}(\dvert\psi)
\\&+DH(\nu_{1})(\psi)(\dvert\nu_{1}-\dvert\nu_{2})
\\&+DH(\nu_{1})(\psi)(\dvert\nu_{2})-DH(\nu_{2})(\psi)(\dvert\nu_{2})
\end{split}\]
Using Lemma \ref{product lemma} part \ref{product3} on the above, there exists a positive constant $c$ (independent of $\nu_{i}$ and $\psi$) so that the following holds.
\[\begin{split}c\norm{D\dbar(\nu_{1})(\psi)-D\dbar(\nu_{2})(\psi)}_{k,\delta}&\leq \norm{DE(\nu_{1})-DE(\nu_{2})}_{k,\delta}\norm{\psi}_{k,\delta}^{1}
\\&+
\norm{H(\nu_{1})-H(\nu_{2})}_{k,\delta}^{1}\norm{\dvert\psi}_{k,\delta}
\\&+\norm{DH(\nu_{1})}_{k,\delta}^{1}\norm{\psi}_{k,\delta}^{1}\norm{\dvert\nu_{1}-\dvert\nu_{2}}_{k,\delta}
\\&+\norm{DH(\nu_{1})-DH(\nu_{2})}_{k,\delta}^{1}\norm{\dvert\nu_{2}}_{k,\delta}\norm\psi_{k,\delta}^{1}
\end{split}\]

The term $DE$ is $\C\infty1$ and $\Delta_{\mathcal S_{0}}DE=DE$, so Corollary \ref{continuous composition} implies that for a fixed $\nu_{1}$, the term $\norm{DE(\nu_{1})-DE(\nu_{2})}_{k,\delta}$ converges to zero as  $\norm{\nu_{1}-\nu_{2}}_{k,\delta}^{1}\rightarrow 0$. Similarly, Corollary \ref{continuous composition} tells us that $\norm{H(\nu_{1})-H(\nu_{2})}_{k,\delta}^{1}$ and $\norm{DH(\nu_{1})-DH(\nu_{2})}_{k,\delta}^{1}$ converge to zero as $\nu_{2}\rightarrow\nu_{1}$ in $\norm\cdot_{k,\delta}^{1}$. Lemma \ref{composition} implies that $\norm{DH(\nu_{1})}_{k,\delta}^{1}$ is bounded. Therefore, $D\dbar(\nu_{2})\rightarrow D\dbar(\nu_{1})$ as $\nu_{2}\rightarrow \nu_{1}$ in $\norm{\cdot}_{k,\delta}^{1}$, so  $\dbar:X_{k,\delta}\longrightarrow Y_{k,\delta}$ is $C^{1}$ as required.

\stop

\

The formula (\ref{derivative formula}) for $D\dbar$ in the proof above shows that it induces a map on $X_{\delta}(f)$. In other words, the restriction, $D\dbar(\nu)(f)$, of $D\dbar(\nu)$ to the domain of a curve, $f$, only depends on the restriction of $\nu$ to this curve, $\nu(f)$. Use the notation $D\dbar(f)$ to refer to the restriction to $f$ of $D\dbar$ at the zero-section.

\subsection{Results for variations of a single curve}

\

\begin{lemma}\label{fredholm}
For any allowable  $\C\infty1$ family of curves, and  a choice of weight $0< \delta< 1$, the linearization of  $\dbar$ at the zero-section restricted to any curve
\[D\dbar(f):X_{ \delta}(f)\longrightarrow Y_{ \delta}(f)\]
is Fredholm.

\

More generally, given a choice of weight $0<\delta<1$ for every coordinate chart in some allowable collection of coordinate charts on our family of curves, we may use the norms with mixed weights defined on page \pageref{mixed definition}, and 
\[D\dbar(f):X_{\text{mixed } \delta}(f)\longrightarrow Y_{\text{mixed } \delta}(f)\]
is Fredholm.

\end{lemma}

\pf

Note that the norm $\norm\cdot_{ \delta}$  restricted to a smooth component of an exploded curve is equivalent to the norm $L^{p}$ with exponential weight on the cylindrical ends given by $ \delta$, and $\norm\cdot_{ \delta}^{1}$ is equivalent to the $L^{p}$ norm with the same exponential weight on  derivatives plus the sup norm. The mixed version of these norms is similar, but uses different weights at different punctures. So $X_{ \delta}(f)$ or $X_{\text{mixed }\delta}(f)$, restricted to a smooth component, is  the corresponding Sobolev space $L^{1}_{p}$ with exponential weights $ \delta$ at punctures times a finite-dimensional space allowing sections which are asymptotic to constants instead of $0$ at cylindrical ends (restricted to the subspace which vanishes on our marked points if appropriate). From the proof of Theorem \ref{smooth dbar}, we have the following formula for $D\dbar(f)$ in local coordinates: 
\[D\dbar(\nu_{0}(f))(\psi)=DE(\nu_{0}(f))(\psi)+DH(\nu_{0}(f))(\psi)(\dvert\nu_{0}(f))+H(\nu_{0}(f))(\dvert\psi)\]
Here, $\nu_{0}$ indicates the zero-section, so the middle term is $0$. The conditions on $H$ from  (\ref{dbar form}) give that the last term  equals  $\frac12(\dvert\psi+J\circ\dvert\psi\circ j)$. This expression is only valid in local coordinates, as $\dvert\psi$ needs a local trivialization of $f^{*}T_{vert}\hat {\ex B}$ to make sense. To remedy this, choose a $\C\infty1$ connection $\nabla$ on $T_{vert}\hat {\ex B}$. Replacing $\dvert$ with $\nabla_{vert}$ in the above formula modifies the other term so that
\begin{equation}\label{Ddbar formula}D\dbar(\nu_{0}(f))(\psi)=E'(f)(\psi)+\frac12(\nabla\psi+J\circ\nabla\psi\circ j)\end{equation}
where restricted to any smooth component of the domain of $f$,  $E'(f)$ is smooth and decays exponentially on cylindrical ends with all weights less than $1$. Therefore, $D\dbar(\psi)$ is $\frac12(\nabla\psi+J\circ\nabla\psi\circ j)$ plus a compact operator.  
It is well known (and proved in \cite{mcowen}) that for $0<\delta<1$, $\psi\mapsto \frac12(\nabla\psi+J\circ\nabla\psi\circ j)$  is a Fredholm operator on the above weighted Soblev spaces restricted to any smooth component. It follows that $\dbar$ is Fredholm from $X_{\delta}(f)$ to $Y_{\delta}(f)$, and $X_{\text{mixed }\delta}(f)$ to $Y_{\text{mixed }\delta}(f)$. Therefore $D\dbar(f)$ is the Fredholm operator $\psi\mapsto\frac12(\nabla\psi+J\circ\nabla\psi\circ j)$ plus some compact operator, so it is Fredholm.

\stop

\

\begin{lemma}\label{linear regularity} If $\psi\in X_{\delta_{0}}(f)$ and $D\dbar(f) (\psi)$ is $\C\infty1$, then $\psi$ is $\C\infty1$.
\end{lemma}

\pf

Work in a coordinate chart where we can trivialize $(X,J)$. Using the corresponding flat, $J$--preserving connection, equation (\ref{Ddbar formula}) implies that
\[D\dbar(f)(\psi)=E'(\psi)+\dbar\psi\] 
where $\dbar\psi$ is the ordinary $\dbar$ operator applied to $\psi$ considered as a map to $\mathbb C^{n}$, and $E'$ is a $\C\infty1$ tensor vanishing on any edge contained in our coordinate chart.

Restricted to the smooth part of our coordinate chart, $E'$ is a smooth tensor, and $D\dbar(f)(\psi)$ is smooth, so a standard elliptic bootstrapping argument, such as given in \cite{MS}, implies that $\psi$ is smooth restricted to the smooth part of our coordinate chart. 

It remains to analyze the behavior of $\psi$ near edges. Let $\tilde z$ be a holomorphic coordinate on our edge and $\totl{\tilde z}=e^{t+i\theta}$. (As any $\C\infty1$ coordinate change will just change the norm on $X_{\delta}$ in a bounded way, we will assume that $\tilde z$ is one of our originally chosen coordinates.) We already know that $\psi(t,\theta)$ is smooth, and as $\psi\in X_{\delta_{0}}$, Lemma \ref{sobolev} implies that there is some constant $c$ so that 
\[\psi_{(t,\theta)}-\lim_{t\to -\infty}\psi(t,\theta)\leq ce^{\delta_{0} t} \ .\]
We must show the same for all weights $\delta<1$ and all derivatives of $\psi$. As $\psi$ is smooth and bounded as $t\to-\infty$ and $E'$ is  $\C\infty1$ and converges to $0$ as $t\to -\infty$, $E'(\psi)$ is smooth and converges exponentially to $0$  as $t\to-\infty$ with every weight less than $1$. It follows that $\dbar\psi$ is smooth and bounded by $e^{\delta t}$ as $t\to -\infty$ for all weights $\delta<1$. Below, we shall use that $\dbar\psi\in Y_{\delta}$ to show that $\psi\in X_{\delta}$. 

We may assume that our coordinates are valid for all $t\leq 0$, and apply Cauchy's integral formula: 
 \[\begin{split} \psi(t_0,\theta_0)&=
 \frac 1{2\pi}\int_0^{2\pi} \frac{\psi(0,\theta)}{1-e^{t_{0}}e^{i(\theta-\theta_0)}}d\theta
 + \frac 1{2\pi}\int_0^{2\pi}\int_{-\infty}^{0}\frac {\dbar\psi(t,\theta)}{1-e^{t_0-t}e^{i(\theta-\theta_0)}}dtd\theta \end{split}\]
 Therefore, 
 \[\begin{split} \psi(t_0,\theta_0)-\lim_{t\to-\infty}\psi(t,\theta)&=
 \frac 1{2\pi}\int_0^{2\pi} \frac{e^{t_{0}}e^{i(\theta-\theta_0)}\psi(0,\theta)}{1-e^{t_{0}}e^{i(\theta-\theta_0)}}d\theta
\\& + \frac 1{2\pi}\int_0^{2\pi}\int_{-\infty}^{0}\frac {e^{t_0-t}e^{i(\theta-\theta_0)}\dbar\psi(t,\theta)}{1-e^{t_0-t}e^{i(\theta-\theta_0)}}dtd\theta \ . \end{split}\]

The first integral is bounded by a constant times $e^{t_{0}}$ as $t_{0}\to-\infty$. We may bound the second integral by 
\[\norm{\dbar\psi}_{\delta} \frac 1{2\pi}\lrb{\int_0^{2\pi}\int_{-\infty}^{0}\lrb{\frac {e^{t_0-t+\delta t}}{\abs{1-e^{t_0-t}e^{i(\theta-\theta_0)}}}}^{q}dtd\theta}^{\frac 1q}\]
where $\frac 1q+\frac 1p=1$, and $p>2$ so $q<2$. It follows that the integral of $e^{q(t_{0}-t)}\abs{1-e^{t_{0}-t}e^{i\theta-\theta_{0}}}^{-q}$ for $t$ between $t_{0}-1$ and $t_{0}+1$ is finite, so the above integral between $t_{0}-1$ and $t_{0}+1$ is bounded by some constant times $e^{\delta t_{0}}$. For $t<t_{0}-1$, we may bound the integrand by $(2e^{\delta t})^{q}$, and for $t>t_{0}+1$, we may bound the integrand by $(2e^{t_{0}-t+\delta t})^{q}$. It follows that the above expression is bounded by a constant times $e^{\delta t_{0}}$. As this holds for all $\delta<1$, $\psi$ is in $X_{\delta}$ for all $\delta<1$. (Equivalently, as we are working over a single curve, $\psi$ is in  $X_{0,\delta}$ for all $\delta<1$.)

Suppose for induction that $\psi\in X_{k,\delta}$ for all $\delta<1$, and that if $\dbar\psi\in Y_{k,\delta}$ for all $\delta<1$, then $ \psi\in X_{k,\delta}$ for all $\delta<1$.
 The fact that $\psi\in X_{k,\delta}$ implies that  $E'(\psi)\in Y_{k+1,\delta}$, therefore, $\dbar\psi\in Y_{k+1,\delta}$. Let $v$ indicate either $\frac \partial{\partial t}$ or $\frac \partial{\partial \theta}$. We have that $\nabla_{v}\dbar\psi=\dbar\nabla_{v}\psi$, therefore $\dbar\nabla_{v}\psi\in Y_{k,\delta}$ for all $\delta<1$. Our inductive assumption then implies that $\nabla_{v}\psi\in X_{k,\delta}$ for all $\delta<1$, so $\psi\in X_{k+1,\delta}$ for all $\delta<1$.
 
 By induction, $\psi\in X_{k,\delta}$ for all $\delta<1$ and $k$, so $\psi\in \C\infty1$.

\stop

\

\begin{lemma}\label{index}
The index of $D\dbar(f):X_{\delta}(f)\longrightarrow Y_{\delta}(f)$ or $D\dbar(f):X_{\text{mixed }\delta}(f)\longrightarrow Y_{\text{mixed }\delta}(f)$  is
\[2c_{1}-2n(g+s-1)\]
where $c_{1}$  the first Chern number of the vectorbundle, $X$ or $f^{*}(T_{vert}\hat{\ex B},J)$, over the domain of $f$, $g$ is the total genus, and $s$ is the number of marked points where sections of $X_{\delta}$ are required to vanish, and $n$ is the complex rank of $X$. In particular, the index of $D\dbar(f)$ is invariant in connected families of curves.

\end{lemma}

\pf This is just the index of the operator $\psi\mapsto\frac12(\nabla\psi+J\circ\nabla\psi\circ j)$. This index  equals  the sum of the indexes of each smooth component minus $2n$ times the number of internal edges, because a section of $X_* (f)$ is equivalent to a choice of section over each smooth component with the extra requirement that these sections match at internal edges, and a section of $Y_*(f)$ is equivalent to a section over each smooth component, with no extra matching condition because sections of $Y_*(f)$ always vanish on edges. By a smooth component, we mean a connected component of the domain of $f$ minus all its edges. Each such component is a smooth punctured Riemann surface. The Riemann--Roch formula implies that the index on each smooth component, $v$,  is  $2c_{1}-2n(g_v+s_v-1)$ where $c_{1}$ is the first Chern number of the pullback of $(T_{vert}\hat{\ex B},J)$ to that component, $g_v$ is the genus of that component, and $s_v$ is the number of marked points on that component where sections of $X_{\delta}(f)$ are required to vanish. The sum over all components of $c_{1}$ and  $s_v$ does not change in connected families. The contribution of internal edges to the index is absorbed into the contribution of genus: $(1-g)$ is  $\sum_{v}(1-g_v)$ minus the number of internal edges. To see this, note that the total genus\footnote{For exploded curves $\ex C$, we could equivalently define this total genus as $\dim H^1(\ex C)/2$, which is invariant in families, by Proposition 11.4 of \cite{dre}, and equal to the genus of $\totl{\ex C}$ when $\totl{\ex C}$ has no nodes by Corollary 4.2 of \cite{dre}. } of an exploded curve $\ex C$ is the genus of the underlying nodal curve $\totl{\ex C}$. This genus is defined so that it is invariant in connected families, and in particular does not change when nodes are smoothed. Topologically, smoothing a node corresponds to replacing a pair of disks around each node with an annulus. The Euler characteristic, $2-2g$,  of the smoothed curve is therefore the Euler characteristic  of $\totl{\ex C}$ minus a pair of disks around each nodal point, or $\sum_v(2-2g_v)$ minus the number of disks removed.      In conclusion, the total index is $2c_{1}-2n(g+s-1)$, and therefore invariant in connected families.

\stop

\begin{remark}[Chern classes]The index formula $2c_{1}-2n(g+s-1)$ involves the total first Chern number $c_{1}$ of the pullback $(T_{vert}\hat {\ex B},J)$ to a curve. As in the case of smooth manifolds,  Chern classes of  $(\hat {\ex B},J)$ may be defined as differential forms in $\Omega^{*}(\hat{\ex B})$ using the Chern-Weil construction of Chern classes. In particular, choose a smooth Hermition metric and connection on $(T_{vert}\hat{\ex B},J)$. The curvature of this connection is a smooth two-form on $\hat {\ex B}$ with values in $\mathbb C$-linear endomorphisms of $(T_{vert}\hat{\ex B},J)$. In a local frame, this curvature two-form is defined from some smooth functions simply by exterior differentiation, addition and multiplication, so it must vanish on $\mathbb R$-nil vectors.   The Chern forms are obtained by taking symmetric polynomials in this curvature two form. In particular, the first Chern form is obtained by multiplying the trace of this curvature two-form by $i/2\pi$.  The proof that these Chern forms are closed real differential forms is entirely local, so following the same construction gives closed real differential forms which vanish on $\mathbb R$-nil vectors.  Therefore, using the notation of \cite{dre},  these Chern forms are closed differential forms in $\Omega^{*}(\hat{\ex B})$ which represent a Chern class in $H^{*}(\hat{\ex B})$. The total Chern number $c_{1}$ of a curve is then the integral over the curve of the first Chern form of $(T_{vert}\hat {\ex B},J)$.\end{remark}

\begin{prop}\label{curve regularity} In the case that our family $\hat{\ex C}\longrightarrow \ex F$ consists of a single curve, if  $\nu\in X_{\frac 12}$ has  $\dbar\nu\in\C\infty1$, then  $\nu\in \C\infty1$.

\end{prop}

\pf

Standard elliptic bootstrapping, as explained in \cite{MS}, gives that locally on smooth components of the domain, $\nu$ is smooth if $\dbar\nu$ is smooth.  It remains to analyze these solutions in a neighborhood of edges in order to verify that they are $\C\infty1$. As a $\C\infty1$ coordinate change will not affect the property of $\nu$ being in $X_{\frac 12}$, we may use holomorphic coordinates $\totl{\tilde z}=e^{t+i\theta}$ on a neighborhood of a puncture in a smooth component of the domain.

 Recall that $\dbar$ is in the form (\ref{dbar form}), so in a local coordinate chart, 
\[\dbar{\nu}=E(\nu)+\frac 12\lrb{H'(\nu)(d\nu)+J\circ H'(\nu)(d\nu\circ j))} \]
where $H'$ is a $\C\infty1$ map to the bundle of invertible $\mathbb R$-linear endomorphisms of $X$.

Suppose that $\nu\in X_{k,\delta_{0}}$. Then $\norm{E(\nu)}_{k+1,\delta}$ is finite for all $\delta<1$, and the inverse of $H'(\nu)(\cdot)$ gives a linear map bounded in the $\norm\cdot_{k+1,\delta}$ norm. Let $J(\nu)$ indicate the pullback of $J$ over the map $H'(\nu)$, and more generally, let $J(\cdot)$ indicate the pullback of $J$ over the map $H'(\cdot)$. Of course, $J(\cdot)$ is $\C\infty1$, so Lemma \ref{composition} implies that $\norm{J(\nu)}^{1}_{k,\delta}$ is finite if $\nu\in X_{k,\delta}$. Applying the inverse of $H'(\nu)(\cdot)$ to $\dbar\nu-E(\nu)$ gives the folowing.
\[\text{ If }\nu\in X_{k,\delta_{0}} \text{, then }\norm{\frac 12(d\nu+J(\nu)\circ d\nu\circ j)}_{k+1,\delta}  \text{ is finite  for all }\delta<1\ .\] 
Let $S$ be the stratum of our coordinate chart corresponding to the edge where $t\rightarrow -\infty$, and let $\dbar_{0}$ indicate the  standard $\dbar$ operator corresponding to $e_{S}\lrb{J(\nu)}$. We may rewrite the above displayed statement as follows.
\[\text{ If }\nu\in X_{k,\delta_{0}}\text{, then }\norm{\dbar_{0}\nu+\frac 12(\Delta_{S}J(\nu))\circ d\nu\circ j}_{k+1,\delta}  \text{  is finite  for all }\delta<1\ .\]

Let $L^{p}_{k}(D(r))$ indicate the $L^{p}$ norm with $k$ derivatives on a disk of radius $r$ around a point $(t,\theta)$.

Suppose that $\nu\in X_{\frac 12}$. Then as $t\rightarrow \infty$, $\norm{d\nu}_{L^{p}(D(1))}\leq e^{\frac t2}$, and Lemma \ref{sobolev} implies that $\abs{\Delta_{S}J(\nu)}\leq ce^{\frac t2}$. Therefore, $\norm{\dbar_{0}\nu}_{\delta}$ is  finite for all $\delta<1$. The argument of Lemma \ref{linear regularity} then implies that $\nu\in X_{\delta}$ for all $\delta<1$.

If $\nu\in X_{k,\delta}$, the fact that $\norm{J(\nu)}_{L^{p}_{k+2}D(1)}$ is bounded as $t\to -\infty$ allows us to use the following inequality from Proposition B.4.9 in \cite{MS}:  there exists some positive constant $c$ so that for $t$ large enough, 
\[\norm{\Delta_{S}\nu}_{L^{p}_{k+2}(D(\frac 12))}\leq c(\norm{d\nu+J(\nu)\circ d\nu\circ j}_{L^{p}_{k+1}(D(1))}+\norm{\Delta_{S}\nu}_{L^{p}_{k+1}(D(1))})\ .\]

It follows that if $\nu\in X_{k,\delta}$, then as $t\to-\infty$ 
\[\norm{\Delta_{S}\nu}_{L^{p}_{k+2}D(\frac 12)}\leq  e^{t\delta}\ .\]
 Therefore, $\nu$ is in $X_{k+1,\delta'}$ for all $\delta'<\delta$ if $\nu\in X_{k,\delta}$. As we have already established that $\nu\in X_{0,\delta}$ for all $\delta<1$, it follows that $\nu\in X_{k,\delta}$ for all $k$ and all $\delta<1$. So $\nu$ is $\C\infty1$ as required.

\stop

\

\subsection{A linear gluing theorem}

\

Throughout this section, we consider a linear operator in the form (\ref{Ddbar form}) over  a family of curves $\hat{\ex C}\longrightarrow \ex F$ with a finite collection of nonintersecting marked-point sections $s_{i}:\ex F\longrightarrow \hat{\ex C}$ staying away from  the edges of curves in $\hat{\ex C}$; sections of $X$ will vanish at the image of $s_{i}$. Throughout,   we shall use $V$ to indicate a finite-rank sub-vectorbundle $V\subset Y^{\infty,\underline 1}$ as in Definition \ref{sub-vectorbundle}.  When we have $V$ over an allowable family of curves, we always assume that $V$ is extendable.

\begin{defn}
Any  finite-rank sub-vectorbundle $V\subset Y^{\infty,\underline 1}$ over an allowable family of curves defines a closed subspace $V_{k, \delta}\subset Y_{k, \delta}$ consisting of all sections  contained in $V(f)$ when restricted  to  $f\in \ex F$. Define the natural projection

\[\pi_{V}:Y_{k, \delta}\longrightarrow\frac{ Y_{k, \delta}}{V_{k, \delta}}:=Y/V_{k, \delta}\ .\] 
\end{defn}
Note that  $V_{k,\delta}$ is a subset of $Y_{k,\delta}$, not $Y$, so $Y/V_{k,\delta}$ is only to be interpreted as a short-hand for $Y_{k,\delta}/V_{k,\delta}$, not as the quotient of $Y$ by $V_{k,\delta}$, or sections of a vectorbundle $Y/V$.
Similarly define the restriction of $\pi_{V}$ to $Y_{\delta}(f)$. In this case, the projection $\pi_{V}$ has finite-dimensional kernel equal to $V$ restricted to $f$, so 
\[ \pi_{V}\circ D\dbar(f):X(f)_{ \delta}\longrightarrow Y/V_{ \delta}(f)\]
is still Fredholm. 

\

\

We now prove a standard gluing theorem. First, we  describe a `gluing' and `cutting' map.

\begin{lemma}\label{gluing}
Given any curve $f$ in $\ex F$, there exists an allowable collection of coordinate charts on neighborhood $U$ of $f$ and for all $f'\in U$,  bounded linear maps
\[G:X_{\delta'}(f)\longrightarrow X_{\delta'}(f')\]
\[G:Y_{\delta'}(f)\longrightarrow Y_{\delta'}(f')\]
\[C: Y_{\delta'}(f')\longrightarrow Y_{\delta'}(f)\]
so that the following holds:
\begin{enumerate}
\item 
Assign  a weight  $\delta_{i}\in\{\delta',\delta+\delta'\}$ to each coordinate chart, but use $\delta'$ on every chart not containing an edge of both $f$ and $f'$.\footnote{The reason we must be cautious about where the higher weight is used is that $\C\infty1$ functions are only uniformly bounded using this higher weight on strata containing edges of curves. This is not a problem if the ${q}$ in Lemma \ref{k comparison} is $1$, or if ${q}(\delta +\delta')<1$.} Then, the following maps are bounded
\[G:X_{\text{mixed }\delta}(f)\longrightarrow X_{\text{mixed }\delta}(f')\]
\[G:Y_{\text{mixed }\delta}(f)\longrightarrow Y_{\text{mixed }\delta}(f')\]
\[C:Y_{\text{mixed }\delta}(f')\longrightarrow Y_{\text{mixed }\delta}(f)\]
 and the norms of such $G$ and $C$ are bounded uniformly independent of  $f'\in U$,  so long as the charts using the higher weight contain an edge of $f$ and $f'$. 
\item $G$ is a left inverse to $C$. In other words, $G\circ C:Y_{\delta'}(f')\longrightarrow Y_{\delta'}(f)\longrightarrow Y_{\delta'}(f')$ is the identity.

\item   For any $\epsilon>0$, there exists an open neighborhood of $f$ so that, for all $f'$ in this neighborhood, so long as the higher weights are only used where $f'$ has an edge,   the following holds.
\begin{enumerate}
\item \label{commute G and dbar} 
\[\norm{D\dbar(f')\circ G \psi-G\circ D\dbar(f)\psi}_{\text{mixed }\delta}\leq \epsilon \norm \psi^{1}_{\text{mixed }\delta}\]
\item\label{GV} Given a finite-rank sub-vectorbundle $V\subset Y^{\infty,\underline 1}(\hat f)$, we may additionally ensure  that for all $v\in V(f)$,
\[\norm{\pi_{V}Gv}_{\text{mixed }\delta}\leq \epsilon\norm v_{\text{mixed }\delta}\ .\] 
\end{enumerate}
\end{enumerate}
\end{lemma}

\

\pf

Choose an allowable collection of coordinate charts lifting a single coordinate chart $U$ on $\ex F$ so that the largest stratum of $\totb U$ contains $f$. Choose these coordinate charts $U_{i}\longrightarrow U$ small enough that we can trivialize both $X$ and $Y$ over them as in Definition \ref{allowable collection}. Note that in each coordinate chart, the trivialization of $Y$ may differ from the trivialization induced on  $X\otimes T^{*}_{vert}\hat{\ex C}$, because $Y(\hat f)\subset X\otimes T^{*}_{vert}\hat{\ex C}$ may not be a constant sub-bundle.

Roughly speaking, $G$ is $(1-\tilde \Delta_{S})$, where $S$ is the stratum of $\totb U$ containing $f$, but we shall describe it more accurately below.
Let $g$ denote a section of the relevant vectorbundle, written using the above coordinates.
Now define the maps $G$ and $C$ in these coordinate charts $\pi_{\ex F}:U_{i}\longrightarrow U$ as follows:

If $U_{i}$ is of type \ref{boring} from the definition of allowable coordinate charts on page \pageref{boring}, $U_{i}$ is the product of an open subset of $\et 11$ with $ U$. The coordinate $\tilde z_{1}$ on $\et 11$ is a coordinate for both $\pi_{\ex F}^{-1}(f)\cap U_{i}$ and for $\pi_{\ex F}^{-1}(f')\cap U_{i}$. In this case, we can define $G(\psi)(\tilde z_{1})=\psi(\tilde z_{1})$ and $C(\psi)(\tilde z_{1})=\psi(\tilde z_{1})$. Note that if $\psi$ vanishes at marked points, then $G(\psi)$ does too, as we chose our coordinates so that the position of marked points is independent of $U$. (This is a condition of being an allowable collection of coordinate charts.) Our assumption that transition maps for our vectorbundle are independent of $U$ ensures that $G$ and $C$ are compatible with coordinate changes between charts of type \ref{boring}.

If, on the other hand, $U_{i}$ is a chart of type \ref{interesting}, then $U_{i}$ is of the form $U_{i}:=\pi^{-1}_{\ex F}\lrb{U}\cap \{\abs{\tilde z_{1}}<c,\abs{\tilde z^{\beta}}<c \}$. Choose some smooth cut-off function $\rho:\et 11\longrightarrow\mathbb R$ satisfying the following.
\[\rho(\tilde z)=\begin{cases}0 \text{ if }\abs {\tilde z}>\frac c2 
\\ 1 \text{ if }\abs {\tilde z}\leq \frac c4\end{cases}\]

 Recall that $\tilde z_{1}\tilde z^{\beta}$ is a coordinate function on $U$, and coordinates for $\pi_{\ex F}^{-1}(f)$ are given by $\tilde z_{1}\in \et 11$ and $\tilde z^{\beta}\in \et 11$ so that $\tilde z_{1}\tilde z^{\beta}=\tilde z_{1}\tilde z^{\beta}(f)$, and $\abs {\tilde z_{1}}<c, \abs {\tilde z^{\beta}}<c$. We have assumed  that $\totb {\tilde z_{1}\tilde z^{\beta}(f)}> 0$ by assuming that $f$ is in the largest stratum of $\totb U$. In our coordinates, $\psi$ is some constant $x$ on the edge of $\pi_{\ex F}^{-1}(f)$, (where $\totl{\tilde z_{1}}=0=\totl{\tilde z^{\beta}}$). Define $\psi_{+}(\tilde z_{1}):=\psi(\tilde z_{1})$ when $\abs{\tilde z_{1}}>\abs{\tilde z^{\beta}}$, and extend $\psi_{+}$ to be $x$ elsewhere. Similarly define $\psi_{-}(\tilde z^{\beta})=\psi(\tilde z^{\beta})$ for $\abs{\tilde z^{\beta}}>\abs{\tilde z^{\alpha}}$, and extend $\psi_{-}$ to be $x$ everywhere else.

Note that $\pi_{\ex F}^{-1}(f')\cap U_{i}$ has the same coordinates $\tilde z_{1}$ and $\tilde z^{\beta}$, except $\tilde z_{1}\tilde z^{\beta}=\tilde z_{1}\tilde z^{\beta}(f')$. Define $G$ in these coordinate charts as follows.

\[G(\psi)(\tilde z_{1},\tilde z^{\beta}):=\rho(\tilde z^{\beta})(\psi_{+}(\tilde z_{1})-x)+\rho({\tilde z_{1}})(\psi_{-}(\tilde z^{\beta})-x)+x\]

Our assumptions about the transition maps between allowable collections of coordinate charts from Definiton \ref{allowable collection} ensure that $G$ is compatible with our transition maps,  so we have a globally defined map $G$.  Note also that if $\psi$ is a $\C\infty1$ section of $X(f)$ or $Y(f)$, $G(\psi)$ actually defines a $\C\infty 1$ extension of $\psi$ to neighborhood of $f$.

The norm of $G$ as a map from $X_{\text{mixed }\delta}(f)$ to $X_{\text{mixed }\delta}(f')$ is bounded independent of $f'$ in $U$. Indeed,  $\norm{\psi_{+}(\tilde z_{1})}^{1}_{\text{mixed }\delta}\leq\norm {\psi}^{1}_{\text{mixed }\delta}$ because $w_{0}(\tilde z_{1})$ is smaller at the fiber over $f$ than the fiber at $f'$, and similarly,  $\norm{\psi_{-}(\tilde z^{\beta})}^{1}_{\text{mixed }\delta}\leq \norm{\psi}^{1}_{{\text{mixed }\delta}}$. If our coordinate chart contains edges of both $f$ and $f'$, the cutoff function $\rho$ plays no role in the definition of $G$, because it is always $1$ where it multiplies a nonzero function. Otherwise, we always use the lower weight $\delta_{0}$ in our coordinate chart, and when we multiply $(\psi_{+}(\tilde z_{1})-x)$ by $\rho(\tilde z^{\beta})$, we increase its $\norm{\cdot}^{1}_{\delta_{0}}$--norm by a factor of at most $2\norm{\rho(\tilde z^{\beta})}^{1}_{\delta_{0}}$, which is finite because $\rho(\tilde z^{\beta})$ is smooth. Therefore, the norm of $G$ as a map from $X_{\text{mixed }\delta}(f)$ to $X_{\text{mixed }\delta}(f')$ is bounded by $1+ 4\norm{\rho(\tilde z^{\beta})}^{1}_{\delta_{0}}+4\norm{\rho(\tilde z_{0})}^{1}_{\delta_{0}}$.

 The argument bounding $G$ as a map from $Y_{\text{mixed }\delta}$ is similar, however, to define  $G:=Y_{\text{mixed }\delta}(f)\longrightarrow Y_{\text{mixed }\delta}(f')$ we used a trivialization of the vectorbundle $Y$ over a coordinate chart which does not necessarily coincide with the trivialization of $T_{vert}^{*}\hat{\ex C}\otimes X$ used to define the $Y_{\text{mixed } \delta}$ norm.
This does not matter as Lemma \ref{product lemma} implies that any extendable $\C\infty1$ change of coordinates on our vectorbundle gives an equivalent norm to $Y_{\delta'}$, and restricted to curves in our chart containing an edge, it gives an equivalent norm to $Y_{\delta+\delta'}$. Therefore, $G:=Y_{\text{mixed }\delta}(f)\longrightarrow Y_{\text{mixed }\delta}(f')$ is also uniformly bounded for all $f'$ in $U$ so long as the charts using the higher weight contain an edge of $f'$ and $f$.

\

Define the cutting map $C$ in these coordinates by 
\[C(\psi)(\tilde z_{1},\tilde z^{\beta})=\psi(\tilde z_{1})\text{ if }\abs{\tilde  z_{1}}>\sqrt{\abs {\tilde z_{1}\tilde z^{\beta}(f')}} \text{ and }\totl{\tilde z_{1}}\neq 0\]
\[C(\psi)(\tilde z_{1},\tilde z^{\beta})=\psi(\tilde z^{\beta})\text{ if }\abs{\tilde  z^{\beta}}\geq\sqrt{\abs {\tilde z_{1}\tilde z^{\beta}(f')}}\text{ and }\tot{\tilde z^{\beta}}\neq0\]
\[C(\psi)(\tilde z_{1},\tilde z^{\beta})=0 \text{ everywhere else.}\]

Recall that, within coordinate charts of type \ref{interesting}, $\abs{\tilde z_{1}\tilde z^{\beta}}<\frac {c^{2}}{16}$. Therefore  $G\circ C$ is the identity, because the cutoff functions involved are $1$ in the relevant regions.  Our  assumptions on transition maps within allowable collections of coordinate charts ensure that $C$ is defined independent of coordinates.
Observe also that if our trivialization of $Y(\hat f)$ in a coordinate chart coincided with the trivialization used to calculate norms,  $C:Y_{\text{mixed }\delta}(f')\longrightarrow Y_{\text{mixed }\delta}(f)$ would have norm bounded by $2$ because $w_{0}=\abs{\totl{\tilde z_{1}}}+\abs{\totl{\tilde z^{\beta}}}$, so after accounting for the different trivializations, the norm of $C$ is uniformly  bounded independent of $f'$ in $U$ so long as charts using the higher weight contain an edge of $f'$ and $f$. We have now verified the first three items of our lemma.

\

Recall that  $D\dbar(f)$ is in the following form.
\[D\dbar (f)(\psi)(\tilde z_{1},\tilde z^{\beta})=E'(\tilde z_{1},\tilde z^{\beta},f)(\psi(\tilde z_{1},\tilde z^{\beta}))+\frac 12(\dvert\psi+J\circ\dvert\psi\circ j)(\tilde z_{1},\tilde z^{\beta},f)\]
In the above, we write things as depending on $(\tilde z_{1},\tilde z^{\beta})$ even though these coordinates are related over $f$. Below, we shall sometimes use only one of these coordinates when we wish to emphasize dependence on that coordinate. In the above expression, $E'$ is some linear map which depends in a $\C\infty1$ way on position in $U_{i}$ and disappears on the edges of fibers of $\pi_{\ex F} $, and $J$ and $j$ also depend in a $\C\infty1$ way on position in $U_{i}$. Let us bound  the expression $D\dbar(f')\circ G \psi-G\circ D\dbar(f)\psi$  in the interesting case that $U_{i}$ is a chart of type \ref{interesting}. 

\

This calculation is complicated by the difference between the coordinates in which $D\dbar f$ is expressed and the coordinates in which the gluing map on $Y_{\text{mixed }\delta}$ is defined. In particular let $A(f,\tilde z_{1},\tilde z^{\beta},f)(\cdot)$ indicate the coordinate change, so in our current coordinates,
\[G:Y_{\text{mixed }\delta}(f)\longrightarrow Y_{\text{mixed }\delta}(f')\]
is equal to 
\[G(\nu)(\tilde z_{1},\tilde z^{\beta}):=A^{-1}(f',\tilde z_{1},\tilde z^{\beta})\lrb{\rho(\tilde z^{\beta})A(f,\tilde z_{1})\nu_{+}(\tilde z_{1})+\rho({\tilde z_{1}})A(f,\tilde z^{\beta})\nu_{-}(\tilde z^{\beta})}\ .
\]
So, 
\begin{equation}\label{llo}\begin{split}D\dbar(f')\circ G \psi&-G\circ D\dbar(f)\psi=
(\id-A^{-1}(f',\tilde z_{1})A(f,\tilde z_{1}))\rho(\tilde z^{\beta})(D\dbar(f)\psi)_{+}
\\&+(\id-A^{-1}(f',\tilde z^{\beta})A(f,\tilde z^{\beta}))\rho(\tilde z_{1})(D\dbar(f)\psi)_{-}
\\ & +D\dbar(f')\circ G(\psi)- \rho(\tilde z^{\beta})(D\dbar(f)\psi)_{+}-\rho(\tilde z_{1})(D\dbar(f)\psi)_{-}  \ .\end{split}\end{equation}
As $(\id-A^{-1}(f',\tilde z_{1})A(f,\tilde z_{1}))\rho(\tilde z^{\beta})$ is linear,  $\C\infty1$ and vanishes when $f'=f$, the $\norm{\cdot}_{\text{mixed }\delta}$ norm of the righthand side of the top line can be bounded by $\epsilon\norm\psi^{1}_{\text{mixed }\delta}$ by choosing $f'$ close enough to $f$, and the second line may be bounded similarly. It therefore suffices to bound the last line.

The  expression we need to bound is linear in $\psi$. Consider the case that  $\psi$ is a constant $x$. Then,  $G(x)=x$, and so our expression is 
\[E'(\tilde z_{1},\tilde z^{\beta},f')(x)-\rho(\tilde z^{\beta})E'(\tilde z_{1},f)_{+}(x)-\rho(\tilde z_{1})E'(\tilde z^{\beta},f)_{-}(x)\ .\] 
This is $\C\infty1$ and vanishes on the edges of all curves, and vanishes on $f$, so its $\norm\cdot_{\text{mixed }\delta}$ norm may be forced to be as small as we like compared to $x$ by choosing $f'$ close to $f$ and requiring that $f'$ always have an edge in any coordinate using the higher weight.

Using the above, we may reduce to the case that $\psi$ restricted to the edge of $f$ in our coordinate chart is $0$, in which case,  
\[G(\psi)=\rho(\tilde z^{\beta})\psi_{+}+\rho(\tilde z_{1})\psi_{-}\ .\]
Now, the last line of equation (\ref{llo}), breaks into 
\begin{equation}\label{bthis}
D\dbar(f')(\rho(\tilde z^{\beta})\psi_{+})- \rho(\tilde z^{\beta})(D\dbar(f)\psi)_{+}
\end{equation}
and a similar expression involving $\rho(\tilde z_{1})$ and $\psi_{-}$. To bound the last line of (\ref{llo}), it therefore suffices to bound the above expression (\ref{bthis}). Expand (\ref{bthis}) as follows:
\begin{equation}\label{bthis'}\begin{split}
&\rho(\tilde z^{\beta})(E'(\tilde z_{1},f')-E'(\tilde z_{1},f))(\psi_{+}(\tilde z_{1}))
\\ &+\frac 12\lrb{d_{vert}\lrb{\rho(\tilde z^{\beta})\psi_{+}(\tilde z_{1})}+J\circ d_{vert}\lrb{\rho(\tilde z^{\beta})\psi_{+}(\tilde z_{1})}\circ j)(\tilde z_{1},f')}
\\&-\frac {\rho(\tilde z^{\beta})}2(d_{vert}\psi_{+}(\tilde z_{1})+ J\circ d_{vert}\psi_{+}(\tilde z_{1})\circ j)(\tilde z_{1},f) 
\end{split}\end{equation}
As $\rho(\tilde z^{\beta})(E'(\tilde z_{1},f')-E'(\tilde z_{1},f))$ is $\C\infty1$ and vanishes on $f$ and the edges of curves, the top line of the above expression (\ref{bthis'}) may have its $\norm\cdot_{\text{mixed }\delta }$ norm bounded by $\epsilon\norm\psi_{\text{mixed }\delta}^{1}$ so long as we choose $f'$ close enough to $f$ and require that $f'$ has an edge contained in the charts on which the higher weight is used. It therefore suffices to bound the last two lines of the expression (\ref{bthis'}).

\

We may expand these last two lines of (\ref{bthis'}) multiplied by 2 as follows:
\begin{equation}\label{bthis2}\begin{split}
&d\rho(\tilde z^{\beta})\psi_{+}(\tilde z_{1})+(d\rho(\tilde z^{\beta})\circ j) (J\psi_{+}(\tilde z_{1}))
 \\& \rho(\tilde z^{\beta})\lrb{d_{vert}\psi_{+}(\tilde z_{1})+J\circ d_{vert}\psi_{+}(\tilde z_{1})\circ j)(\tilde z_{1},f')}
\\&- \rho(\tilde z^{\beta})\lrb{d_{vert}\psi_{+}(\tilde z_{1})+ J\circ d_{vert}\psi_{+}(\tilde z_{1})\circ j)(\tilde z_{1},f)}
\end{split}
\end{equation}

The first line of the above expression may be bounded as follows:  Note $\abs{d \rho(\tilde z^{\beta})}$ is bounded, and is supported  on the region $\{\frac c4<\abs{\tilde z^{\beta}}<\frac c2\}$. By choosing $\abs{\tilde z_{1}\tilde z^{\beta}(f')}$ small enough by choosing $f'$ close to $f$, we can ensure that $\tilde z_{1}$ is as small as we like in the above region, thus making  $\norm{\psi_{+}(\tilde z_{1})}_{\text{mixed }\delta}$  here as small as we like compared to $\norm \psi_{\text{mixed }\delta}^{1}$ because the weight involved in the calculation the  norm  here (at $f'$) involves a term $\abs{\totl{\tilde z^{\beta}}}\geq \frac c4$ and Lemma \ref{sobolev} gives us a bound on the size of $\psi_{+}(\tilde z_{1})$ by a uniform constant times $\abs{\totl{\tilde z_{1}}}^{\delta'}\norm{\psi}^{1}_{\delta'}$. It follows that we can  make the $\norm\cdot_{\text{mixed }\delta}$ norm of the first line as small as we like compared to $\norm\psi^{1}_{\text{mixed }\delta}$.

Therefore, it remains to bound the last two lines of the expression (\ref{bthis2}). These may be expanded as follows:

\[\begin{split}
&\rho(\tilde z^{\beta})\lrb{(J(\tilde z_{1},f')-J(\tilde z_{1},f))\circ d_{vert}\psi_{+}(\tilde z_{1})\circ j(\tilde z_{1},f'))}
\\&+\rho(\tilde z^{\beta})\lrb{ J(\tilde z_{1},f)\circ d_{vert}\psi_{+}(\tilde z_{1})\circ (j(\tilde z_{1},f')-j(\tilde z_{1},f)}\end{split}\]

As $(J(\tilde z_{1},f')-J(\tilde z_{1},f))$ and $ (j(\tilde z_{1},f')-j(\tilde z_{1},f))$ are $\C\infty1$ and vanish when $f'=f$, the $\norm\cdot_{\text{mixed }\delta}$ norm of the above  expression may be bounded by $\epsilon\norm\psi_{\text{mixed }\delta}^{1}$ by choosing $f'$ close enough to $f$ and requiring that the coordinate charts which use the higher weights contain an edge of $f'$.

Therefore, by choosing our neighborhood of $f$ small enough, we can achieve item \ref{commute G and dbar} from our lemma on charts of type \ref{interesting}. On the other charts of type \ref{boring}, proving the same thing is similar, but easier as it involves no cutoff functions and just analysis of how $A$, $E'$, $j$ and $J$ depend on position. 

\

We now just need to prove item \ref{GV} from our lemma. To see this, consider $\C\infty1$ sections of  $V$, $v_{1},\dotsc,v_{n}$ so that $\{v_{i}(f)\}$ is a basis for $V(f)$. Then $G(v_{i}(f))-v_{i}$ is  $\C\infty1$ and vanishes on the domain of $f$ and all strata in $\mathcal S_{0}$. Therefore, we can get $\norm {G(v_{i}(f))(f')-v_{i}(f')}_{\text{mixed }\delta}$ as small as we like by choosing $f'$  close to $f$. This in turn bounds $\norm{\pi_{V}G(v_{i}(f))}_{\text{mixed }\delta}$, and item \ref{GV} follows from the linearity of $G$:
\[\norm{\pi_{V}G(v)}_{\text{mixed }\delta}\leq \epsilon \norm v_{\text{mixed }\delta}\]

\stop

\
 
\begin{thm}\label{invertible1}
Suppose that  for some curve $f$ in $\hat f$,
\[ D\dbar(f):X^{\infty,\underline 1}(f)\longrightarrow Y^{\infty,\underline 1}(f)\]
is injective and has image complementary to $V(f)\subset Y^{\infty,\underline 1}(f)$ for some finite-rank sub-vectorbundle $V\subset Y^{\infty,\underline 1}(\hat f)$.

Then for any $0<\delta<1$, there exists an open neighborhood $U$ of   $f\in\ex F$ covered by an allowable collection of coordinate charts, so that, for all $f'\in U$, 
\[\pi_{V}\circ D\dbar(f'):X_{\delta\&\delta'}(f')\longrightarrow Y/V_{\delta\&\delta'}(f')\]
and 
\[\pi_{V}\circ D\dbar(f'):X_{\delta'}(f')\longrightarrow Y/V_{\delta'}(f')\]
are invertible, and the norm of the inverse is uniformly bounded independent of $f'\in U$. 


\end{thm}

\pf

Recall that $X_{\delta\&\delta'}(f')$, as defined  on page \pageref{equivalent norm}, uses the weight $w_{0}^{-\delta'}$ in coordinate charts not containing an edge of  $f'$ and the weight  $w_{0}^{-\delta-\delta'}$ in coordinate charts containing an edge of $f'$.

Choose an allowable collection of coordinate charts containing $f$ so that Lemma \ref{gluing} holds. Let $f'$ be some curve in the image of this allowable collection of coordinate charts. Assign a weight $\delta'$ to any coordinate chart that does not contain an edge of $f'$, and a weight $\delta+\delta'$ to any coordinate chart that contains a edge of $f'$. Consider the corresponding norm with these mixed weights.

\

Lemma \ref{linear regularity} implies that the kernel of $\pi_{V}\circ D\dbar(f)$ consists of $\C\infty1$ sections, so our assumptions and the fact that $\pi_{V}\circ D\dbar(f):X_{\text{mixed }\delta}(f)\longrightarrow Y/V_{\text{mixed }\delta}(f)$ is Fredholm imply that $\pi_{V}\circ D\dbar(f)$ has a bounded inverse.

Use the notation 
\[Q:=\lrb{\pi_{V}\circ D\dbar(f)}^{-1}\circ \pi_{V}\ .\]
 Assume that $U$ was chosen small enough that  the gluing and cutting maps, $G$ and $C$ from Lemma \ref{gluing}, satisfy the following conditions:
\begin{enumerate} 

\item $\norm{D\dbar(f)}$, $\norm Q$,  $\norm G$,   
and $\norm C$ are  bounded by some $M$.
\item \label{fjie1}
\[\norm{D\dbar(f')\circ G \psi-G\circ D\dbar(f)\psi}_{\text{mixed }\delta}\leq \epsilon \norm \psi^{1}_{\text{mixed }\delta}\]
\item\label{fjie2} For all $v\in V(f)$,
\[\norm{\pi_{V}Gv}_{\text{mixed }\delta}\leq \epsilon\norm v_{\text{mixed }\delta}\ .\]
\end{enumerate}
By choosing $U$ small enough, we may make the $\epsilon$ in the conditions above as small as we like while keeping $M$ the same.

\

Now consider the map $Q(f'):Y_{\text{mixed }\delta}(f')\longrightarrow X_{\text{mixed }\delta}(f')$ defined by

\[Q(f'):=G\circ Q\circ C \]
By exchanging  $ D\dbar(f')\circ G$ with $G\circ D\dbar(f)$ using the inequality (\ref{fjie1}), and using that $G\circ C$ is the identity, and then  using that $ D\dbar(f)\circ Q (C\nu)-C\nu\in V(f)$ and the inequality (\ref{fjie2}), we get the following: 
\begin{equation}\label{inverse inequality}\begin{split}\norm{\pi_{V}\lrb{D\dbar(f')\circ Q(f')\nu-\nu }}_{\text{mixed }\delta}&=
\norm{\pi_{V}\lrb{D\dbar(f')\circ G\circ Q\circ C\nu-\nu }}_{\text{mixed }\delta}
\\&\hspace {-3cm}\leq \epsilon \norm{QC\nu}^{1}_{\text{mixed }\delta}
+\norm{\pi_{V}\lrb{G\lrb{ D\dbar(f)\circ Q\circ C\nu-C\nu}}}_{\text{mixed }\delta}
\\&\hspace {-3cm}\leq\epsilon M^{2}\norm\nu_{\text{mixed }\delta}+\epsilon(M^{3}+M)\norm\nu_{\text{mixed }\delta}
\end{split}\end{equation}

As $V(f')$ is $n$-dimensional, there exists a linear inclusion  
\[i_{V}:Y/V_{\text{mixed }\delta}(f')\longrightarrow Y_{\text{mixed }\delta}(f')\]  with norm bounded by $2^{n}$ so that $ \pi_{V}\circ i_{V}$ is the identity. To see this for the case $n=1$, use the Hahn Banach theorem to  choose a linear functional $L$ on $Y_{\text{mixed }\delta}(f')$ so that $\abs{Lv}=\norm v_{\text{mixed }\delta}$ for $v\in V(f)$, and $\norm L=1$. Then the obvious inclusion with image $\ker L$ is bounded by $2$. Repeating this argument $n$ times gives the $n$-dimensional case. 

In particular, if we choose $\epsilon$ small enough, the above inequality (\ref{inverse inequality}) tells us that
\[\norm {\pi_{V}D\dbar(f')\circ (Q(f')\circ i_{V})-Id}<\frac 12\]
so a right inverse to $\pi_{V}D\dbar(f')$ is given by 
\[(Q(f')\circ i_{V})\lrb{\pi_{V}D\dbar(f')\circ (Q(f')\circ i_{V})}^{-1}\]
which is bounded by $M2^{n+1}$.

 Lemma \ref{index} gives that the index of $D\dbar$ does not change in connected families, so  this right inverse must be a genuine inverse to $\pi_{V}D\dbar$. We may repeat this argument for the different mixed norms corresponding to whether  edges of $f'$ are in our coordinate charts to obtain a uniform bound for the inverse of  $\pi_{V}\circ D\dbar(f'):X_{\delta\&\delta'}(f')\longrightarrow Y/V_{\delta\&\delta'}(f')$ for all $f'$ in $U$.

The argument for $\pi_{V}\circ D\dbar(f'):X_{\delta'}(f)\longrightarrow Y/V_{\delta'}(f')$ is identical, except the weight $\delta'$ is used in all coordinate charts. 

\stop

\subsection{Linear regularity results in families}
\label{linear section}

\

\begin{prop}\label{invertible2} 
If $\hat f$ is covered by an allowable collection of coordinate charts, and, for all $f$ in   $\hat f$, the maps \[\pi_{V}\circ D\dbar(f):X_{\delta'}(f)\longrightarrow Y/V_{\delta'}(f)\]
\[\pi_{V}\circ D\dbar(f):X_{\delta\&\delta'}(f)\longrightarrow Y/V_{\delta\&\delta'}(f)\]
are invertible with a uniformly bounded inverse, then
\[\pi_{V}\circ D\dbar:X_{0, \delta}\longrightarrow Y/V_{0, \delta}\] has a bounded inverse.
\end{prop}

\pf

We need to bound $\norm{\phi}^{1}_{0, \delta}$ in terms of $\norm{\pi_{V}\circ D\dbar \phi}_{0, \delta}$. Our assumption that the maps $D\dbar(f)$ have uniformly bounded inverses implies that  there exists some constant $c$ independent of $\phi$ so that the following holds.
\begin{equation}\label{ii1}\norm{\phi}^{1}_{\delta'}\leq c\norm{\pi_{V}\circ D\dbar\phi}_{\delta'}\end{equation}
\begin{equation}\label{ii2}\norm{\phi}^{1}_{\delta\&\delta'}\leq c\norm{\pi_{V}\circ D\dbar\phi}_{\delta\&\delta'}\end{equation}
When there are $N$ charts in our allowable collection of coordinate charts,  the constant $c$ is $4N$ times our assumed uniform bound,  because $\norm{\phi}^{1}_{\delta'}\leq 4N\sup_{f}\norm{\phi(f)}^{1}_{\delta'}$, and $\norm{\pi_{V}\nu}_{\delta'}\geq \sup_{f}\norm{\pi_{V}\nu(f)}_{\delta'}$, and analogous bounds hold using the  $\norm{\cdot}_{\delta\&\delta'}$ norm.

We will use the equivalent form of $\norm\cdot_{0,\delta}^{1}$ from Lemma \ref{equivalent norm} and Lemma \ref{global norm}  because this equivalent form  involves only weights $w_{\mathcal S}$ with $\dvert w_{\mathcal S}=0$. So, in a local coordinate chart $D\dbar w_{\mathcal S}\phi=w_{\mathcal S} D\dbar \phi$. If $\pi_{V}D\dbar \tilde\Delta_{\mathcal S}\phi$ was also equal to $\tilde\Delta_{\mathcal S}\pi_{V}D\dbar\phi$, then the above inequalities (\ref{ii1}) and (\ref {ii2})  would be adequate to quickly prove our lemma. The main part of the following proof  estimates the extent to which this fails to hold. 

\

Restrict attention to the set of coordinate charts $ U_{i}$ over a single chart $U$ in our allowable collection of coordinate charts. 
Let $\mathcal S$ indicate a set of strata  in $\totb{U}$, and $\tilde{\mathcal S}$ indicate the lift of $\mathcal S$ to a set of strata in $\totb{\tilde U_{i}}$.  Lemma \ref{global norm} on page \pageref{global norm} tells us that the norm $\norm\phi_{0,\delta}^{1}$ is equivalent to the norm
\[\norm\phi^{1}_{\delta'}+\max_{\mathcal S}\norm{w^{-\delta}_{\mathcal S}\tilde \Delta_{\mathcal S}\phi}^{1}_{\delta\&\delta'}\ .\]

 For induction, assume that for any  collection of strata $I$ in $\totb{U}$ with fewer than $\abs{\mathcal S}$ members,
 \begin{equation}\label{inductive}\norm{w^{-\delta}_{ I}\tilde \Delta_{I}\phi}^{1}_{\delta\&\delta'} \leq c\norm{\pi_{V}D\dbar\phi}_{0, \delta} \ .\end{equation}
 The case when $I$ has no members is satisfied because of the inequality (\ref{ii2}) and Lemma \ref{global norm}. In what follows, we  bound $\norm{w^{-\delta}_{\mathcal S}\tilde \Delta_{\mathcal S}\phi}^{1}_{\delta\&\delta'}$ by bounding $\norm{D\dbar w^{-\delta}_{\mathcal S}\tilde \Delta_{\mathcal S}\phi}_{\delta\&\delta'}$.

\

Work in a local coordinate chart where we can use the local form (\ref{Ddbar form}) of $D\dbar$.
\[D\dbar \phi=E'(\phi)+H(\dvert \phi)\]
The important facts are that $E'$ is a $\C\infty1$ tensor so that $\Delta_{\mathcal S_{0}}E'=E'$, and $H$ is a $\C\infty1$ tensor.

For any stratum $S$,
\[\Delta_{S}(D\dbar\phi)=D\dbar \Delta_{S}\phi+(\Delta_{S} E')(e_{S}\phi) +(\Delta_{S}H)(e_{S}\dvert\phi)\ .\]
Similarly, for the lift, $\tilde{\mathcal S}$ of our set of strata $\mathcal S$ in $\totb{U}$ to our coordinate chart,
\begin{equation}\label{commute Delta}\begin{split}D\dbar&\Delta_{\tilde {\mathcal S}}\phi= \Delta_{\tilde{\mathcal S}}(D\dbar\phi)
-\sum_{\emptyset\neq I\subset \tilde{\mathcal S} }\Delta_{I} E'\lrb{  e_{I}\Delta_{(\tilde{\mathcal S}-I)}\phi} +(\Delta_{I}H)\lrb{  e_{I}\Delta_{(\tilde {\mathcal S}-I)}\dvert \phi}\ .\end{split}\end{equation}

 We now bound the $\norm{w_{\mathcal S}^{-\delta}\cdot}_{\delta\&\delta'}$ norm of terms in the above sum. Each of the terms in (\ref{commute Delta}) vanish on all strata in $\tilde {\mathcal S}$, so if $\nu$ indicates one of the above terms, $\Delta_{\ti S}\nu=\nu$. For such $\nu$, the following inequality holds:
 \begin{equation}\label{compare}\norm{w^{-\delta}_{\mathcal S}\Delta_{\ti S}\nu}_{\delta\&\delta'}\leq c\norm{w^{-\delta}_{\ti S\cup \mathcal S_{0}}\Delta_{\ti S}\nu}_{\delta'}\end{equation}
 In the above, $\mathcal S_{0}$ indicates the set of strata on which $w_{0}$ vanishes. To see why the inequality (\ref{compare}) holds, note that it holds on strata $T$ so that $T^{c}=\emptyset$ because then the norm on the  lefthand side is just the $\norm\cdot_{\delta'}$ norm. The inequality (\ref{compare}) also holds trivially on strata in $\ti S$, because  both sides are $0$ restricted to strata in $\ti S$. It remains to show (\ref{compare}) holds on strata $T$ so that $T\notin\ti S$ and $T^{c}\neq\emptyset$. It suffices to show that in this case $ e_{T}w_{\ti S\cup\mathcal S_{0}}$ is bounded by a constant times $e_{T}w_{\ti S}w_{0}$, which is inequality (\ref{claim}) from the proof of Lemma \ref{equivalent norm} on page \pageref{claim}. Therefore  the inequality (\ref{compare})  holds.

 \
 
 We shall now bound terms involving $H$ in (\ref{commute Delta}).

\[\begin{split}\norm{w^{-\delta}_{\mathcal S}(\Delta_{I}H)\lrb{  e_{I}\Delta_{(\ti S-I)}\dvert \phi}}_{\delta\&\delta'}
&\leq c\norm{w^{-\delta}_{\ti S\cup \mathcal S_{0}}(\Delta_{I}H)\lrb{  e_{I}\Delta_{(\ti S-I)}\dvert \phi}}_{\delta'}
\\&\leq c'\norm{(w^{-\delta}_{I}\Delta_{I}H)\lrb{e_{I} w_{(\ti S\cup \mathcal S_{0})-I}^{-\delta} \Delta_{(\ti S-I)}\dvert \phi}}_{\delta'}
\\&\leq c''\norm{e_{I} w_{(\ti S\cup\mathcal S_{0})-I}^{-\delta} \Delta_{(\ti S-I)}\dvert \phi}_{\delta'}
\end{split}\]

We wish to use our inductive hypothesis to estimate this last term, but $\ti S-I$ may not be a lifted set of strata. To remedy this, note that \begin{equation}\label{ei}e_{I}\Delta_{\ti S-I}\Delta_{\mathcal S_{0}}=e_{I}\Delta_{\ti S'}\Delta_{S_{0}}\end{equation} where $\ti S'$ is the largest lifted collection of strata which is a subset of $\ti S-I$. Therefore $e_{I}\Delta_{\ti S-I}\dvert\phi=e_{I}\Delta_{\ti S'}\dvert \phi$. We also have that 
\begin{equation}\label{ei2}e_{I}w_{(\ti S\cup\mathcal S_{0})-I}^{-\delta}\leq c e_{I}w_{\ti S'}^{-\delta}w_{0}^{-\delta} \end{equation}  because if $\zeta$ is a smooth monomial from $w_{\ti S'}w_{0}$, $e_{I}\zeta$  must vanish on $I^{c}$, and therefore must  vanish on $(\ti S\cup \mathcal S_{0})-I$ and be dominated by $w_{(\ti S\cup\mathcal S_{0})-I}$. Therefore,
\[\norm{e_{I} w_{(\ti S\cup\mathcal S_{0})-I}^{-\delta} \Delta_{(\ti S-I)}\dvert \phi}_{\delta'}
 \leq c\norm{e_{I} w_{\ti S'}^{-\delta} \Delta_{\ti S'}\dvert \phi}_{\delta+\delta'}\ .
\]

If $I$ is not a lifted set of strata, then the  righthand side of the above inequality is bounded by  a constant times $\norm{ w_{\ti S'}^{-\delta} \Delta_{(\ti S')}\dvert \phi}_{\delta\&\delta'}$. So if $I$ is not a lifted set of strata, 

\begin{equation}\label{H0}\norm{w^{-\delta}_{\mathcal S}(\Delta_{I}H)\lrb{  e_{I}\Delta_{(\ti S-I)}\dvert \phi}}_{\delta\&\delta'}\leq c\norm{ w_{\ti S'}^{-\delta} \Delta_{\ti S'}\dvert \phi}_{\delta\&\delta'}\ .\end{equation}

On the other hand, if $I$ is a lifted set of strata, then $\ti S-I$ is lifted, so $\ti S'=\ti S-I$, and we can derive the above inequality (\ref{H0}) using Lemma \ref{product lemma} as follows:

\[\begin{split}\norm{w^{-\delta}_{\mathcal S}(\Delta_{I}H)\lrb{  e_{I}\Delta_{(\ti S-I)}\dvert \phi}}_{\delta\&\delta'}
&\leq c\norm{(w^{-\delta}_{I}\Delta_{I}H)\lrb{e_{I} w_{\ti S-I}^{-\delta} \Delta_{(\ti S-I)}\dvert \phi}}_{\delta\&\delta'}
\\&\leq c'\norm{e_{I} w_{\ti S'}^{-\delta} \Delta_{\ti S'}\dvert \phi}_{\delta\&\delta'}
\\&\leq c'\norm{ w_{\ti S'}^{-\delta} \Delta_{\ti S'}\dvert \phi}_{\delta\&\delta'}
\end{split}\]

\

\

Now, we need to bound the terms involving $E'$ in (\ref{commute Delta}).

\[\begin{split}\norm{w^{-\delta}_{\ti S}(\Delta_{I} E')\lrb{  e_{I}\Delta_{(\ti S-I)}\phi} }_{\delta\&\delta'}
&\leq c\norm{w^{-\delta}_{\ti S\cup \mathcal S_{0}}(\Delta_{I} E')\lrb{  e_{I}\Delta_{(\ti S-I)}\phi} }_{\delta'}\text{ using (\ref{compare})} 
\\&\leq c'\norm{(w^{-\delta}_{ I\cup\mathcal  S_{0}}\Delta_{I\cup \mathcal S_{0}} E')\lrb{e_{I} w_{\ti S-I}^{-\delta} \Delta_{(\ti S-I)}\phi}}_{\delta'} 
\\&\leq c''\sup\abs{e_{I}w^{-\delta}_{(\ti S-I)}\Delta_{(\ti S-I)}\phi}
\end{split}\]

We must  bound $\sup\abs{e_{I}w^{-\delta}_{(\ti S-I)}\Delta_{(\ti S-I)}\phi}$.   If $I$ is a lifted set of strata and $\ti S'=\ti S-I$, then the above inequality immediately gives
\begin{equation}\label{D0}\norm{w^{-\delta}_{\ti S}(\Delta_{I} E')\lrb{  e_{I}\Delta_{(\ti S-I)}\phi} }_{\delta\&\delta'}
\leq c \norm{w_{\ti S'}^{-\delta}\Delta_{\ti S'}\phi}^{1}_{\delta\&\delta'}\ .
\end{equation}
On the other hand, if $I$ is not a lifted set of strata, 
 there exists some stratum $S\in \ti S-I$ and $\bar S\in S^{c}$ so that $\bar S\in I$. Then we can use Lemma \ref{sobolev} and the observation that $e_{I}w^{-\delta}_{\ti S-I}$ is bounded by a constant times $e_{I}w^{-\delta}_{\ti S'}w_{0}^{-\delta}$ (which follows from (\ref{ei2})) to get the following: 
\[\begin{split}\sup\abs{e_{I}w^{-\delta}_{(\ti S-I)}\Delta_{(\ti S-I)}\phi}
\leq & c\sup\abs{e_{I}w^{-\delta}_{\ti S'}w_{0}^{-\delta}\Delta_{(\ti S-I)}\phi}
\\&=c\sup e_{I}w_{0}^{-\delta}\abs{e_{\bar S}\Delta_{S}\lrb{e_{I}w^{-\delta}_{\ti S'}\Delta_{(\ti S-I)}\phi}}
\\&\leq c'\norm{e_{\bar S}\dvert \lrb{e_{I}w^{-\delta}_{\ti S'}\Delta_{(\ti S-I)}\phi}}_{\delta+\delta'}
\\&\ = c'\norm{e_{I}w^{-\delta}_{\ti S'}\Delta_{(\ti S-I)}\dvert \phi}_{\delta+\delta'}
\\&\ = c'\norm{e_{I}w^{-\delta}_{\ti S'}\Delta_{\ti S'}\dvert \phi}_{\delta+\delta'}\text{ using (\ref{ei})}
\\&\ \leq c''\norm{w^{-\delta}_{\ti S'}\Delta_{\ti S'}\dvert \phi}_{\delta\&\delta'}\text{ as }I^{c}\neq\emptyset\ .
\end{split}\]

Therefore, inequality (\ref{D0}) holds for all nonempty sets of strata $I\subset \ti S$. Using the inequalities (\ref{H0}) and (\ref{D0}) along with equation (\ref{commute Delta}) gives

\begin{equation} \label{commute Delta 2}
\norm{D\dbar w^{-\delta}_{\ti S} \Delta_{\ti S}\phi-w_{\ti S}^{-\delta} \Delta_{\ti S}D\dbar\phi}_{ \delta\&\delta'}\leq  c\sum_{\ti S'\subsetneq\ti S}\norm{w_{\ti S'}^{-\delta}\Delta_{\ti S'}\phi }^{1}_{\delta\&\delta}\ .\end{equation}

If $\ti S$ consists entirely of strata $S$ so that $S^{c}=\emptyset$, then $\tilde \Delta_{\mathcal S}=\Delta_{\tilde S}$, and we obtain the inequality below. Otherwise,  we may apply the estimates of Lemma \ref{tilde estimates}  to the above inequality to exchange $\Delta_{\ti S}$ with $\tilde \Delta_{S}$ and get error terms involving $\tilde \Delta_{\mathcal S'}$ where $\mathcal S'\subsetneq\mathcal S$, and also apply the estimates of Lemma \ref{tilde estimates} to the  righthand side of (\ref{commute Delta 2}) to get
 \begin{equation} \label{commute Delta 3}
\norm{D\dbar w^{-\delta}_{\mathcal S} \tilde \Delta_{\mathcal S}\phi-w_{\ti S}^{-\delta}\tilde \Delta_{\mathcal S}D\dbar\phi}_{ \delta\&\delta'}\leq  c\sum_{\mathcal S'\subsetneq\mathcal S}\norm{w_{\mathcal S'}^{-\delta}\tilde\Delta_{\mathcal S'}\phi }^{1}_{\delta\&\delta}\ .\end{equation}

The  righthand side of (\ref{commute Delta 3}) is bounded by $c\norm{\pi_{V}D\dbar\phi}_{0, \delta}$ by our inductive assumption (\ref{inductive}). As (\ref{commute Delta 3}) applies on all the coordinate charts $U_{i}$ over $U$, we can  apply $\pi_{V}$ to the  lefthand term in  (\ref{commute Delta 3}), use the triangle inequality, and rearrange to get

\begin{equation}\label{commute Delta 4}\begin{split}\norm{\pi_{V}D\dbar w^{-\delta}_{\mathcal S} \tilde\Delta_{\mathcal S}\phi}_{ \delta\&\delta'}\leq & 
\norm{\pi_{V}w_{\mathcal S}^{-\delta} \tilde \Delta_{\mathcal S}D\dbar\phi}_{\delta\& \delta'}
\\ &+ c\norm{\pi_{V}D\dbar\phi}_{0, \delta}\ .\end{split}\end{equation}

 We now estimate this  righthand side of this inequality by a constant times $\norm{\pi_{V}D\dbar\phi}_{0, \delta}$. (This estimation is trivial if $V$ has rank $0$, but it will take us over two pages.) After that, the inequality (\ref{ii2}) can be applied to complete the inductive argument.

Let $v'$ be any section of $V$ with $\norm{v'}_{0,\delta}$ finite so that 
\begin{equation}\label{def v'}2\norm{\pi_{V}D\dbar \phi}_{0, \delta}\geq \norm{D\dbar \phi-v'}_{0, \delta}\ .\end{equation}
 Then, the following holds.
\begin{equation}\label{v'1}\begin{split}
\norm{\pi_{V}w_{\mathcal S}^{-\delta} \tilde\Delta_{\mathcal S}D\dbar\phi}_{\delta\& \delta'}
&:=\inf_{v}\norm{w_{\mathcal S}^{-\delta} \tilde\Delta_{\mathcal S}D\dbar\phi-v}_{\delta\& \delta'}
\\&\leq \norm{w_{\mathcal S}^{-\delta} \tilde\Delta_{\mathcal S}(D\dbar\phi-v')}_{\delta\& \delta'}
+\inf_{v}\norm{w_{\mathcal S}^{-\delta}\tilde\Delta_{\mathcal S}v'-v}_{\delta\& \delta'}
\\&\leq c\norm{\pi_{V}D\dbar\phi}_{0, \delta} +\inf_{v}\norm{w_{\mathcal S}^{-\delta}\tilde\Delta_{\mathcal S}v'-v}_{\delta\& \delta'}
\\&\ \ = c\norm{\pi_{V}D\dbar\phi}_{0, \delta} +\norm{\pi_{V}w_{\mathcal S}^{-\delta}\tilde\Delta_{\mathcal S}v'}_{ \delta\&\delta'}
\end{split}\end{equation}

\

We must show that $\norm{\pi_{V}w_{\mathcal S}^{-\delta}\tilde\Delta_{\mathcal S}v'}_{ \delta\&\delta'}$ can be bounded in terms of  $\norm{\pi_{V}D\dbar\phi}_{0, \delta}$.

\

 In our coordinate chart $U$,  choose a trivialization of $V$ using a basis of $\C\infty1$ sections $\{v_{i}\}$. We then have 
\[v'=\sum_{i}f_{i}v_{i}\]
where $\{f_{i}\}$ is some collection of fiberwise-constant real-valued functions. Choose $v_{i}$  so that for some constant $c$, 
\begin{equation}\label{abs assumption}\abs {g_{i}}\leq c\norm {\sum_{i}g_{i}v_{i}}_{ \delta'}\ .\end{equation}
for any collection of functions $g_{i}$ on $U$.

 Using the identity $\Delta_{S}fg=(\Delta_{S}f)g+(e_{S}f)\Delta_{S}g$ repeatedly, we get 
\begin{equation}\label{missed}\Delta_{\ti S}v'=\sum_{I\subset \ti S,i}\lrb{\Delta_{I}e_{(\ti S-I)}f_{i}}
\lrb{\Delta_{(\ti S-I)}v_{i} }\ .\end{equation}
As $f_{i}$ is fiberwise-constant,   $\Delta_{S}e_{\bar S}f_{i}=0$ when $\bar S\in S^{c}$. Therefore, we can rewrite the above expression using only lifted sets of strata. 
\begin{equation}\label{tilde}\Delta_{\ti S}v'=\sum_{ I\subset \mathcal S,i}\lrb{\Delta_{\tilde I}e_{(\ti S-\tilde I)}f_{i}}
\lrb{\Delta_{(\ti S-\tilde I)}v_{i} }\end{equation}
We need to estimate $\tilde \Delta_{\mathcal S}v'$; this  equals  $\Delta_{\ti S}v'$ unless we're in a coordinate chart of type \ref{interesting} with coordinates including $ \abs{\tilde z_{1}}<c$ and $\abs{\tilde z^{\beta}}<c$.  Using  notation from Definition \ref{tilde delta}, write  $I^{+}$ for the subset of $\tilde I$ where $\totb{\tilde z_{1}}=\e 0$, and $I^{-}$ for the subset where $\totb{\tilde z^{\beta}}=\e 0$. Then 
\[\tilde \Delta_{\mathcal S}=\Delta_{\ti S}-\check\rho^{+}(\Delta_{\ti S}-\Delta_{\mathcal S^{+}})-\check\rho^{-}(\Delta_{\ti S}-\Delta_{\mathcal S^{-}})\]
where $\check \rho^{+}$ and $\check\rho^{-}$ are some cutoff functions. As $f_{i}$ is constant on fibers, $e_{I^{+}}f=e_{\tilde I}f_{i}=e_{I}f_{i}$, where in $e_{I}f_{i}$,  we consider $f_{i}$ as a function on $U$. Similarly,  $\Delta_{\tilde I}f_{i}=\Delta_{I^{+}}f_{i}=\Delta_{I^{-}}f_{i}=\Delta_{I}f_{i}$, so $\Delta_{\tilde I}f_{i}=\tilde \Delta_{I}f$. As every stratum in $\mathcal S$ has a unique lift to $\mathcal S^{+}$, equation (\ref{missed}) implies that
\begin{equation}\label{plus}\Delta_{\mathcal S^{+}}v'=\sum_{I\subset \mathcal S,i}\lrb{\Delta_{ I}e_{(\mathcal S- I)}f_{i}}
\lrb{\Delta_{(\mathcal S- I)^{+}}v_{i}}\end{equation}
and similarly, 
\begin{equation}\label{minus}\Delta_{\mathcal S^{-}}v'=\sum_{I\subset \mathcal S,i}\lrb{\Delta_{ I}e_{(\mathcal S- I)}f_{i}}
\lrb{\Delta_{(\mathcal S- I)^{-}}v_{i}}\end{equation}
so using equations (\ref{tilde}), (\ref{plus}) and (\ref{minus}) on a coordinate chart of type \ref{interesting}, we get
\begin{equation}\label{mnb}\tilde \Delta_{\mathcal S}v'=\sum_{I\subset\mathcal S,i}\lrb{\Delta_{ I}e_{(\mathcal S- I)}f_{i}}
\lrb{\tilde\Delta_{(\mathcal S- I)}v_{i}}\ .\end{equation}
We've proved the above equation (\ref{mnb}) for coordinate charts of type \ref{interesting}, and on all other coordinate charts, equation (\ref{mnb}) is equivalent to equation (\ref{tilde}), so (\ref{mnb}) is valid for all our coordinate charts over $U$.

We can apply equation (\ref{mnb}) along with inequality (\ref{abs assumption}) to bound $\abs{w_{\mathcal S}^{-\delta}\Delta_{\mathcal S}f_{i}}$ as follows:
\[\begin{split}\sup\abs{w^{-\delta}_{\mathcal S}\Delta_{\mathcal S}f_{i}}&\leq c\norm{\sum_{i}\lrb{w^{-\delta}_{\mathcal S}\Delta_{\mathcal S}f_{i}}v_{i} }_{\delta'}
\\&\leq c\norm{w_{\mathcal S}^{-\delta}\tilde\Delta_{\mathcal S}v'}_{\delta'}+c\sum_{I\subsetneq\mathcal S,j}\norm {w_{\mathcal S}^{-\delta}\lrb{\Delta_{ I}e_{(\mathcal S- I)}f_{j}}
\lrb{\tilde\Delta_{(\mathcal S- I)}v_{j}}}_{\delta'}
\\&\leq c\norm{w_{\mathcal S}^{-\delta}\tilde\Delta_{\mathcal S}v'}_{\delta'}+c'\sum_{I\subsetneq\mathcal S,j}\norm {\lrb{w_{I}^{-\delta}\Delta_{ I}e_{(\mathcal S- I)}f_{j}}
\lrb{w_{\mathcal S-I}^{-\delta}\tilde\Delta_{(\mathcal S- I)}v_{j}}}_{\delta'}
\\&\leq c\norm{w_{\mathcal S}^{-\delta}\tilde\Delta_{\mathcal S}v'}_{\delta'}+c''\sum_{I\subsetneq\mathcal S,j}\sup\abs {w_{I}^{-\delta}\Delta_{ I}e_{(\mathcal S- I)}f_{j}}
\end{split}\]
Using the above inequality again on the terms with fewer strata gives the following inequality.
\begin{equation}\label{bound abs}\sup\abs{w^{-\delta}_{\mathcal S}\Delta_{\mathcal S}f_{i}}\leq c\sum_{I\subseteq \mathcal S}\norm{w_{ I}^{-\delta}\tilde\Delta_{ I}v'}_{\delta'}\end{equation}
As the sections $v_{i}$ are $\C\infty1$, the norm $\norm{w^{-\delta}_{ I} \tilde\Delta_{ I}v_{i} }_{\delta\&\delta'}$ is bounded. Now bound terms on the  righthand side of equation (\ref{mnb}) as follows:  
\[\begin{split}\norm{\pi_{V}w^{-\delta}_{\mathcal S}\lrb{\Delta_{ I}e_{(\mathcal S- I)}f_{i}}
\lrb{\tilde \Delta_{(\mathcal S- I)}v_{i} }}_{\delta\& \delta'}\hspace{-3cm}&
\\&\leq c\norm{\pi_{V}\lrb{w^{-\delta}_{ I}\Delta_{ I}e_{(\mathcal S- I)}f_{i}}
\lrb{w^{-\delta}_{(\mathcal S- I)}\tilde\Delta_{(\mathcal S- I)}v_{i} }}_{\delta\& \delta'}
\\&\leq c\sup\abs{w^{-\delta}_{ I}\Delta_{ I}e_{(\mathcal S- I)}f_{i}}
\norm{\pi_{V}w^{-\delta}_{(\mathcal S- I)}\tilde \Delta_{(\mathcal S- I)}v_{i} }_{\delta\& \delta'} 
\\&\leq c'\sum_{I'\subseteq I}\norm{w_{ I'}^{-\delta}\tilde \Delta_{ I'}v'}_{ \delta'}\text{ using (\ref{bound abs})}
\end{split}\]
When considering $\pi_{V}$ of equation (\ref{mnb}), we can remove the term where $I=\mathcal S$, because $(\Delta_{\mathcal S}f_{i})v_{i}\in V$. In the remaining terms, $f_{i}$ is always being acted on by $\Delta_{ I}$ where $\abs{I}<\abs{\mathcal S}$. Therefore, using equation (\ref{mnb}) and the above inequality, 
\[\begin{split}\norm{\pi_{V}w^{-\delta}_{\mathcal S}v'}_{\delta\&\delta'}&\leq c\sum_{I\subsetneq \mathcal S}\norm{w_{ I}^{-\delta}\tilde \Delta_{ I}v'}_{ \delta'}
\\& \leq c'\lrb{\sum_{I\subsetneq\mathcal S}\norm{w_{ I}^{-\delta}\tilde \Delta_{ I}D\dbar \phi}_{\delta'}+\norm{\pi_{V}D\dbar \phi}_{0,\delta}} \text{ using (\ref{def v'})}
\\&\leq c''\lrb{\sum_{I\subsetneq\mathcal S}\norm{ w_{ I}^{-\delta}\tilde\Delta_{ I} \phi}^{1}_{\delta\&\delta'} +\norm{\pi_{V}D\dbar \phi}_{0,\delta}} \text{ using (\ref{commute Delta 3})}
\\&\leq c'''\norm{\pi_{V}D\dbar \phi}_{0,\delta} \text{ using (\ref{inductive}).}
\end{split}\]

Using the above inequality with equations (\ref{commute Delta 4}) and (\ref{v'1}), we get
\[\norm{\pi_{V}D\dbar w_{\mathcal S}^{-\delta}\tilde \Delta_{\mathcal S}\phi}_{\delta\& \delta'}
\leq c\norm{\pi_{V}D\dbar\phi}_{0, \delta}\ .\]
Applying inequality (\ref{ii2}) to the  lefthand side of the above gives
\[\norm{ w_{\mathcal S}^{-\delta}\tilde \Delta_{\mathcal S}\phi}^{1}_{\delta\& \delta'}
\leq c\norm{\pi_{V}D\dbar\phi}_{0, \delta}\ .\]
This completes the inductive argument that the above expression holds for all sets of strata $\mathcal S$ in $\totb U$. Therefore, using the equivalent norm from Lemma \ref{global norm}, we've proved the required inequality on the set of coordinate charts over $U$.
\begin{equation}\label{required inequality}\norm{\phi}^{1}_{0, \delta}\leq c\norm{\pi_{V}D\dbar\phi}_{0, \delta}\end{equation}
It follows that inequality (\ref{required inequality}) holds on our entire family, as our entire family is covered by a finite collection of allowable coordinate charts.
We shall use the above inequality (\ref{required inequality}) to see that $\pi_{V}D\dbar:X_{0,\delta}\longrightarrow Y/V_{0,\delta}$ has a bounded inverse, determined by the individual inverses to $\pi_{V}D\dbar (f):X_{0,\delta}(f)\longrightarrow Y/V_{0,\delta}(f)$ by
\[\lrb{\lrb{\pi_{V}D\dbar}^{-1}(\theta)}(f):=\lrb{\pi_{V}D\dbar(f)}^{-1}(\theta(f))\ .\]

Let us check that $\lrb{\pi_{V}D\dbar}^{-1}(\theta)$ actually is in $X_{0,\delta}$. Approximate $\theta$ by a $\C\infty1$ section $\theta'$. Then, for each $f$, Lemma \ref{linear regularity} implies that $\lrb{\pi_{V}D\dbar(f)}^{-1}(\theta(f))$ is $\C\infty1$. We can then choose a $\C\infty1$ extension $\phi_{f}$ of $\lrb{\pi_{V}D\dbar(f)}^{-1}(\theta(f))$ to a neighborhood of $f$, (using, for example, the map $G$ from Lemma \ref{gluing}).  So, $\pi_VD\dbar \phi_f$ is $\C\infty1$, and agrees with $\theta'$ over $f$. This implies that $\norm{\pi_VD\dbar\phi_f-\theta'}_{0,\delta}$ can be made arbitrarily small by restricting to a suitably small neighborhood of $f$. We can patch together such $\phi_{f}$ using a partition of unity (on $\ex F$) to define a section $\phi'$ in $X^{\infty, \underline 1}(\hat f)$ with $\norm{\pi_{V}D\dbar (\phi')-\theta'}_{0,\delta}$  as small as we like. There therefore exists a sequence of sections $\phi'$ in  $X^{\infty,\underline 1}$ with $\pi_{V}D\dbar (\phi')$ converging to $\theta$ in $Y/V_{0,\delta}$. The bound (\ref{required inequality}) then implies that this sequence $\phi'$ converges in  $X_{0,\delta}$, and converges to $\lrb{\pi_{V}D\dbar}^{-1}(\theta)$.

\stop 

\

\begin{prop}\label{invertible3}
If 
$\pi_{V}D\dbar:X_{0, \delta}\longrightarrow Y/V_{0, \delta}$
 has a bounded inverse, then
\[\pi_{V}D\dbar:X_{k, \delta}\longrightarrow Y/V_{k, \delta}\]
 has a bounded inverse.
\end{prop}

\

\pf 

This proposition is a tautology for $k=0$. Assume for induction that the corresponding inequality for the case $k-1$ holds, so
\[\norm{\phi}^{1}_{k-1,\delta}\leq c\norm{\pi_{V}D\dbar\phi}_{k-1, \delta}\ .\]
Suppose that $w$ is an extendable $\C\infty1$ vectorfield on the total space of our family, $\hat{\ex C}$, and that $w$  is the lift of some vectorfield on the base $\ex F$. In a coordinate chart,  
\begin{equation}\label{formula}\begin{split}
\nabla_{w}(D\dbar\phi)&=\nabla_{w}\lrb{E'(\phi)+H(\dvert\phi)}
\\&=D\dbar(\nabla_{w}\phi)+(\nabla_{w}E')(\phi)+(\nabla_{w}H)(\dvert \phi)+H([\dvert,\nabla_{w}]\phi) \end{split}\end{equation}
where $E'$ and $H$ are $\C\infty1$ tensors so that $\Delta_{\mathcal S_{0}}E'=E'$.
Note that $[\dvert,\nabla_{w}]$ is an  extendable  $\C\infty1$ first-order operator that involves derivatives only in the vertical direction,  because $w$ is the lift of a vectorfield on $\ex F$.  Therefore, there exists a constant $c$, depending on $w$, so that
 \[\norm{D\dbar\nabla_{w}\phi-\nabla_{w}\lrb{D\dbar \phi}}_{k-1,\delta}\leq c\norm{\phi}^{1}_{k-1,\delta}\ .\]
  We have proved the above inequality for a single coordinate chart. As our family is covered by a finite collection of such coordinate charts, the above inequality holds globally. Therefore, taking $\pi_{V}$ of the  lefthand side gives the following inequality where the constant $c$ depends on $w$:
  \begin{equation}\label{hjlk1}\norm{\pi_{V}D\dbar\nabla_{w}\phi}_{k-1,\delta}\leq  \norm{\pi_{V}\nabla_{w}(D\dbar\phi)}_{k-1, \delta}+c\norm{\phi}^{1}_{k-1, \delta}\end{equation}
We wish now to bound $\norm{\pi_{V}\nabla_{w}(D\dbar\phi)}_{k-1, \delta}$ by  $\norm{\pi_{V}D\dbar \phi}_{k,\delta}+\norm{\phi}^{1}_{k-1, \delta}$ times something depending on $w$  (a particularly easy task if $V$ is a rank-zero vectorbundle). Choose some section $v$ of $V$ with $\norm v_{k,\delta}$ finite so that 
\begin{equation}\label{hjlk4}2 \norm{\pi_{V}D\dbar\phi}_{k, \delta}\geq \norm{D\dbar\phi-v}_{k, \delta}\ .\end{equation}
 Therefore there exists some constant $c$ depending on $w$ so that
\begin{equation}\label{hjlk}\begin{split}c\norm{\pi_{V}D\dbar\phi}_{k, \delta}&\geq \norm{\nabla_{w}(D\dbar\phi -v)}_{k-1,\delta}
\\&\geq \norm{\pi_{V}\nabla_{w}(D\dbar\phi)}_{k-1, \delta}-\norm{\pi_{V}\nabla_{w}v}_{k-1, \delta}\ .\end{split}\end{equation}
Work in some finite collection of extendable coordinate patches on $\ex F$ on which $V$ is a trivial vectorbundle, and choose some basis $v_{1}\dotsc v_{n}$ of $\C\infty1$ sections of $V$ so that we can write any section of $V$ as $v=\sum f_{i}v_{i}$ where $\norm {f_{i}}_{k-1, \delta}\leq \norm{v}_{k-1, \delta}$. Note that
\[\pi_{V}\nabla_{w}(f_{i}v_{i})=\pi_{V}f_{i}\nabla_{w}v_{i}\ .\]
Therefore we can bound the right-most term in (\ref{hjlk}) using Lemma \ref{product lemma} as follows:
 \begin{equation}\label{hjlk3}\begin{split}\norm{\pi_{V}\nabla_{w}v}_{k-1, \delta}
 &\leq \norm{\sum_{i}f_{i}\nabla_{w}v_{i}}_{k-1,\delta}
 \\&\leq c\sum_{i}\norm{f_{i}}_{k-1,\delta}\norm{\nabla_{w}v_{i}}^{1}_{k-1,\delta}
 \\
 &\leq c'\norm{v}_{k-1, \delta}
 \\&\leq c'(\norm{D\dbar\phi}_{k-1,\delta}+2\norm{\pi_{V}D\dbar\phi}_{k,\delta}) \text{ using (\ref{hjlk4})}
 \\&\leq c''(\norm{\phi}^{1}_{k-1,\delta}+\norm{\pi_{V}D\dbar\phi}_{k,\delta})\end{split}\end{equation}
 
Using the inequality (\ref{hjlk3}) for the right-most term in (\ref{hjlk}) and rearranging gives 
\[\norm{\pi_{V}D\dbar \nabla_{w}\phi}_{k-1,\delta}\leq c(\norm{\pi_{V}D\dbar \phi}_{k, \delta}+\norm{\phi}^{1}_{k-1, \delta})\ .\]
Then, using our inductive hypothesis on the leftmost and rightmost terms above,
\[\norm{\nabla_{w}\phi}^{1}_{k-1,\delta}\leq c\norm{\pi_{V}D\dbar\phi}_{k, \delta}\ .\]
In the above inequality, $c$ depends on $w$, which is any extendable $\C\infty1$ vectorfield  lifted from a vectorfield on $\ex F$. As we can span all vectorfields on $\hat{\ex C}$ by real functions times a finite collection of such $w$, it follows that 
\begin{equation}\label{k bound}\norm{\phi}^{1}_{k,\delta}\leq c\norm{\pi_{V}D\dbar\phi}_{k, \delta}\end{equation}
where $c$ is independent of $\phi$. By induction on $k$, the bound \ref{k bound} holds for all $k$ (with different constants for different $k$).

We still need to prove that the bound (\ref{k bound}) implies that the restriction, to $Y/V_{k,\delta}$, of our low regularity inverse   is the inverse we want. Suppose that $\theta$ is in $Y^{\infty,\underline 1}$. Let us argue that $\lrb{\pi_{V}D\dbar}^{-1}(\pi_{V}\theta)$ is in $X^{\infty,\underline 1}$, by first arguing that all its derivatives exist. For each $f$, Lemma \ref{gluing} implies that $\lrb{\pi_{V}D\dbar(f)}^{-1}(\pi_{V}\theta(f))$ is $\C\infty1$. Using the map $G$ from Lemma \ref{gluing}, we can extend $\lrb{\pi_{V}D\dbar(f)}^{-1}(\pi_{V}\theta(f))$ to a $\C\infty1$ section $\phi_{0}$, so $D\dbar\phi_{0}$ plus some $\C\infty1$ section $v_{0}$ of $V$  agrees with $\theta$ over $f$. For showing that derivatives exist, we can assume without losing generality that $\ex F=\mathbb R^{n}$, with coordinates $x_{i}$, and $f$ at $0$. We can approximate $\theta$ to order $k$ over $0$ by 
\[D\dbar\lrb{\sum_{\abs\alpha\leq k} (x^{\alpha}\phi_{\alpha})}+\sum x^{\alpha}v_{\alpha}=\sum_{\abs \alpha\leq k}x^{k}(D\dbar\phi_{\alpha}+v_{\alpha})\] where, for each multi-index $\alpha$,   $\phi_{\alpha}$ is a section of $X^{\infty,\underline 1}$ and $v_{\alpha}$ is a $\C\infty1$ section of $V$ so that over $0$, 
\[D^{\alpha} \lrb{\theta -\sum_{\abs {\alpha'}<\abs\alpha}x^{k}(D\dbar\phi_{\alpha'}+v_{\alpha'})}=\abs\alpha!\lrb{D\dbar\phi_{\alpha}+v_{\alpha}}\ .\]
Above,  $D^{\alpha}$ indicates the partial derivative corresponding to the multi-index $\alpha$ --- this partial derivative is well defined independent of the lifts of $\frac\partial{\partial x_{i}}$ or connection,  because $D^{\alpha}$ is acting on something that vanishes to order $(\abs\alpha-1)$ over $0$.   
We can make $\norm{\pi_{V}(\theta-D\dbar\sum x^{\alpha}\phi_{\alpha})}_{k,\delta}$, and therefore $\norm{\lrb{\pi_{V}D\dbar}^{-1}(\pi_{V}\theta)-\sum x^{\alpha}\phi_{\alpha}}^{1}_{k,\delta}$, as small as we like by restricting to a sufficiently small neighorhood of $0$. It follows that at $0$, the derivatives of $\lrb{\pi_{V}D\dbar}^{-1}(\pi_{V}\theta)$ exist up to order $k$. Repeating the argument for different $k$ and at different points gives that all derivatives of $\lrb{\pi_{V}D\dbar}^{-1}(\pi_{V}\theta)$ exist. The bounds (\ref{k bound}) together with Lemma 7.8 of \cite{iec} then imply that $\lrb{\pi_{V}D\dbar}^{-1}(\pi_{V}\theta)$ is $\C\infty\delta$. Using Theorem \ref{invertible1} and Proposition \ref{invertible2},  we can locally repeat this argument for any $0<\delta<1$, to get that $\lrb{\pi_{V}D\dbar}^{-1}(\pi_{V}\theta)$ is in fact $\C\infty1$. 
 
 Now we know that $\lrb{\pi_{V}D\dbar}^{-1}$ sends $\C\infty 1$ sections to $\C\infty1$ sections, it follows that $\lrb{\pi_{V}D\dbar}^{-1}$ is our required inverse. For any $\theta\in Y/V_{k,\delta}$, we can choose a sequence of $\C\infty1$ sections $\theta'$ converging to $\theta$, then $\lrb{\pi_{V}D\dbar}^{-1}(\theta')$ is a sequence of $\C\infty1$ sections converging in $X_{k,\delta}$ to ${\lrb{\pi_{V}D\dbar}^{-1}\theta}$.

\stop

\begin{cor}\label{inverse regularity}Suppose that $V\subset Y^{\infty,\underline 1}$ is a finite-rank sub-vectorbundle, so that for every curve $f$ in our family, $D\dbar(f):X^{\infty,\underline 1}(f)\longrightarrow Y^{\infty,\underline1}(f)$ is injective and complementary to $V(f)$. Then, for each $\theta\in Y^{\infty,\underline 1}$, there exists a unique $\nu\in X^{\infty,\underline 1}$ and $\C\infty1$ section $v$ of $V$ so that 
\[D\dbar\nu=\theta+v\ .\]
\end{cor}

\pf

 Lemmas \ref{fredholm} and \ref{linear regularity} imply that $\pi_{V}\circ D\dbar(f):X_{\delta}\longrightarrow Y/V_{\delta}$ is invertible for all $f$, so there is a unique $\nu(f)\in X_{\delta}(f)$ and $v(f)\in V(f)$ so that $D\dbar\nu(f)=\theta(f)+v(f)$. Lemma \ref{linear regularity} implies that $\nu(f)$ is $\C\infty1$, and it remains to verify that the corresponding global section $\nu$ of $X$ is $\C\infty1$; this may be done locally.  Theorem \ref{invertible1}, and Propositions \ref{invertible2} and \ref{invertible3} imply that locally, $\pi_VD\dbar:X_{k,\delta}\longrightarrow Y/V_{k,\delta}$ has a bounded inverse for all $k\in \mathbb N$ and $\delta<1$.  Because $\nu(f)$ is the unique solution to $\pi_VD\dbar(f)\nu(f)=\pi_V(\theta(f))$, the inverse image of $\pi_V(\theta)$ must coincide with $\nu$, so $\nu$ is locally in $X_{k,\delta}$ for all $k\in\mathbb N$ and $\delta<1$. Therefore   $\nu$ is locally $\C\infty1$.

\stop

\
%
%
%

We can now prove Theorem \ref{linear theorem}, stated on page \pageref{linear theorem}. First, $D\dbar(f):X^{\infty,\underline 1}(f)\longrightarrow Y^{\infty,\underline 1}(f)$ has closed image, finite-dimensional kernel and cokernel, and the required index. This follows from lemmas \ref{fredholm}, \ref{linear regularity}, and \ref{index}. It remains to show that, if $D\dbar(f)$ is transverse to $V(f)\subset Y^{\infty,\underline 1}(f)$,  there exists a neighborhood  $\ex F'$ of $ f\subset \ex F$ so that $D\dbar$ surjects onto $Y^{\infty,\underline 1}(\ex F')/V$ and $D\dbar^{-1}(V)$ restricted to $X^{\infty,\underline1}(\ex F')$ is a finite-rank sub-vectorbundle in the sense of Definition \ref{sub-vectorbundle}.

Choose sections $\{s_{i}\}$ of $\hat{\ex C}\longrightarrow \ex F$ so that $D\dbar(f)^{-1}(V(f))=0$, when the domain of $D\dbar$ is restricted to the subspace of $X^{\infty,\underline 1}$ vanishing at the image of all $s_{i}$. On some neighborhood $\ex F'$ of $f\in \ex F$,  extend $V$ to a finite-rank sub-vectorbundle $V'\subset Y^{\infty,\underline 1}(\ex F')$ so that $V'(f)$ is complementary to  the image of (the restricted) $D\dbar(f)$. Then Proposition \ref{invertible2} implies that, by shrinking our neighborhood $\ex F'$, we may assume that the same holds for all $f'$ in $\ex F'$. 

 Choose a finite-rank sub-vectorbundle $K'\subset X^{\infty,\underline 1}(\ex F')$ that is a complement to the subspace of sections vanishing at $s_{i}$.  We can parametrize $K'$  by the value of its sections at $s_{i}$.  Corollary \ref{inverse regularity} allows us to parametrize $D\dbar^{-1}(V)$ by $K'$   as follows: for any section $\nu$ of $K'$, Corollary \ref{inverse regularity} implies that  there exists a unique section $\nu_{0}\in X^{\infty,\underline 1}(\ex F')$ that vanishes at the image of all $s_{i}$ so that $D\dbar\nu=\dbar\nu_{0}\mod V'$. We may modify $K'$ to $K''$ by replacing sections $\nu$ of $K'$ with sections  $\nu-\nu_{0}$ of $K''$. Such sections are still parametrized by their values at $s_{i}$, and the map $\nu\mapsto\nu-\nu_{0}$ is linear with respect to multiplication by $\C\infty1(\ex F',\mathbb R)$, so $K''$ is again a finite-rank sub-vectorbundle of $X^{\infty,\underline 1}(\hat f')$. The uniqueness part of Corollary \ref{inverse regularity} implies that $K''=D\dbar^{-1}V'$. 
 
 Now $D\dbar$ restricted to $K''$ is a $\C\infty1$ map of vectorbundles  $K''\longrightarrow V'$.  By assumption, this map is transverse to $V\subset V'$ at $f$, therefore by shrinking $\ex F'$, we may assume that it is transverse to $V$ everywhere in $\ex F'$, so $D\dbar^{-1}(V)$ is a finite-rank $\C\infty1$ sub-vectorbundle of  $K''$. Corollary \ref{inverse regularity} has already told us that $D\dbar$ surjects onto $Y^{\infty,\underline 1}(\ex F')/V'$, and we have now assumed that $D\dbar$ restricted to $K''$ surjects onto $V'/V$. Therefore, $D\dbar$ surjects onto $Y^{\infty,\underline1}(\ex F')/V$. This completes the proof of Theorem \ref{linear theorem}.

\subsection{Regularity for the $\dbar$ equation on variations of a family}

\

\begin{thm}\label{family regularity}
Suppose that $V\subset Y^{\infty,\underline 1}(\ex F)$ is a finite-rank sub-vectorbundle (Definition \ref{sub-vectorbundle}) and at  $f\in \ex F$,   $\dbar f\in V(f)$, and
\[ D\dbar(f):X^{\infty,\underline 1}(f)\longrightarrow Y^{\infty,\underline 1}(f)\]
is injective and has image complementary to $V(f)$. Then, there exists a neighborhood $\ex F '$ of $f\in \ex F$, and a neighborhood $O$ of $0$ in $X_{0,\delta}(\ex F')$ so that:
\begin{enumerate}
\item\label{fr1} Given any $f'\in \ex F'$,  and $\nu\in O$,
\[\pi_{V}D\dbar(\nu(f')) :X_{\delta'}(f')\longrightarrow Y/V_{\delta'}(f')\]
is invertible.
\item \label{fr3}
The map \[\pi_{V}\dbar:O\longrightarrow Y/V_{0,\delta}\] is a homeomorphism onto a neighborhood of $0\in Y/V_{0,\delta}$.

\item\label{fr2} Given any $f'\in \ex F'$, let $O(f')$ be the restriction of $O$ to $f'$. Then  the map \[\pi_{V}\dbar(f'):O(f')\longrightarrow Y/V_{0,\delta}(f')\] is also a homeomorphism onto a neighborhood of $0\in Y/V_{0,\delta}(f')$.

 \item\label{fr4} For any $\nu\in O$,  $\nu$ has regularity $\C\infty1$ if $\pi_{V}\dbar\nu$ has regularity $\C\infty1$.  In particular, there is a unique $\C\infty1$ solution to the equation $\pi_{V}\dbar\nu=0 $ over $\ex F'$.

\end{enumerate}
  
\end{thm}

\

\pf

 Apply Theorem \ref{invertible1}  and Proposition \ref{invertible2} to see that there exists an allowable collection of coordinate charts covering an open neighborhood of $f$ so that, in these coordinates, there is a bounded inverse to $\pi_{V}\circ D\dbar:X_{0,\delta}\longrightarrow Y/V_{0,\delta}$.

Theorem \ref{smooth dbar} implies that $\pi_{V}\circ \dbar:X_{0,\delta}\longrightarrow Y/V_{0,\delta}$ is $C^{1}$, so we may choose a neighborhood $O$ of $0\in X_{0,\delta}$ so that for $\nu\in O$,
\begin{equation}\label{Ddbar estimate}\norm{\lrb{\pi_{V}\circ D\dbar}^{-1}\circ\pi_{V}\circ D\dbar(\nu)(\cdot)-\id}<\frac 12\ .
\end{equation}

 Therefore, $\pi_{V}\circ \dbar$ is a homeomorphism from $O$ to an open subset of $Y/V_{0,\delta}$.  As $\pi_V\dbar$ of the zero section is a $\C\infty1$ section which vanishes at $f$,   $\pi_{V}\circ \dbar(O)$ contains some section that vanishes on a neighborhood of $f$, so,  by restricting to a small enough  neighborhood of $f$,  the image under $\pi_{V}\dbar$ of  $O$ will contain $0 \in Y/V_{0,\delta}$. 
We have therefore proved item \ref{fr3} of our theorem. As the construction may be restricted to any curve $f'$, item \ref{fr2} also holds, and restricting our estimate on the inverse of $\pi_{V}D\dbar$ to $\nu(f')$ gives item \ref{fr1}.

%
%

\

It remains to prove item \ref{fr4}. We must show  that, if $\nu\in O$  and $\pi_{V}\dbar \nu=\theta$ is $\C\infty1$,  then $\nu$ is $\C\infty1$. Do this locally around any curve  $g$ in $\ex F'$ as follows: First, as $\pi_{V}\dbar \nu(g)$ is $\C\infty1$,  $\dbar\nu(g)\in \C\infty1$ and Proposition \ref{curve regularity} implies that $\nu(g)$ is $\C\infty1$. Choose any $\C\infty1$ extension  of $\nu(g)$ to an open neighborhood of $g$, (such an extension can be constructed using the map $G$ from Lemma \ref{gluing}), and call this extension $\nu_{0}$.

Consider the following Newton iteration scheme:
\[\nu_{n+1}:=\nu_{n} -\lrb{\pi_{V}D\dbar(\nu_{n})}^{-1}\lrb{\pi_{V}\dbar \nu_{n}-\theta}\]
On a sufficiently small neighborhood $U$ of $g$, the estimate (\ref{Ddbar estimate}) implies that  $(\pi_{V}D\dbar(\nu_{0}))^{-1}$  exists and is a bounded operator. Moreover, by choosing $U$ small, we may ensure that $\norm{\pi_{V}\dbar\nu_{0}-\theta}_{0,\delta}$ is as small as we like. Therefore, as $D\dbar:X_{0,\delta}\longrightarrow Y_{0,\delta}$ is $C^{1}$, for a sufficiently small $U$,  the above Newton iteration scheme will converge in $X_{0,\delta}(U)$ and remain in (the restriction to $U$ of) our chosen neighborhood $O$.

 We shall now prove that, if $\nu_{n}$ is $\C\infty1$, then $\nu_{n+1}$ is $\C\infty1$ as well. As $\pi_{V}D\dbar(\nu_{n}):X_{0,\delta}\longrightarrow Y/V_{0,\delta}$ is invertible, it follows that for all $f'$ in $\ex F'$ , $D\dbar(\nu_{n}(f')):X_\delta(f')\longrightarrow Y/V_{\delta}(f')$ is invertible. It follows that $D\dbar(\nu_{n})(f'):X^{\infty,\underline 1}(f')\longrightarrow Y^{\infty,\underline1}(f') $ is injective and complementary to $V(f)$. Corollary \ref{inverse regularity} then implies that  $\nu_{n+1}$ is $\C\infty1$, so by induction $\nu_{n}$ is $\C\infty1$ for all $n$.

 Now we shall see that this Newton iteration scheme will converge in $X_{k,\delta}$ for all  $k$ and all $\delta<1$ when restricted to a small enough  open neighborhood of $g$ (dependent on $k$ and $\delta$).

\

Claim: The size of   $\norm{\pi_{V}\dbar\nu_{n}-\theta}_{n,\delta}$ can be made arbitrarily small by restricting to a suitably small  open neighborhood of  the curve over $g$. More precisely, $\norm{\pi_{V}\dbar\nu_{n}-\theta}_{n,\delta}$ can be made arbitrarily small by restricting to a small enough neighborhood of the curve over $g$ while using the same metric and coordinate chart choices in the definition of $\norm{\cdot}_{n,\delta}$. 

\

Consider $\dbar\nu_{n}$ as a section of the bundle $Y$.  As $\nu_{n}$ and $\theta$ are $\C\infty1$, the above claim is equivalent to $\dbar\nu_{n}$ being tangent to order $n$ over the domain of  $g$ to  some $\C\infty1$ section $\theta'$ with  $\pi_{V}\theta'=\theta$. 

We shall prove the above claim by induction --- it holds for $n=0$, now assume that it holds for some $n$. Therefore, $\norm{\nu_{n+1}-\nu_{n}}^{1}_{n,\delta}$ can be made arbitrarily small by restricting to a suitably small open neighborhood of the domain of $g$, i.e. $\nu_{n}$ and $\nu_{n+1}$ are tangent to order $n$ on the domain of $g$. Define an operator $(D^{n+1})$ as follows: The domain of $(D^{n+1})$ consists of $\C\infty1$ sections of $Y$ tangent to $\dbar\nu_{n}$ to order $n$ over $g$. Define $(D^{n+1})\theta$ to be the derivative to order $n+1$ of the section $\theta$, restricted to  the domain of $g$. With this domain, $(D^{n+1})$ is an affine operator --- so $(D^{n+1})(\lambda \theta_1+(1-\lambda)\theta_2)=\lambda (D^{n+1})(\theta_1)+(1-\lambda)(D^{n+1})(\theta_2)$ for any constant $\lambda$, and $\theta_i$ in the domain of $(D^{n+1})$. The operator   $(D^{n+1})\dbar$ is also an affine operator when restricted to $\nu\in \C\infty1$ so that $\nu$ is tangent to $\nu_{n}$ to order $n$ over the domain of $g$. 
Therefore, for such $\nu$, using that $(D^{n+1})\dbar$ is affine gives  
\[(D^{n+1})\dbar(\nu_n+\lambda(\nu-\nu_n))=(D^{n+1})\dbar \nu_n+\lambda((D^{n+1})\dbar\nu-(D^{n+1})\dbar \nu_n)
\ .\]
Then, as $(D^{n+1})$ is affine, it commutes with taking first-order approximations, so taking $(D^{n+1})$ of the approximation of $\dbar(\nu_n+\lambda(\nu-\nu_n))$ at $\lambda=0$ gives 
\[(D^{n+1})(\dbar \nu_n+\lambda D\dbar(\nu_{n})(\nu-\nu_{n})))=(D^{n+1})\dbar \nu_n+\lambda((D^{n+1})\dbar\nu-(D^{n+1})\dbar \nu_n)\ .\]
Setting $\lambda=1$ in the above expression gives \[(D^{n+1})\dbar\nu=(D^{n+1})(\dbar\nu_{n} +D\dbar(\nu_{n})(\nu-\nu_{n}))\] so
\[\begin{split}(D^{n+1})\dbar\nu_{n+1}&=(D^{n+1})(\dbar\nu_{n}+D\dbar(\nu_{n})(\nu_{n+1}-\nu_{n}))
\\&=(D^{n+1})\lrb{\dbar\nu_{n}-D\dbar(\nu_{n})\lrb{\pi_{V}D\dbar(\nu_{n})}^{-1}\lrb{\pi_{V}\dbar \nu_{n}-\theta}}
\\&=(D^{n+1})(\theta')\end{split}\]
where $\theta'$ is a $\C\infty1$ section of $Y$ so that
\[\pi_{V}\theta'=\pi_{V}\dbar\nu_{n}-\pi_{V}\dbar\nu_{n}+\theta=\theta\ .\]
In other words, $\dbar\nu_{n+1}$ will be tangent to  order $n+1$ over the domain of $g$ to some $\theta'$ so that $\pi_{V}\theta'=\theta$. By induction, this is true for all $n$, and our claim has been proven.

\

Apply Theorem \ref{invertible1}, Proposition \ref{invertible2}, Proposition \ref{invertible3} and Theorem \ref{smooth dbar} and the above claim to infer that for any  $\delta<1$, restricted to a small enough neighborhood of the curve over $g$, the above Newton iteration scheme converges in $X_{k,\delta}$ to our solution $\nu$ to $\pi_{V}\dbar \nu=\theta$. This implies that, for all  $k$ and $\delta<1$, there exists some  open neighborhood of $g$ so that our solution $\nu$ is in $C^{k,\delta}$ restricted to this neighborhood. Repeating the argument around any point gives that $\nu$ is $\C\infty1$.

\stop

\bibliographystyle{plain}
\bibliography{ref.bib}

\begin{thebibliography}{10}

\bibitem{acgw}
Dan Abramovich and Qile Chen.
\newblock Stable logarithmic maps to {D}eligne-{F}altings pairs {II}.
\newblock {\em The Asian Journal of Mathematics}, 18(3):465--488, 2014.

\bibitem{CLW}
Bohui Chen, An-Min Li, and Bai-Ling Wang.
\newblock Virtual neighborhood technique for pseudo-holomorphic spheres.
\newblock arXiv:1306.3276.

\bibitem{Chen}
Qile Chen.
\newblock Stable logarithmic maps to {D}eligne-{F}altings pairs {I}.
\newblock {\em Ann. of Math. (2)}, 180(2):455--521, 2014.

\bibitem{FO}
Kenji Fukaya and Kaoru Ono.
\newblock Arnold conjecture and {G}romov-{W}itten invariant.
\newblock {\em Topology}, 38(5):933--1048, 1999.

\bibitem{GSlogGW}
Mark Gross and Bernd Siebert.
\newblock Logarithmic {G}romov-{W}itten invariants.
\newblock {\em J. Amer. Math. Soc.}, 26(2):451--510, 2013.

\bibitem{polyfold0}
H.~Hofer.
\newblock A general {F}redholm theory and applications.
\newblock In {\em Current developments in mathematics, 2004}, pages 1--71. Int.
  Press, Somerville, MA, 2006.

\bibitem{polyfold1}
H.~Hofer, K.~Wysocki, and E.~Zehnder.
\newblock A general {F}redholm theory. {I}. {A} splicing-based differential
  geometry.
\newblock {\em J. Eur. Math. Soc. (JEMS)}, 9(4):841--876, 2007.

\bibitem{polyfoldint}
H.~Hofer, K.~Wysocki, and E.~Zehnder.
\newblock Integration theory on the zero sets of polyfold {F}redholm sections.
\newblock {\em Math. Ann.}, 346(1):139--198, 2010.

\bibitem{polyfoldgw}
Helmut Hofer, Kris Wysocki, and Eduard Zehnder.
\newblock A general {F}redholm theory. {III}. {F}redholm functors and
  polyfolds.
\newblock {\em Geom. Topol.}, 13(4):2279--2387, 2009.

\bibitem{polyfoldsc}
Helmut Hofer, Kris Wysocki, and Eduard Zehnder.
\newblock sc-smoothness, retractions and new models for smooth spaces.
\newblock {\em Discrete Contin. Dyn. Syst.}, 28(2):665--788, 2010.

\bibitem{polyfold2}
Helmut Hofer, Krzysztof Wysocki, and Eduard Zehnder.
\newblock A general {F}redholm theory. {II}. {I}mplicit function theorems.
\newblock {\em Geom. Funct. Anal.}, 19(1):206--293, 2009.

\bibitem{IonelGW}
Eleny-Nicoleta Ionel.
\newblock G{W} invariants relative to normal crossing divisors.
\newblock {\em Adv. Math.}, 281:40--141, 2015.

\bibitem{Tian-Li}
Jun Li and Gang Tian.
\newblock Virtual moduli cycles and {G}romov-{W}itten invariants of general
  symplectic manifolds.
\newblock In {\em Topics in symplectic {$4$}-manifolds ({I}rvine, {CA}, 1996)},
  First Int. Press Lect. Ser., I, pages 47--83. Int. Press, Cambridge, MA,
  1998.

\bibitem{Liu-Tian}
Gang Liu and Gang Tian.
\newblock Floer homology and {A}rnold conjecture.
\newblock {\em J. Differential Geom.}, 49(1):1--74, 1998.

\bibitem{mcowen}
Robert~B. Lockhart and Robert~C. McOwen.
\newblock Elliptic differential operators on noncompact manifolds.
\newblock {\em Ann. Scuola Norm. Sup. Pisa Cl. Sci. (4)}, 12(3):409--447, 1985.

\bibitem{McDuff}
Dusa McDuff.
\newblock The virtual moduli cycle.
\newblock In {\em Northern {C}alifornia {S}ymplectic {G}eometry {S}eminar},
  volume 196 of {\em Amer. Math. Soc. Transl. Ser. 2}, pages 73--102. Amer.
  Math. Soc., Providence, RI, 1999.

\bibitem{MS}
Dusa McDuff and Dietmar Salamon.
\newblock {\em $J$-holomorphic Curves and Symplectic Topology}, volume~52 of
  {\em Colloquium Publications}.
\newblock American Mathematical Society, 2004.

\bibitem{tpgf}
Brett Parker.
\newblock Gluing formula for {G}romov-{W}itten invariants in a triple product.
\newblock \href{http://arxiv.org/abs/1511.0779}{arXiv:1511.0779}.

\bibitem{scgp}
Brett Parker.
\newblock Notes on exploded manifolds and a tropical gluing formula for
  {G}romov-{W}itten invariants.
\newblock \href{http://arxiv.org/abs/1605.00577}{arXiv:1605.00577}.

\bibitem{tec}
Brett Parker.
\newblock Tropical enumeration of curves in blowups of the projective plane.
\newblock \href{http://arxiv.org/abs/1411.5722}{arXiv:1411.5722}.

\bibitem{dre}
Brett Parker.
\newblock De {R}ham theory of exploded manifolds.
\newblock \href{http://arxiv.org/abs/1003.1977}{arXiv:1003.1977}, 2011.

\bibitem{iec}
Brett Parker.
\newblock Exploded manifolds.
\newblock {\em Adv. Math.}, 229:3256--3319, 2012.
\newblock \href{http://arxiv.org/abs/0910.4201}{arXiv:0910.4201}.

\bibitem{elc}
Brett Parker.
\newblock Log geometry and exploded manifolds.
\newblock {\em Abh. Math. Sem. Hamburg}, 82:43--81, 2012.
\newblock \href{http://arxiv.org/abs/1108.3713}{arxiv:1108.3713}.

\bibitem{evc}
Brett Parker.
\newblock Holomorphic curves in exploded manifolds: Kuranishi structure.
\newblock \href{http://arxiv.org/abs/1301.4748}{arXiv:1301.4748}, 2013.

\bibitem{uts}
Brett Parker.
\newblock Universal tropical structures for curves in exploded manifolds.
\newblock \href{http://arxiv.org/abs/1301.4745}{arXiv:1301.4745}, 2013.

\bibitem{tropicalIonel}
Brett Parker.
\newblock On the value of thinking tropically to understand {I}onel's {GW}
  invariants relative normal crossing divisors.
\newblock \href{http://arxiv.org/abs/1407.3020}{arXiv:1407.3020}, 2014.

\bibitem{cem}
Brett Parker.
\newblock Holomorphic curves in exploded manifolds: compactness.
\newblock {\em Adv. Math.}, 283:377--457, 2015.
\newblock \href{http://arxiv.org/abs/0911.2241}{arXiv:0911.2241}.

\bibitem{vfc}
Brett Parker.
\newblock Holomorphic curves in exploded manifolds: virtual fundamental class.
\newblock \href{http://arxiv.org/abs/1512.05823}{arXiv:1512.05823}, 2015.

\bibitem{3d}
Brett Parker.
\newblock Three dimensional tropical correspondence formula.
\newblock {\em Communications in Mathematical Physics}, 352(2):791--819, 2017.
\newblock \href{https://arxiv.org/abs/1608.02306}{arXiv:1608.02306}.

\bibitem{gfgw}
Brett Parker.
\newblock Tropical gluing formulae for {G}romov-{W}itten invariants.
\newblock \href{http://arxiv.org/abs/1703.05433}{arXiv:1703.05433}, 2017.

\bibitem{Ruanvirtual}
Yongbin Ruan.
\newblock Virtual neighborhoods and pseudo-holomorphic curves.
\newblock In {\em Proceedings of 6th {G}\"okova {G}eometry-{T}opology
  {C}onference}, volume~23, pages 161--231, 1999.

\end{thebibliography}
\end{document}